\input amstex
\documentstyle{amsppt}
%
\catcode`@=11
\redefine\output@{%
  \def\break{\penalty-\@M}\let\par\endgraf
  \ifodd\pageno\global\hoffset=105pt\else\global\hoffset=8pt\fi  
  \shipout\vbox{%
    \ifplain@
      \let\makeheadline\relax \let\makefootline\relax
    \else
      \iffirstpage@ \global\firstpage@false
        \let\rightheadline\frheadline
        \let\leftheadline\flheadline
      \else
        \ifrunheads@ 
        \else \let\makeheadline\relax
        \fi
      \fi
    \fi
    \makeheadline \pagebody \makefootline}%
  \advancepageno \ifnum\outputpenalty>-\@MM\else\dosupereject\fi
}
\catcode`\@=\active
\nopagenumbers
\def\negskp{\hskip -2pt}
\def\Img{\operatorname{Im}}
\def\Ker{\operatorname{Ker}}
\def\const{\operatorname{const}}
\def\MatGrSL{\operatorname{SL}}
\def\MatGrSO{\operatorname{SO}}
\def\MatGrO{\operatorname{O}}
\def\HSym{\operatorname{HSym}}
\def\Alpha{\operatorname{A}}
\def\compos{\,\raise 1pt\hbox{$\sssize\circ$} \,}
\def\idop{\operatorname{\bold{id}}}
\def\msum#1{\operatornamewithlimits{\sum^#1\!{\ssize\ldots}\!\sum^#1}}
\def\msums#1#2{\operatornamewithlimits{\sum^#1\!{\ssize\ldots}\!\sum^#2}}
\accentedsymbol\tx{\tilde x}
\accentedsymbol\tbZ{\kern3pt\tilde{\kern -3pt\bar Z}\vphantom{Z}}
\accentedsymbol\tbGamma{\kern0pt\tilde{\kern -0pt\bar\Gamma}%
\vphantom{\Gamma}}
\accentedsymbol\vnabla{\nabla\kern -7pt\raise 5pt\vbox{\hrule width 5.5pt}
\kern 1.7pt}
\accentedsymbol\bPsi{\kern 1pt\overline{\kern -1pt\boldsymbol\Psi
\kern -1pt}\kern 1pt}
\accentedsymbol\bbd{\kern 2pt\bar{\kern -2pt\bold d}}
\accentedsymbol\bd{\kern 2pt\bar{\kern -2pt d}}

\def\blue#1{#1}
\catcode`#=11\def\diez{#}\catcode`#=6
\catcode`_=11\def\podcherkivanie{_}\catcode`_=8
\catcode`~=11\catcode`~=\active
\def\mycite#1{\cite{\blue{#1}}\immediate\special{ps:
     ShrHPSdict begin /ShrBORDERthickness 0 def}}
\def\myciterange#1#2{\cite{\blue{#2}}\immediate\special{ps:
     ShrHPSdict begin /ShrBORDERthickness 0 def}}
\def\mytag#1{%
    \tag#1}
\def\mythetag#1{\thetag{\blue{#1}}\immediate\special{ps:
     ShrHPSdict begin /ShrBORDERthickness 0 def}}
\def\myrefno#1{\no#1}
\def\myhref#1#2{\blue{#2}\immediate\special{ps:
     ShrHPSdict begin /ShrBORDERthickness 0 def}}
\def\myEarXivlink{\myhref{http://arXiv.org}{http:/\negskp/arXiv.org}}
\def\myGeoCities{\myhref{http://www.geocities.com}{GeoCities}}
\def\mytheorem#1{\csname proclaim\endcsname{Theorem #1}}
\def\mythetheorem#1{\blue{#1}\immediate\special{ps:
     ShrHPSdict begin /ShrBORDERthickness 0 def}}
\def\mylemma#1{\csname proclaim\endcsname{Lemma #1}}
\def\mythelemma#1{\blue{#1}\immediate\special{ps:
     ShrHPSdict begin /ShrBORDERthickness 0 def}}
\def\mydefinition#1{\definition{Definition #1}}
\def\mythedefinition#1{\blue{#1}\immediate\special{ps:
     ShrHPSdict begin /ShrBORDERthickness 0 def}}

\font\tencyr=wncyr10
\font\sevencyr=wncyr7
\font\sevencyi=wncyi7
\pagewidth{360pt}
\pageheight{606pt}
\hphantom{a}\vskip 0pt
\line{\kern 10pt\includegraphics{spfun01.eps}\hss
\vtop{\hsize=270pt\baselineskip=2pt plus 1pt minus 1pt 
\lineskip=2pt plus 1pt minus 1pt 
\lineskiplimit=-0.5pt\sevencyi
Odnazhdy starik-egiptyanin pri\-zval k sebe svoego
starshego syna i govorit emu:\vphantom{!}\newline --- Synok, ya otkroyu 
tebe strashnuyu ta\ae nu. Zhre\-cy iz hrama solnca Amon-Ra --- ne
bogi\,!\newline 
--- Otec, kak mozhno\,! Oni zhe govoryat so zvezdami i predskazyvayut 
zatmeniya\,! --- izumilsya syn. No starik prodolzhil nevozmutimo:\newline 
--- Eshch\oe{} v det\vphantom{a}stve ya ta\ae kom pronik v ih hram i 
videl, chto oni edyat i spravlyayut nuzhdu sovsem tak zhe, kak my s
tobo\ae.\newline\vphantom{!}
\newline
\vphantom{a}\kern 10pt\sevencyr (Naveyano vremenem.
Egipet, drevnee carstvo.)}}
\vskip -120pt
\topmatter
\title
Spinor functions of spinors\\
and the concept of extended spinor fields.
\vskip 80pt
\endtitle
\author
R.~A.~Sharipov
\endauthor
\address Rabochaya street 5, 450003 Ufa, Russia
\endaddress
\email \vtop to 30pt{\hsize=280pt\noindent
\myhref{mailto:R\podcherkivanie Sharipov\@ic.bashedu.ru}
{R\_\hskip 1pt Sharipov\@ic.bashedu.ru}\newline
\myhref{mailto:r-sharipov\@mail.ru}
{r-sharipov\@mail.ru}\newline
\myhref{mailto:ra\podcherkivanie sharipov\@lycos.com}{ra\_\hskip 1pt
sharipov\@lycos.com}\vss}
\endemail
\urladdr
\vtop to 20pt{\hsize=280pt\noindent
\myhref{http://www.geocities.com/r-sharipov}
{http:/\negskp/www.geocities.com/r-sharipov}\newline
\myhref{http://www.freetextbooks.boom.ru/index.html}
{http:/\negskp/www.freetextbooks.boom.ru/index.html}\vss}
\endurladdr
\abstract
    Spinor fields depending on tensor fields and other spinor fields are
considered. The concept of extended spinor fields is introduced and the theory of differentiation for such fields is developed.
\endabstract
\endtopmatter
\loadbold
\loadeufb
\TagsOnRight
\document
\vskip -20pt

\rightheadtext{Spinor functions of spinors \dots}
\head
1. Introduction.
\endhead
    The {\it space-time} is a stage for all events in general relativity.
Saying the space-time, we understand a smooth $4$-dimensional manifold $M$
equipped with a pseudo-Euclidean metric $\bold g$ of the Minkowski-type 
signature $(+,-,-,-)$. Word lines of particles, tangent spaces of $M$ and
their light cones, local charts and local coordinates, tensors and tensor
fields in $M$, the metric connection $\Gamma$ with the components 
$$
\Gamma^k_{ij}=\sum^3_{s=0}\frac{g^{ks}}{2}\left(\frac{\partial g_{sj}}
{\partial x^i}+\frac{\partial g_{is}}{\partial x^j}
-\frac{\partial g_{ij}}{\partial x^s}\right)\!,
$$
the covariant differentiation $\nabla$ determined by $\Gamma$, and some 
other basic concepts of general relativity are assumed to be known to
the reader. The free textbooks \myciterange{1}{1--\,4}\linebreak are 
recommended as an introductory material for getting acquainted with that
basic concepts. The concepts of spinors and spinor fields are less 
commonly known. Therefore, we consider them in full details in the next 
eight sections~2 through 8. These sections form the preliminary part of 
the present paper.\par
     The concept of an extended tensor field in its full generality was
introduced in paper \mycite{5}. It arises if one consider a tensor-valued
function 
$$
\hskip -2em
\bold X=\bold X(p,\bold T[1],\ldots,\bold T[Q])
\mytag{1.1}
$$
with one point argument $p\in M$ and several tensorial arguments 
$\bold T[1],\,\ldots,\,\bold T[Q]$. The paper \mycite{6} illustrates 
how extended tensor fields are applied in the theory of continuous
media. When passing to the quantum area, one should also take into 
account the specifically quantum phenomenon of {\it spin} which does 
not reduce to the pure rotation. Therefore, we need to add some 
spin-tensors $\bold S[1],\,\ldots,\,\bold S[J]$ to the arguments of 
the tensor-valued function $\bold X$ in \mythetag{1.1}
$$
\pagebreak
\hskip -2em
\bold X=\bold X(p,\bold S[1],\ldots,\bold S[J],\bold T[1],\ldots,
\bold T[Q])
\mytag{1.2}
$$
and extend the theory to the spin-tensorial functions with the same
kind arguments:
$$
\hskip -2em
\bold U=\bold U(p,\bold S[1],\ldots,\bold S[J],\bold T[1],\ldots,
\bold T[Q]).
\mytag{1.3}
$$
This is the main purpose of the present paper. We reproduce the results 
of \mycite{5} as applied to the extended tensor fields \mythetag{1.2} and 
extended spin-tensorial fields \mythetag{1.3}.
\head
2. Spinors and spin-tensors.
\endhead
    Let $M$ be the space-time. Denote by $SM$ some smooth $2$-dimensional 
complex vector bundle over $M$. We shall call it the {\it spinor} bundle. 
Its smooth global and local sections are called {\it spinor fields}. Apart
from $SM$, we consider the tangent bundle $TM$. For the sake of uniformity 
we denote the canonical projections for both bundles $TM$ and $SM$ with the
same symbol $\pi$:
$$
\xalignat 2
&\hskip -2em
\pi\!:\,TM\to M,
&&\pi\!:\,SM\to M.
\mytag{2.1}
\endxalignat
$$
Let $p\in M$ be a point of the space-time $M$. We denote by $T_p(M)$ and 
$S_p(M)$ the {\it fibers} of the bundles $TM$ and $SM$ over the point $p$.
In other words, $T_p(M)$ and $S_p(M)$ are the total preimages of the point
$p$ under the mappings \mythetag{2.1}. Note that $T_p(M)$ is a 
$4$-dimensional vector space over the real numbers $\Bbb R$, while $S_p(M)$
is a $2$-dimensional vector space over the complex numbers $\Bbb C$.\par
    Let's denote by $T^*_p(M)$ and $S^*_p(M)$ the conjugate spaces for
$T^*_p(M)$ and $S^*_p(M)$ respectively. The elements of $T_p(M)$ are called 
{\it vectors}, the elements of $T^*_p(M)$ are called {\it covectors}. 
Similarly, the elements of $S_p(M)$ are called {\it spinors}, the elements 
of $S^*_p(M)$ are called {\it cospinors}. Let's introduce the following spaces:
$$
\align
&\hskip -2em
T^r_s(p,M)=\overbrace{T_p(M)\otimes\ldots\otimes T_p(M)}^{\text{$r$
times}}\otimes\underbrace{T^*_p(M)\otimes\ldots\otimes 
T^*_p(M)}_{\text{$s$ times}},
\mytag{2.2}\\
&\hskip -2em
S^\alpha_\beta(p,M)=\overbrace{S_p(M)\otimes\ldots\otimes 
S_p(M)}^{\text{$\alpha$ times}}\otimes
\underbrace{S^*_p(M)\otimes\ldots\otimes 
S^*_p(M)}_{\text{$\beta$ times}}.
\mytag{2.3}
\endalign
$$
The elements of $T^r_s(p,M)$ are called {\it tensors of the type 
$(r,s)$}. Similarly, the elements of $S^\alpha_\beta(p,M)$ are called 
{\it spin-tensors of the type $(r,s)$}. The space \mythetag{2.2} is
a linear space over the real numbers, while \mythetag{2.3} is a complex
linear space. Apart from \mythetag{2.2} and \mythetag{2.3}, one can
consider their tensor product
$$
\hskip -2em
S^\alpha_\beta T^r_s(p,M)=S^\alpha_\beta(p,M)\otimes T^r_s(p,M).
\mytag{2.4}
$$
The elements of the tensor product \mythetag{2.4} are also called 
{\it spin-tensors}. Note that $S^\alpha_\beta T^r_s(p,M)$ is
a complex linear space. Remember also that for any real linear space 
$V$ the complex linear space $\Bbb C\otimes V$ is defined, it is called
the {\it complexification} of $V$. Let's denote by $\Bbb CT_p(M)$ the
complexification of the tangent space:
$$
\pagebreak
\hskip -2em
\Bbb CT_p(M)=\Bbb C\otimes T_p(M).
\mytag{2.5}
$$
Then let's use \mythetag{2.5} in the following tensor product:
$$
\hskip -2em
\Bbb CT^r_s(p,M)=\overbrace{\Bbb CT_p(M)\otimes\ldots\otimes
\Bbb CT_p(M)}^{\text{$r$times}}\otimes\underbrace{\Bbb CT^*_p(M)
\otimes\ldots\otimes \Bbb CT^*_p(M)}_{\text{$s$ times}}.
\mytag{2.6}
$$
The tensor product \mythetag{2.4} now can be written as
$$
\hskip -2em
S^\alpha_\beta T^r_s(p,M)=S^\alpha_\beta(p,M)\otimes\Bbb CT^r_s(p,M).
\mytag{2.7}
$$
This representation \mythetag{2.7} of $S^\alpha_\beta T^r_s(p,M)$ is
absolutely equivalent to \mythetag{2.4}.\par
For the conjugate spaces $\Bbb CT^*_p(M)$ and $T^*_p(M)$ in 
\mythetag{2.6} and \mythetag{2.2} and for the tensor products 
\mythetag{2.6} and \mythetag{2.2} themselves we have the following
obvious equalities:
$$
\xalignat 2
&\hskip -2em
\Bbb CT^*_p(M)=\Bbb C\otimes T^*_p(M),
&\Bbb CT^r_s(p,M)=\Bbb C\otimes T^r_s(p,M).
\mytag{2.8}
\endxalignat
$$
Due to \mythetag{2.8} we have the semilinear isomorphism of complex
conjugation:
$$
\hskip -2em
\tau\!:\,\Bbb CT^r_s(p,M)\to \Bbb CT^r_s(p,M).
\mytag{2.9}
$$
It is introduced by the formula $\tau(\alpha\otimes\bold X)=
\overline{\alpha}\otimes\bold X$, where $\bold X\in T^r_s(p,M)$,
$\alpha\in\Bbb C$, and $\overline{\alpha}$ is the conjugate number
for $\alpha$.\par
     The complexified tangent spaces \mythetag{2.5} are glued into
a complex vector bundle $\Bbb CT\!M$. It is called the {\it complexified
tangent bundle}. Its local and global sections are called {\it complexified
vector fields}. However, below we shall often call them vector fields
for the sake of brevity.
\head
3. Semilinear functions and Hermitian conjugate spaces.
\endhead
    Let $V$ be some finite-dimensional linear space over the field of
complex numbers. Its dual space $V^*$ by definition is the set of linear 
functions $f=f(\bold v)$, where $\bold v\in V$. This means 
that each $f\in V^*$ satisfies the following two linearity conditions:
\roster
\item\quad $f(\bold v_1+\bold v_2)=f(\bold v_1)+f(\bold v_2)$
\ for any two vectors \ $\bold v_1,\bold v_2\in V$;
\item\quad $f(\alpha\,\bold v)=\alpha\,f(\bold v)$
\ for any \ $\bold v\in V$ and for any $\alpha\in\Bbb C$.
\endroster
In geometric terminology linear functions $f\in V^*$ are called
{\it covectors}, their values, when applied to a vector $\bold v\in V$,
are written in the form of a scalar product:
$$
\hskip -2em
f(\bold v)=(f,\,\bold v).
\mytag{3.1}
$$
Semilinear functions in $V$ are defined by the following two conditions:
\roster
\item\quad $f(\bold v_1+\bold v_2)=f(\bold v_1)+f(\bold v_2)$
\ for any two vectors \ $\bold v_1,\bold v_2\in V$;
\item\quad $f(\alpha\,\bold v)=\bar\alpha\,f(\bold v)$
\ for any \ $\bold v\in V$ and for any $\alpha\in\Bbb C$.
\endroster
By means of $\bar\alpha$ in \therosteritem{2} above we denote the
conjugate complex number $\bar\alpha=x-iy$ for a complex number
$\alpha=x+iy$. Semilinear functions form another linear vector space
$V^{\sssize\dagger}$, it is called the {\it Hermitian conjugate space\/}
or {\it Hermitian dual space\/} for $V$. Like in \mythetag{3.1}, the
value of a function $f\in V^{\sssize\dagger}$ is denoted as a scalar
product:
$$
\hskip -2em
f(\bold v)=\langle\bold v,\,f\rangle.
\mytag{3.2}
$$
However, \mythetag{3.2} is a Hermitian scalar product. For this reason
we used the angular brackets for it. The elements of $V^{\sssize\dagger}$
are called {\it conjugate covectors}.\par
     Each linear function $f\in V^*$ can be converted to a semilinear 
function $\bar f\in V^{\sssize\dagger}$ by passing to conjugate numbers
in its values:
$$
\hskip -2em
\bar f(\bold v)=\overline{f(\bold v)}
\mytag{3.3}
$$
Conversely, if $f\in V^{\sssize\dagger}$, then $\bar f\in V^*$. Thus, 
the formula \mythetag{3.3} defines the canonical semilinear isomorphism:
$V^*\rightleftarrows V^{\sssize\dagger}$. We shall use the same symbol
$\tau$ for both semilinear mappings implementing this isomorphism:
$$
\xalignat 2
&\hskip -2em
\tau\!:\,V^*\to V^{\sssize\dagger},
&&\tau\!:\,V^{\sssize\dagger}\to V^*.
\mytag{3.4}
\endxalignat
$$
Let's denote by $V^{*\sssize\dagger}$ the conjugate dual for the dual 
space $V^*$. It is clear that $V^{{\sssize\dagger}*}$ is canonically 
isomorphic to $V^{*\sssize\dagger}$, so that one can treat them as the 
same space. Indeed, for any vector $\bold v\in V$ the following function 
is defined:
$$
\hskip -2em
w_{\bold v}(f)=\langle\bold v,\,f\rangle\text{, \ where \ }
f\in V^{\sssize\dagger}.
\mytag{3.5}
$$
The formula \mythetag{3.5} defines the canonical semilinear isomorphism
$V^{{\sssize\dagger}*}\rightleftarrows V$:
$$
\xalignat 2
&\hskip -2em
\tau\!:\,V\to V^{{\sssize\dagger}*},
&&\tau\!:\,V^{{\sssize\dagger}*}\to V.
\mytag{3.6}
\endxalignat
$$
Similarly, for any vector $\bold v\in V$ the following function 
is defined:
$$
\hskip -2em
w_{\bold v}(f)=\overline{(f,\,\bold v)}\text{, \ where \ }
f\in V^*.
\mytag{3.7}
$$
The formula \mythetag{3.7} defines the canonical semilinear isomorphism
$V^{*\sssize\dagger}\rightleftarrows V$:
$$
\xalignat 2
&\hskip -2em
\tau\!:\,V\to V^{*\sssize\dagger},
&&\tau\!:\,V^{*\sssize\dagger}\to V.
\mytag{3.8}
\endxalignat
$$
All canonical isomorphisms \mythetag{3.4}, \mythetag{3.6}, and
\mythetag{3.8} can be extended to various tensor products of $V$,
$V^*$, $V^{\sssize\dagger}$, and $V^{*\sssize\dagger}$. In coordinate
representation they are expressed as passing to the conjugate numbers
in the coordinates of tensors.
\head
4. Conjugate spinors and spin-tensors.
\endhead
     The construction of the Hermitian conjugate space described in
the previous section~3 can be applied to the spinor spaces $S_p(M)$,
where $p\in M$. As a result we get the spaces $S^{\sssize\dagger}_p(M)$ 
and $S^{{\sssize\dagger}*}_p(M)=S^{*\sssize\dagger}_p(M)$ in addition to
$S_p(M)$ and $S^*_p(M)$. These new spaces are glued into smooth complex
bundles $S^{\sssize\dagger}\!M$ and $S^{{\sssize\dagger}*}\!M$ over the
base manifold $M$. They are called the {\it Hermitian conjugate spinor
bundles}.\par
     Like in \mythetag{2.3}, one can construct a tensor product
of multiple copies of Hermitian conjugate spinor spaces 
$S^{\sssize\dagger}_p(M)$ and $S^{{\sssize\dagger}*}_p(M)$:
$$
\hskip -2em
\bar S^\nu_\gamma(p,M)=\overbrace{S^{{\sssize\dagger}*}_p(M)\otimes
\ldots\otimes S^{{\sssize\dagger}*}_p(M)}^{\text{$\nu$ times}}
\otimes\underbrace{S^{\sssize\dagger}_p(M)\otimes\ldots\otimes 
S^{\sssize\dagger}_p(M)}_{\text{$\gamma$ times}}.
\mytag{4.1}
$$     
Then, using \mythetag{4.1}, one can extend the tensor product \mythetag{2.7}
as follows:
$$
\hskip -2em
S^{\,\alpha}_\beta\bar S^{\nu}_\gamma T^r_s(p,M)=S^{\,\alpha}_\beta(p,M)
\otimes\bar S^\nu_\gamma(p,M)\otimes\Bbb CT^r_s(p,M).
\mytag{4.2}
$$
Elements of the spaces \mythetag{4.2} are also called {\it spin-tensors}. The 
spaces \mythetag{4.2} with $p$ running over $M$ constitute a smooth
complex bundle over the base manifold $M$. Local and global smooth sections
of this bundle are called {\it spin-tensorial fields}.
Upon choosing some local chart $U$ of $M$ over which the spinor bundle
$SM$ trivializes one can express any spin-tensorial field $\bold Y$ in
the coordinate form
$$
\hskip -2em
Y^{i_1\ldots\,i_\alpha\bar i_1\ldots\,\bar i_\nu h_1\ldots\,h_r}_{j_1
\ldots\,j_\beta\bar j_1\ldots\,\bar j_\gamma k_1\ldots\, k_s}
=Y^{i_1\ldots\,i_\alpha\bar i_1\ldots\,\bar i_\nu h_1\ldots\,h_r}_{j_1
\ldots\,j_\beta\bar j_1\ldots\,\bar j_\gamma k_1\ldots\, k_s}(x^0,
x^1,x^2,x^3),
\mytag{4.3}
$$
where $x^0,\,x^1,\,x^2,\,x^3$ are local coordinates of a point $p\in M$.
The barred indices $\bar i_1,\,\ldots\,\,\bar i_\nu $ and $\bar j_1,\,
\ldots,\,\bar j_\gamma$ in \mythetag{4.3} are treated as separate variables
independent of $i_1,\,\ldots,\,i_\alpha$ and $j_1,\ldots,\,j_\beta$. In
some books (see \mycite{7}) dotted indices are used instead of barred
indices.\par
     The semilinear isomorphism \mythetag{2.9} and the semilinear
isomorphisms \mythetag{3.4}, \mythetag{3.6}, and  \mythetag{3.8} with
$V=S_p(M)$ induce the following semilinear isomorphisms of the spin-tensorial 
spaces \mythetag{4.2} and corresponding bundles:
$$
\align
&\hskip -2em
\tau\!:\,
S^{\,\alpha}_\beta\bar S^{\nu}_\gamma T^r_s(p,M)
\to S^{\nu}_\gamma\bar S^{\,\alpha}_\beta T^r_s(p,M),
\mytag{4.4}\\
\vspace{2ex}
&\hskip -2em
\tau\!:\,S^\alpha_\beta\bar S^\nu_\gamma T^r_s\!M
\to S^\nu_\gamma\bar S^\alpha_\beta T^r_s\!M.
\mytag{4.5}
\endalign
$$
In local coordinates the isomorphisms \mythetag{4.4} and \mythetag{4.5}
are expressed by the formula
$$
\hskip -2em
Z^{i_1\ldots\,i_\nu\bar i_1\ldots\,\bar i_\alpha h_1\ldots\,h_r}_{j_1
\ldots\,j_\gamma\bar j_1\ldots\,\bar j_\beta k_1\ldots\, k_s}=
\overline{Y^{\bar i_1\ldots\,\bar i_\alpha i_1\ldots\, i_\nu h_1\ldots
\,h_r}_{\bar j_1\ldots\,\bar j_\beta j_1\ldots\, j_\gamma k_1\ldots\, 
k_s}},
\mytag{4.6}
$$
where $\bold Z=\tau(\bold Y)$. The formula \mythetag{4.6} means 
that $\tau$ acts as the complex conjugation upon the components 
of spin-tensors simultaneously exchanging barred and non-barred 
spinor indices in upper and lower positions. 
\head
5. Local trivializations of the spinor bundle.
\endhead
    The spinor bundle $SM$, as it was introduced above in section~2,
is not the actual spinor bundle. It is yet an arbitrary $2$-dimensional
complex bundle over the space-time $M$. In order to get the actual 
spinor bundle we should relate $SM$ to $TM$ through the metric tensor
$\bold g$, which is the basic structure of the space-time manifold $M$.
\par
    Let $p_0$ be some point of $M$ and let $U$ be a local chart covering
that point $p_0$. Then any point $p\in U$ is represented by its local
coordinates 
$$
\hskip -2em
x^0=x^0(p),\quad x^1=x^1(p),\quad x^2=x^2(p),\quad x^3=x^3(p)
\mytag{5.1}
$$
and any tangent vector $\bold v\in T_pM$ is represented by its components
$v^0,\,v^1,\,v^2,\,v^3$:
$$
\hskip -2em
\bold v=v^0\,\frac{\partial}{\partial x^0}+v^1\,\frac{\partial}
{\partial x^1}+v^2\,\frac{\partial}{\partial x^2}+v^3\,\frac{\partial}
{\partial x^3}.
\mytag{5.2}
$$
The coordinate vectors
$$
\hskip -2em
\bold E_0=\frac{\partial}{\partial x^0},\quad 
\bold E_1=\frac{\partial}{\partial x^1},\quad 
\bold E_2=\frac{\partial}{\partial x^2},\quad 
\bold E_3=\frac{\partial}{\partial x^3}.
\mytag{5.3}
$$
in the expansion \mythetag{5.2} form the {\it holonomic moving frame}
associated with the local coordinates \mythetag{5.1} in $U$. The pair
$q=(p,\bold v)$ is a point of the tangent bundle $TM$ and $p=\pi(q)$,
see the formula \mythetag{2.1} above. Therefore, the quantities 
$$
\alignat 4
&\hskip -2em
x^0=x^0(q),\quad &&x^1=x^1(q),\quad &&x^2=x^2(q),\quad &&x^3=x^3(q),\\
\vspace{-1.4ex}
\mytag{5.4}\\
\vspace{-1ex}
&\hskip -2em
v^0=v^0(q),\quad &&v^1=v^1(q),\quad &&v^2=v^2(q),\quad &&v^3=v^3(q)
\endalignat
$$
determine the trivialization of $TM$ over $U$ in the holonomic frame
\mythetag{5.3}.
     Suppose that $\boldsymbol\Upsilon_0,\,\boldsymbol\Upsilon_2,
\,\boldsymbol\Upsilon_2,\,\boldsymbol\Upsilon_4$ are some smooth 
vector fields in $U$ being linearly independent at each point $p\in U$.
In this case we can write the expansion similar to \mythetag{5.2}:
$$
\hskip -2em
\bold v=v^0\,\boldsymbol\Upsilon_0+v^1\,\boldsymbol\Upsilon_1
+v^2\,\boldsymbol\Upsilon_2+v^3\,\boldsymbol\Upsilon_3.
\mytag{5.5}
$$
Replacing \mythetag{5.2} by the expansion \mythetag{5.5}, we say that
the quantities \mythetag{5.4} determine the trivialization of $TM$ over
$U$ in the {\it non-holonomic frame} $\boldsymbol\Upsilon_0,\,\boldsymbol
\Upsilon_1,\,\boldsymbol\Upsilon_2,\,\boldsymbol\Upsilon_3$.\par
     In the holonomic frame \mythetag{5.3} the metric tensor $\bold g$
is represented by a square $4\times 4$ matrix. In general case of a 
non-flat space-time $M$ this matrix is non-diagonal:
$$
\hskip -2em
g_{ij}=g_{j\,i}=\Vmatrix g_{00} &g_{01} &g_{02} &g_{03}\\
g_{01} &g_{11} &g_{12} &g_{13}\\ g_{02} &g_{12} &g_{22} &g_{23}\\
g_{03} &g_{13} &g_{23} &g_{33}\endVmatrix
\mytag{5.6}
$$
This matrix \mythetag{5.6} can be diagonalized at any fixed point $p_0
\in M$, but it 
cannot be diagonalized in a whole neighborhood of the point $p_0$ if
we use only holonomic frames, i\.\,e\. frames \mythetag{5.3} associated
with some local coordinates $x^0,\,x^1,\,x^2,\,x^3$. In the case of 
non-holonomic frames one can easily prove the following theorem.
\mytheorem{5.1} For any local chart $U\subset M$ and for any point 
$p_0\in U$ there is some smaller neighborhood $\tilde U\subset U$ 
of the point $p_0$ and there is some smooth non-holonomic frame
$\boldsymbol\Upsilon_0,\,\boldsymbol\Upsilon_1,\,\boldsymbol\Upsilon_2,
\,\boldsymbol\Upsilon_3$ within $\tilde U$ such that the metric tensor
$\bold g$ is represented by the standard diagonal matrix of the Minkowski-type 
metric:
$$
\hskip -2em
g_{ij}=g_{j\,i}=\Vmatrix\format \l&\quad\r&\quad\r&\quad\r\\
1 &0 &0 &0\\ 0 &-1 &0 &0\\ 0 &0 &-1 &0\\ 0 &0 &0 &-1\endVmatrix.
\mytag{5.7}
$$
\endproclaim
\mydefinition{5.1} A frame of four smooth vector fields $\boldsymbol
\Upsilon_0,\,\boldsymbol\Upsilon_1,\,\boldsymbol\Upsilon_2,\,
\boldsymbol\Upsilon_3$ in which the metric tensor $\bold g$ is given 
by the matrix \mythetag{5.7} is called an {\it orthonormal frame}.
\enddefinition
\mytheorem{5.2} Any two orthonormal frames of the space-time 
$\boldsymbol\Upsilon_0,\,\boldsymbol\Upsilon_1,\,\boldsymbol\Upsilon_2,
\,\boldsymbol\Upsilon_3$ and $\tilde{\boldsymbol\Upsilon}_0,\,
\tilde{\boldsymbol\Upsilon}_1,\,\tilde{\boldsymbol\Upsilon}_2,\,
\tilde{\boldsymbol\Upsilon}_3$ at the same point $p\in M$ are related 
to each other by some Lorentzian transition matrices $S\in\MatGrO(1,3)$
and $T\in\MatGrO(1,3)$:
$$
\xalignat 2
&\hskip -2em
\tilde{\boldsymbol\Upsilon}_i=\sum^3_{j=0}S^j_i\,\boldsymbol\Upsilon_j,
&&\boldsymbol\Upsilon_i=\sum^3_{j=0}T^j_i\,\tilde{\boldsymbol\Upsilon}_j.
\mytag{5.8}
\endxalignat
$$
The matrices $S$ and $T$ in the formulas \mythetag{5.8} are inverse 
to each other: $T=S^{-1}$.
\endproclaim
   Apart from the metric tensor $\bold g$, the space-time manifold 
$M$ carries other two structures: the {\it orientation} and the 
{\it polarization}. The orientation in $M$ means that we can distinguish
{\it right oriented frames} and {\it left oriented frames}. If two
frames $\boldsymbol\Upsilon_0,\,\boldsymbol\Upsilon_1,\,\boldsymbol
\Upsilon_2,\,\boldsymbol\Upsilon_3$ and $\tilde{\boldsymbol\Upsilon}_0,
\,\tilde{\boldsymbol\Upsilon}_1,\,\tilde{\boldsymbol\Upsilon}_2,\,
\tilde{\boldsymbol\Upsilon}_3$ are of the same orientation, both right
or both left, then the determinant of the transition matrix $S$ in 
\mythetag{5.8} is positive: $\det S>0$. Otherwise, for two frames 
of different orientation $\det S<0$. Remember that for any Lorentzian
matrix $S\in\MatGrO(1,3)$ we have $\det S=\pm\,1$. Lorentzian
matrices with $\det S=1$ form a subgroup in the Lorentz group
$\MatGrO(1,3)$. This is the {\it special Lorentz group\/} $\MatGrSO(1,3)$.
\par
     The polarization is a geometric structure of the space-time $M$ that
distinguishes {\tencyr\char '074}the Future{\tencyr\char '076} and 
{\tencyr\char '074}the Past{\tencyr\char '076} (see \mycite{4}). Remember
that a time-like vector is a tangent vector $\bold v\in T_p(M)$ at a point
$p\in M$ such that $g(\bold v,\bold v)>0$ with respect to the Minkowski
metric $\bold g$. At each point $p\in M$ the time-like vectors form the
interior of two cones (they are called the {\it light cones}). The
polarization marks one of these two cones as {\tencyr\char '074}the Future
cone{\tencyr\char '076} and the other as {\tencyr\char '074}the Past 
cone{\tencyr\char '076}. The polarization of the space-time $M$ is a smooth
structure. This means that for any smooth parametric curve $p=p(\tau)$ in
$M$ there is a smooth vector-valued function $\bold v=\bold v(\tau)$ on this
curve such that all its values $\bold v(\tau)$ are time-like vectors from
the Future cone. A Lorentzian matrix $S$ is called an {\it orthochronous}
matrix if $S^0_0>1$. Orthochronous Lorentzian matrices form a subgroup
in $\MatGrO(1,3)$. This subgroup is denoted as $\MatGrO^+(1,3)$. The {\it
special orthochronous Lorentz group\/} is defined as
$$
\hskip -2em
\MatGrSO^+(1,3)=\MatGrSO(1,3)\cap\MatGrO^+(1,3).
\mytag{5.9}
$$
Suppose that $\boldsymbol\Upsilon_0,\,\boldsymbol\Upsilon_1,\,\boldsymbol
\Upsilon_2,\,\boldsymbol\Upsilon_3$ and $\tilde{\boldsymbol\Upsilon}_0,
\,\tilde{\boldsymbol\Upsilon}_1,\,\tilde{\boldsymbol\Upsilon}_2,
\,\tilde{\boldsymbol\Upsilon}_3$ are two orthonormal frames. Then
$\boldsymbol\Upsilon_0$ and $\tilde{\boldsymbol\Upsilon}_0$ both are
time-like unit vectors. If they belong to the same light cone (both are 
in the the Future cone or both are in the Past cone), then the Lorentz
matrices $S$ and $T$ relating the frames $\boldsymbol\Upsilon_0,
\,\boldsymbol\Upsilon_1,\,\boldsymbol\Upsilon_2,\,\boldsymbol\Upsilon_3$
and $\tilde{\boldsymbol\Upsilon}_0,\,\tilde{\boldsymbol\Upsilon}_1,
\,\tilde{\boldsymbol\Upsilon}_2,\,\tilde{\boldsymbol\Upsilon}_3$ in 
\mythetag{5.8} are orthochronous Lorentzian matrices. Let's consider the
following Pauli matrices:
$$
\xalignat 3
&\boldsymbol\sigma_1=\Vmatrix 0 & 1\\1 & 0\endVmatrix
&&\boldsymbol\sigma_2=\Vmatrix 0 & -i\\i & 0\endVmatrix
&&\boldsymbol\sigma_3=\Vmatrix 1 & 0\\0 & -1\endVmatrix
\qquad
\mytag{5.10}
\endxalignat
$$
We complement the Pauli matrices \mythetag{5.10} with the unit matrix
$\boldsymbol\sigma_0$:
$$
\hskip -2em
\boldsymbol\sigma_0=\Vmatrix 1 & 0\\0 & 1\endVmatrix.
\mytag{5.11}
$$
Using the matrices \mythetag{5.10} and \mythetag{5.11}, one can \pagebreak
map each $4$-vector $\bold w\in\Bbb R^4$ to a Hermitian $2\times 2$
matrix by means of the following formula:
$$
\hskip -2em
\Vmatrix w^0\\ w^1\\ w^2\\ w^3\endVmatrix\mapsto h(\bold w)
=\sum^3_{m=1}w^m\,\boldsymbol\sigma_m=\Vmatrix w^0+w^3 & w^1-i\,w^2\\
\vspace{2ex} w^1+i\,w^2 & w^0-w^3\endVmatrix.
\mytag{5.12}
$$
Conversely, each hermitian matrix $\bold h$ can be mapped back to
a $4$-vector $\bold w\in\Bbb R^4$:
$$
\hskip -2em
\Vmatrix h^{11} & h^{12}\\ 
\vspace{2ex}\overline{h^{12}} & h^{22}
\endVmatrix\mapsto w(\bold h)=\frac{1}{2}
\Vmatrix h^{11}+h^{22}\\ 
\vspace{2ex} h^{12}+\overline{h^{12}}\\
\vspace{2ex} i\,h^{12}-i\,\overline{h^{12}}\\
\vspace{2ex} h^{11}-h^{22}
\endVmatrix.
\mytag{5.13}
$$
It is easy to see that the above mappings \mythetag{5.12} and
\mythetag{5.13} are inverse to each other. They define an isomorphism 
of two linear spaces 
$$
\hskip -2em
\CD
@>h>>\\
\vspace{-4ex}
\Bbb R^4@.\operatorname{\Bbb H\Bbb C}^{2\times 2}.\\
\vspace{-4.2ex}
@<<w< 
\endCD
\mytag{5.14}
$$
Note that the determinant of the matrix $\bold h=h(\bold w)$ coincides
with the square of the norm of $\bold w$ measured in the diagonal
Minkowski metric \mythetag{5.7}:
$$
\det h(\bold w)=(w^0)^2-(w^1)^2-(w^2)^2-(w^3)^2=g(\bold w,\bold w)=
\sum^3_{i=0}\sum^3_{j=0}g_{ij}\,w^i\,w^j.\quad
\mytag{5.15}
$$
Let $\bold U$ be a complex $2\times 2$ matrix and let $\det\bold U=1$. 
Then due to the isomorphism \mythetag{5.14} we can write the following 
equality:
$$
\hskip -2em
\bold U\,h(\bold w)\,\bold U^{\sssize\dagger}=h(\bold v).
\mytag{5.16}
$$
Here $\bold U^{\sssize\dagger}$ is the Hermitian transposition of 
$\bold U$. Due to the equality $\det\bold U=1$ and due to \mythetag{5.15}
for the components of the vector $\bold v$ in  \mythetag{5.16} we derive
$$
v^i=\sum^3_{j=0}S^i_j\,w^j.
\mytag{5.17}
$$
The components of the Lorentzian matrix $S$ in \mythetag{5.17} can be 
explicitly expressed through the components of the matrix $\bold U\in 
\MatGrSL(2,\Bbb C)$ in \mythetag{5.16} (see \mycite{8}). Here is the
formula for $S^0_0$. It is derived by direct calculations:
$$
\hskip -2em
S^0_0=\frac{|U^1_1|^2+|U^1_2|^2+|U^2_1|^2+|U^2_2|^2}{2}>0.
\mytag{5.18}
$$
We shall not write the explicit expressions for other components of
the matrix $S$, however, using them, we calculate the determinant of $S$:
$$
\pagebreak
\hskip -2em
\det S=(\det\bold U)^2\,(\det\bold U^{\sssize\dagger})^2=1. 
\mytag{5.19}
$$
From \mythetag{5.18} and \mythetag{5.19} we get $S\in\MatGrSO^+(1,3)$ 
(see \mythetag{5.9}). Thus, we have a mapping:
$$
\hskip -2em
\varphi\!:\,\MatGrSL(2,\Bbb C)\to\MatGrSO^+(1,3,\Bbb R).
\mytag{5.20}
$$
\mytheorem{5.3} The mapping $\varphi$ in \mythetag{5.20} is a group
homomorphism. Its kernel consists of the following two matrices
$$
\xalignat 2
&\hskip -2em
\boldsymbol\sigma_0=\Vmatrix 1 & 0\\0 & 1\endVmatrix,
&&-\boldsymbol\sigma_0=\Vmatrix -1 & 0\\0 & -1\endVmatrix.
\mytag{5.21}
\endxalignat
$$
Topologically, the mapping $\varphi$ is a two-sheeted non-ramified
covering of the $6$-dimen\-sional real manifold $\MatGrSO^+(1,3,\Bbb R)$
by the other $6$-dimensional real manifold $\MatGrSL(2,\Bbb C)$.
\endproclaim
More details concerning the theorem~\mythetheorem{5.3} and its proof can 
be found in \mycite{9}. This theorem is a base for defining the spinor
bundle over the space-time $M$. Remember that any bundle is determined 
by its local trivializations and by transition functions relating any
two trivializations with intersecting domains (see \mycite{10}). In the 
case of the tangent bundle $TM$ of the space-time manifold $M$ any local
chart $U\subset M$ equipped with an orthonormal frame
$\boldsymbol\Upsilon_0,\,\boldsymbol\Upsilon_1,\,
\boldsymbol\Upsilon_2,\,\boldsymbol\Upsilon_3$ determines a local
trivialization of the bundle $TM$ over $U$. Indeed, a point $q=(p,\bold v)$
of the open set $\pi^{-1}(U)$ in this case is associated with $8$ numbers 
$$
x^0,\ x^1,\ x^2,\ x^3,\ v^0,\ v^1,\ v^2,\ v^2,
$$
where $x^0,\,x^1,\,x^2,\,x^3$ are the local coordinates of the point $p$ in
the local chart $U$ and $v^0,\,v^1,\,v^2,\,v^3$ are the coefficients in the 
expansion \mythetag{5.5} for the vector $\bold v$. Such a trivialization of
$TM$ is called an {\it orthonormal trivialization}. The
theorem~\mythetheorem{5.1} says that orthonormal trivializations of $TM$
do exist and their domains $\pi^{-1}(U)$ cover $TM$. Therefore, the orthonormal trivializations are sufficient for defining the bundle $TM$ in
whole. Any two orthonormal trivializations of the tangent bundle are related
to each other by transition functions:
$$
\xalignat 2
&\hskip -2em
v^i=\sum^3_{j=0}S^i_j\,\tilde v^j,
&&\tilde v^i=\sum^3_{j=0}T^i_j\,v^j.
\mytag{5.22}
\endxalignat
$$
Due to theorem~\mythetheorem{5.2} the transition functions \mythetag{5.22}
are given by the components of two mutually inverse Lorentzian matrices
$S=S(p)$ and $T=T(p)$.\par
     An orthonormal frame $\boldsymbol\Upsilon_0,\,\boldsymbol\Upsilon_1,
\,\boldsymbol\Upsilon_2,\,\boldsymbol\Upsilon_3$ is called a positively
polarized right orthonormal frame if it is right frame in the sense of
the orientation in the space-time $M$ and if $\boldsymbol\Upsilon_0$ is a
time-like unit vector in the Future cone in the sense of the polarization 
in the space-time $M$. By means of the proper choice of sign
$\boldsymbol\Upsilon_i\to\pm\boldsymbol\Upsilon_i$ one can transform any
orthonormal frame into a positively polarized right orthonormal frame.
Therefore, the positively polarized right orthonormal trivializations of
$TM$ are sufficient to describe completely the tangent bundle $TM$. The
transition functions \mythetag{5.22} in this case are given by mutually
inverse Lorentzian matrices 
$$
\xalignat 2
&\hskip -2em
S(p)\in\MatGrSO^+(1,3,\Bbb R),
&&T(p)\in\MatGrSO^+(1,3,\Bbb R).
\mytag{5.23}
\endxalignat
$$
\mydefinition{5.2} A two-dimensional complex vector bundle $SM$ over the
four-dimensional real space-time manifold $M$ is called the {\it spinor
bundle} if for each positively polarized right orthonormal trivialization
$(U,\,\boldsymbol\Upsilon_0,\,\boldsymbol\Upsilon_1,\,\boldsymbol\Upsilon_2,
\,\boldsymbol\Upsilon_3)$ of the tangent bundle $TM$ there is a  
trivialization $(U,\,\boldsymbol\Psi_1,\,\boldsymbol\Psi_2)$ of $SM$ such 
that for any two positively polarized right orthonormal trivializations of 
the tangent bundle $TM$ related to each other by means of the formulas 
$$
\xalignat 2
&\hskip -2em
\tilde{\boldsymbol\Upsilon}_i=\sum^3_{j=0}S^j_i\,\boldsymbol\Upsilon_j,
&&\boldsymbol\Upsilon_i=\sum^3_{j=0}T^j_i\,\tilde{\boldsymbol\Upsilon}_j
\mytag{5.24}
\endxalignat
$$
the associated trivializations $(U,\,\boldsymbol\Psi_1,\,
\boldsymbol\Psi_2)$ and $(\tilde U,\,\tilde{\boldsymbol\Psi}_1,\,
\tilde{\boldsymbol\Psi}_2)$ of $SM$ are related to each other by means
of the formulas
$$
\xalignat 2
&\hskip -2em
\tilde{\boldsymbol\Psi}_i=\pm\sum^2_{j=1}\goth S^j_i\,
\boldsymbol\Psi_j,
&&\boldsymbol\Psi_i=\pm\sum^2_{j=1}\goth T^j_i\,
\tilde{\boldsymbol\Psi}_j,
\mytag{5.25}
\endxalignat
$$
where $S=\varphi(\goth S)$, $T=\varphi(\goth T)$, and $\varphi$ is the
group homomorphism \mythetag{5.20}.
\enddefinition
    The definition~\mythedefinition{5.2} is self-consistent due to the
theorem~\mythetheorem{5.3}. However, there is the uncertainty in sign
in the formulas \mythetag{5.25} because of the nontrivial kernel 
\mythetag{5.21} of the group homomorphism \mythetag{5.20}. Therefore,
there could be some global topological obstructions for the existence 
of the spinor bundle $SM$ over a particular manifold $M$. They are
discussed in \mycite{8}. Below we shall assume that the topology of the
actual physical space-time $M$ is such that the spinor bundle $SM$
over $M$ introduced by the definition~\mythedefinition{5.2} does exist
and is unique up to an isomorphism.
\subhead A remark\endsubhead The trivialization $(U,\,\boldsymbol\Psi_1,
\,\boldsymbol\Psi_2)$ of the spinor bundle $SM$ which is canonically
associated with a positively polarized right orthonormal trivialization
$(U,\,\boldsymbol\Upsilon_0,\,\boldsymbol\Upsilon_1,\,\boldsymbol\Upsilon_2,
\,\boldsymbol\Upsilon_3)$ of the tangent bundle $TM$ in the 
definition~\mythedefinition{5.2} is not unique. There are exactly two
trivializations $(U,\,\boldsymbol\Psi_1,\,\boldsymbol\Psi_2)$ and
$(U,\,-\boldsymbol\Psi_1,\,-\boldsymbol\Psi_2)$ associated with each
such trivialization $(U,\,\boldsymbol\Upsilon_0,\,\boldsymbol\Upsilon_1,\,
\boldsymbol\Upsilon_2,\,\boldsymbol\Upsilon_3)$ of $TM$.\par
     Let's return to the coordinate representations of spinors
and spin-tensors that we discussed in section~4. Below we shall assume 
that all coordinate representations of spin-tensorial fields \mythetag{4.3}
are relative to some non-holonomic positively polarized right orthonormal
frame $\boldsymbol\Upsilon_0,\,\boldsymbol\Upsilon_1,\,
\boldsymbol\Upsilon_2,\,\boldsymbol\Upsilon_3$ in $TM$ and its associated
frame $\boldsymbol\Psi_1,\,\boldsymbol\Psi_2$ in $SM$. Let's denote by
$\boldsymbol\eta^0,\,\boldsymbol\eta^1,\,\boldsymbol\eta^2,\,
\boldsymbol\eta^3$ the dual frame for $\boldsymbol\Upsilon_0,\,
\boldsymbol\Upsilon_1,\,\boldsymbol\Upsilon_2,\,\boldsymbol\Upsilon_3$.
This means that $\boldsymbol\eta^0,\,\boldsymbol\eta^1,\,
\boldsymbol\eta^2,\,\boldsymbol\eta^3$ are four covectorial fields
such that 
$$
\hskip -2em
\bigl(\boldsymbol\eta^i,\,\boldsymbol\Upsilon_j\bigr)
=\eta^i(\boldsymbol\Upsilon_j)=\delta^i_j.
\mytag{5.26}
$$
Similarly, let's denote by $\boldsymbol\vartheta^{\,1},\,
\boldsymbol\vartheta^{\,2}$ the dual frame for $\boldsymbol\Psi_1,
\,\boldsymbol\Psi_2$. Here we have
$$
\hskip -2em
\bigl(\boldsymbol\vartheta^{\,i},\,\boldsymbol\Psi_j\bigr)
=\vartheta^{\,i}(\boldsymbol\Psi_j)=\delta^i_j
\mytag{5.27}
$$
(see \mythetag{3.1} for comparison). And finally, let's denote
$$
\xalignat 2
&\hskip -2em
\bPsi_i=\tau(\boldsymbol\Psi_i),
&&\overline{\boldsymbol\vartheta}^{\,i}=\tau(\boldsymbol\vartheta^{\,i}),
\mytag{5.28}
\endxalignat
$$
where by $\tau$ we designate the semilinear mappings considered above 
(see \mythetag{2.9}, \mythetag{3.4}, \mythetag{3.6}, \mythetag{3.8},
\mythetag{4.4}, \mythetag{4.5}). The frames $\bPsi_1,\,\bPsi_2$ and
$\overline{\boldsymbol\vartheta}^{\,1},\,
\overline{\boldsymbol\vartheta}^{\,2}$ are dual to each other:
$$
\hskip -2em
\bigl(\overline{\boldsymbol\vartheta}^{\,i},\,\bPsi_j\bigr)
=\overline{\vartheta}^{\,i}(\boldsymbol\Psi_j)=\delta^i_j
\mytag{5.29}
$$
Now let's consider the following tensor products:
$$
\align
&\hskip -2em
\boldsymbol\Upsilon^{k_1\ldots\,k_s}_{h_1\ldots\,h_r}
=\boldsymbol\Upsilon_{h_1}\otimes\ldots\otimes\boldsymbol\Upsilon_{h_r}
\otimes\boldsymbol\eta^{k_1}\otimes\ldots\otimes\boldsymbol\eta^{k_s},
\mytag{5.30}\\
&\hskip -2em
\boldsymbol\Psi^{j_1\ldots\,j_\beta}_{i_1\ldots\,i_\alpha}
=\boldsymbol\Psi_{i_1}\otimes\ldots\otimes\boldsymbol\Psi_{i_\alpha}
\otimes\boldsymbol\vartheta^{\,j_1}\otimes\ldots\otimes
\boldsymbol\vartheta^{\,j_\beta},
\mytag{5.31}\\
&\hskip -2em
\bPsi^{\bar j_1\ldots\,\bar j_\gamma}_{\bar i_1\ldots\,\bar i_\nu}
=\bPsi_{\bar i_1}\otimes\ldots\otimes\bPsi_{\bar i_\nu}\otimes
\overline{\boldsymbol\vartheta}^{\,j_1}\otimes\ldots\otimes
\overline{\boldsymbol\vartheta}^{\,j_\gamma}.
\mytag{5.32}
\endalign
$$
Using \mythetag{5.30}, \mythetag{5.31}, and \mythetag{5.32}, we introduce
the following tensor product:
$$
\hskip -2em
\boldsymbol\Psi^{j_1\ldots\,j_\beta\bar j_1\ldots\,\bar j_\gamma
k_1\ldots\, k_s}_{i_1\ldots\,i_\alpha\bar i_1\ldots\,\bar i_\nu
h_1\ldots\,h_r}
=\boldsymbol\Psi^{j_1\ldots\,j_\beta}_{i_1\ldots\,i_\alpha}
\otimes\bPsi^{\bar j_1\ldots\,\bar j_\gamma}_{\bar i_1\ldots\,
\bar i_\nu}\otimes\boldsymbol\Upsilon^{k_1\ldots\,k_s}_{h_1\ldots\,h_r}.
\mytag{5.33}
$$
Then the coordinate representation \mythetag{4.3} of a spin-tensorial field
$\bold Y$ means that it is represented as an expansion in the basis of
spin-tensorial fields \mythetag{5.33}:
$$
\bold Y=\msum{2}\Sb i_1,\,\ldots,\,i_\alpha\\ j_1,\,\ldots,\,j_\beta
\endSb\msum{2}\Sb\bar i_1,\,\ldots,\,\bar i_\nu\\ \bar j_1,\,\ldots,\,
\bar j_\gamma\endSb\msum{3}\Sb h_1,\,\ldots,\,h_r\\ k_1,\,\ldots,\,k_s
\endSb
Y^{i_1\ldots\,i_\alpha\bar i_1\ldots\,\bar i_\nu h_1\ldots\,h_r}_{j_1
\ldots\,j_\beta\bar j_1\ldots\,\bar j_\gamma k_1\ldots\, k_s}\ 
\boldsymbol\Psi^{j_1\ldots\,j_\beta\bar j_1\ldots\,\bar j_\gamma
k_1\ldots\, k_s}_{i_1\ldots\,i_\alpha\bar i_1\ldots\,\bar i_\nu
h_1\ldots\,h_r}.\quad
\mytag{5.34}
$$
The coordinate representation \mythetag{4.6} of $\tau$ is based on
the expansion \mythetag{5.34} and on the formulas \mythetag{5.28}
defining $\bPsi_i$ and $\overline{\boldsymbol\vartheta}^{\,i}$.
Due to the duality (biorthogonality) relations \mythetag{5.26},
\mythetag{5.27}, and \mythetag{5.29} we can perform the contraction
of \mythetag{4.3} with respect to any pair within three groups of 
indices in \mythetag{4.3}, i\.\,e\. we can contract $h_m$ with $k_n$,
$i_m$ with $j_n$, and $\bar i_m$ with $\bar j_n$.
\head
6. The spin-metric.
\endhead
    Let's look at the definition~\mythedefinition{5.2} again. The transition
matrices $S$ and $T$ in \mythetag{5.24} belong to the special orthochronous
Lorentz group (see \mythetag{5.23}). This fact reflects three basic structures 
available in the space-time manifold $M$: the metric, the orientation, and 
the polarization. The matrices $\goth S$ and $\goth T$ 
in \mythetag{5.25} are not complex $2\times 2$ matrices of general form. 
They belong to the special linear group:
$$
\xalignat 2
&\hskip -2em
\goth S(p)\in\MatGrSL(2,\Bbb C),
&&\goth T(p)\in\MatGrSL(2,\Bbb C).
\mytag{6.1}
\endxalignat
$$
There is a special structure in $SM$ responsible for \mythetag{6.1}. 
Let's denote by $\bold d$ the skew-symmetric spin-tensor with the 
components
$$
\hskip -2em
d_{ij}=\Vmatrix 0 & 1\\-1 & 0\endVmatrix
\mytag{6.2}
$$
in each trivialization $(U,\,\boldsymbol\Psi_1,\,\boldsymbol\Psi_2)$ of
$SM$ associated with some positively polarized right orthonormal 
trivialization of $TM$. In other words, let's define $\bold d$ as
$$
\hskip -2em
\bold d=\boldsymbol\vartheta^{\,1}\otimes\boldsymbol\vartheta^{\,2}
-\boldsymbol\vartheta^{\,2}\otimes\boldsymbol\vartheta^{\,1}.
\mytag{6.3}
$$
This definition \mythetag{6.3} of $\bold d$ is self-consistent
because of the well-known identities valid for the matrix \mythetag{6.2}
and arbitrary $2\times 2$ matrices $\goth S$ and $\goth T$:
$$
\hskip -2em
\aligned
&\sum^2_{i=1}\sum^2_{j=1}d_{ij}\,\goth S^i_k\,\goth S^j_q=\det\goth S
\cdot d_{kq},\\
&\sum^2_{i=1}\sum^2_{j=1}d_{ij}\,\goth T^i_k\,\goth T^j_q=\det\goth T
\cdot d_{kq}.
\endaligned
\mytag{6.4}
$$
Due to \mythetag{6.1} and \mythetag{6.4} the formula \mythetag{6.3} for
$\bold d$ is compatible with the transition formulas \mythetag{5.25}. The
skew-symmetric spin-tensor $\bold d$ plays the same role for the spinor
bundle $SM$ as the metric tensor $\bold g$ for $TM$. In particular, the
matrix \mythetag{6.2} is used for lowering spinor indices of spin-tensors:
$$
\hskip -2em
Y^{\ldots\,\ldots\,\ldots}_{\ \ldots\,j\,\ldots}
=\sum^2_{i=1}Y^{\ldots\,i\,\ldots}_{\ldots\,\ldots\,\ldots}
\ d_{ij}.
\mytag{6.5}
$$
For this reason $\bold d$ is called the {\it spin-metric tensor}. The {\it
dual spin-metric tensor} is defined by the matrix inverse to \mythetag{6.2}:
$$
\hskip -2em
d^{\,ij}=\Vmatrix 0 & -1\\1 & 0\endVmatrix.
\mytag{6.6}
$$
It is denoted by the same symbol $\bold d$. Therefore, we have the formula
$$
\hskip -2em
\bold d=-\boldsymbol\Psi_1\otimes\boldsymbol\Psi_2
+\boldsymbol\Psi_2\otimes\boldsymbol\Psi_1.
\mytag{6.7}
$$
The dual spin-metric tensor given by the formula \mythetag{6.6} or by the
equivalent formula \mythetag{6.7} is used for raising spinor indices of
spin-tensors:
$$
\hskip -2em
Y^{\ldots\,j\,\ldots}_{\ldots\,\ldots\,\ldots}
=\sum^2_{i=1}Y^{\ldots\,\ldots\,\ldots}_{\ \ldots\,i\,\ldots}
\ d^{\,ij}.
\mytag{6.8}
$$
Applying $\tau$ to \mythetag{6.3} and \mythetag{6.7} we get the {\it
conjugate spin-metric tensors}:
$$
\align
&\hskip -2em
\bbd=\overline{\boldsymbol\vartheta}^{\,1}
\otimes\overline{\boldsymbol\vartheta}^{\,2}
-\overline{\boldsymbol\vartheta}^{\,2}\otimes
\overline{\boldsymbol\vartheta}^{\,1},
\mytag{6.9}\\
&\hskip -2em
\bbd=-\bPsi_1\otimes\bPsi_2+\bPsi_2\otimes\bPsi_1.
\mytag{6.10}
\endalign
$$
The conjugate spin-metric tensors \mythetag{6.9} and \mythetag{6.10}
are represented by the same matrices \mythetag{6.2} and \mythetag{6.6},
but we use the barred $d$ for denoting their components:
$$
\xalignat 2
&\hskip -2em
\bd_{\,\bar i\kern 0.5pt\bar j}=\Vmatrix 0 & 1\\-1 & 0\endVmatrix,
&&\bd^{\,\bar i\kern 0.5pt \bar j}=\Vmatrix 0 & -1\\1 & 0\endVmatrix.
\mytag{6.11}
\endxalignat
$$
The matrices \mythetag{6.11} are used for raising and lowering the
barred spinor indices:
$$
\pagebreak
\align
&\hskip -2em
Y^{\ldots\,\ldots\,\ldots}_{\ \ldots\,\bar j\,\ldots}
=\sum^2_{\bar i=1}Y^{\ldots\,\bar i\,\ldots}_{\ldots\,\ldots\,\ldots}
\ \bd_{\,\bar i\kern 0.5pt\bar j}\,,
\mytag{6.12}\\
&\hskip -2em
Y^{\ldots\,\bar j\,\ldots}_{\ldots\,\ldots\,\ldots}
=\sum^2_{i=1}Y^{\ldots\,\ldots\,\ldots}_{\ \ldots\,\bar i\,\ldots}
\ \bd^{\,\bar i\kern 0.5pt\bar j}.
\mytag{6.13}
\endalign
$$
These formulas \mythetag{6.12} and \mythetag{6.13} are analogs of
the formulas \mythetag{6.5} and \mythetag{6.8}.
\head
7. Tensors to spin-tensors conversion.
\endhead
     As we can see in \mythetag{4.3}, a general spin-tensorial field
$\bold Y$ has three groups of indices: regular spinor indices, barred
spinor indices, and tensorial indices. In this section we shall show
that tensorial indices can be converted into spinor and barred spinor
indices. For this purpose let's return to the formula \mythetag{5.16}.
Let's denote $h=h(\bold v)$, $\tilde h=h(\bold w)$, and let's
substitute $\bold U=\goth S$ into \mythetag{5.16}. Then we find that
\mythetag{5.16} is equivalent to the following equality:
$$
\hskip -2em
h^{i\kern 0.5pt\bar i}=\sum^2_{\bar q=1}\sum^2_{q=1}
\goth S^i_q\ \tilde h^{q\kern 0.5pt\bar q}\ 
\overline{\goth S^{\,\bar i}_{\bar q}}.
\mytag{7.1}
$$
This formula \mythetag{7.1} coincides with the transformation rule for
the components of a spin-tensor with one regular upper index and one 
barred upper index under the change of frame given by the formulas
\mythetag{5.25}. The sign uncertainty in \mythetag{5.25} does not affect
the formula \mythetag{7.1}.\par
     Similarly, if we denote $\bold w=\tilde\bold v$, we can write the
formula \mythetag{5.17} in the form coinciding with the first formula
in \mythetag{5.22}:
$$
\hskip -2em
v^i=\sum^3_{j=0}S^i_j\,\tilde v^j.
\mytag{7.2}
$$
The formula \mythetag{7.2} expresses the transformation rule for the
components of a tensor with one upper index under the change of a frame
given by the formulas \mythetag{5.24}. Due to \mythetag{7.1} and 
\mythetag{7.2} the upper mapping $h$ in \mythetag{5.14} is interpreted
as the mapping
$$
\hskip -2em
h\!:\,T_p(M)\to S_p(M)\otimes S^{{\sssize\dagger}*}_p(M).
\mytag{7.3}
$$
Note that $h(\bold w)$ in \mythetag{5.12} is a Hermitian matrix. For the
spin-tensor $\bold h$ represented by this matrix in \mythetag{7.1} this
yields $\tau(\bold h)=\bold h$, where $\tau$ is the mapping defined by the
formula \mythetag{4.6}. Indeed, we can easily verify that
$$
\hskip -2em
h^{i\kern 0.5pt\bar i}=\overline{h^{\bar i\kern 0.5pt i}}
\mytag{7.4}
$$
Due to \mythetag{7.4} we can take the Hermitian symmetric part of the 
tensor product in the right hand side of \mythetag{7.3} and thus write
\mythetag{5.14} as follows:
$$
\hskip -2em
\CD
@>h>>\\
\vspace{-4ex}
T_p(M)@.\HSym(S_p(M)\otimes S^{{\sssize\dagger}*}_p(M))\\
\vspace{-4.2ex}
@<<w< 
\endCD
\mytag{7.5}
$$
The upper mapping $h$ in \mythetag{7.5} is given by a special spin-tensorial
field $\bold G$ with two upper spinor indices and one lower tensorial index:
$$
\hskip -2em
h^{i\kern 0.5pt\bar i}=\sum^3_{q=0}G^{i\kern 0.5pt\bar i}_q\,w^q.
\mytag{7.6}
$$
Its components $G^{i\kern 0.5pt\bar i}_j$ in \mythetag{7.6} are called 
{\it Infeld-van der Waerden symbols}. The numeric values of these symbols
are calculated through the components of Pauli matrices \mythetag{5.10}
and \mythetag{5.11} because of the formula \mythetag{5.12}:
$$
\xalignat 4
&\hskip -2em
G^{11}_0=1, &&G^{11}_1=0, &&G^{11}_2=0,  &&G^{11}_3=1,\\
&\hskip -2em
G^{12}_0=0, &&G^{12}_1=1, &&G^{12}_2=-i, &&G^{12}_3=0,\\
\vspace{-1.4ex}
\mytag{7.7}\\
\vspace{-1.4ex}
&\hskip -2em
G^{21}_0=0, &&G^{21}_1=1, &&G^{21}_2=i,  &&G^{21}_3=0,\\
&\hskip -2em
G^{22}_0=1, &&G^{22}_1=0, &&G^{22}_2=0,  &&G^{22}_3=-1.\\
\endxalignat
$$
The inverse mapping $w$ in \mythetag{7.5} is also given by some special
spin-tensorial field. We use the same symbol $\bold G$ for this field:
$$
\hskip -2em
w^q=\sum^2_{i=1}\sum^2_{\bar i=1}G^{\,q}_{i\kern 0.5pt\bar i}
\ h^{i\kern 0.5pt\bar i}.
\mytag{7.8}
$$
Its components $G^{\,q}_{i\kern 0.5pt\bar i}$ in \mythetag{7.8} are 
called the {\it inverse Infeld-van der Waerden symbols}. By means of
\mythetag{5.13} one can calculate them in explicit form:
$$
\xalignat 4
&\hskip -2em
G^{\,0}_{11}=\frac{1}{2}, &&G^{\,0}_{12}=0, &&G^{\,0}_{21}=0,  
&&G^{\,0}_{22}=\frac{1}{2},\\
&\hskip -2em
G^{\,1}_{11}=0, &&G^{\,1}_{12}=\frac{1}{2}, &&G^{\,1}_{21}=\frac{1}{2},  
&&G^{\,1}_{22}=0,\\
\vspace{-1.4ex}
\mytag{7.9}\\
\vspace{-1.4ex}
&\hskip -2em
G^{\,2}_{11}=0, &&G^{\,2}_{12}=\frac{i}{2}, &&G^{\,2}_{21}=-\frac{i}{2},  
&&G^{\,2}_{22}=0,\\
&\hskip -2em
G^{\,3}_{11}=\frac{1}{2}, &&G^{\,3}_{12}=0, &&G^{\,3}_{21}=0,  
&&G^{\,3}_{22}=-\frac{1}{2}.
\endxalignat
$$
The following symmetry properties of the Infeld-van der Waerden symbols
are obvious. They are derived from \mythetag{7.7} and \mythetag{7.9}:
$$
\xalignat 2
&\hskip -2em
G^{i\kern 0.5pt\bar i}_q=\overline{G^{\bar i\kern 0.5pt i}_q},
&&G^{\,q}_{i\kern 0.5pt\bar i}=\overline{G^{\,q}_{\bar i\kern 0.5pt i}}.
\mytag{7.10}
\endxalignat
$$
The transformation \mythetag{7.6} can be applied not only to a tensorial
field with one upper index.  We can apply this transformation 
to any spin-tensorial with at least one non-spinor upper index. By analogy
to \mythetag{7.6} we can do it in the following way
$$
\hskip -2em
\hat Y^{\ldots\,i\ldots\,\bar i\ldots}_{\ \ldots\,\ldots\,\ldots}
=\sum^3_{q=0}Y^{\ldots\,q\,\ldots}_{\ldots\,\ldots\,\ldots}
\ G^{i\kern 0.5pt\bar i}_q.
\mytag{7.11}
$$
In \mythetag{7.11} with the use of the Infeld-van der Waerden symbols
\mythetag{7.7} the tensorial index $q$ is converted into the pair of 
spinor indices $i$ and $\bar i$. The spin-tensorial field $\bold Y$ 
in \mythetag{7.11} can be recovered from $\hat\bold Y$ with the use of
the inverse Infeld-van der Waerden symbols \mythetag{7.9} by means
of the formula
$$
\hskip -2em
Y^{\ldots\,q\,\ldots}_{\ldots\,\ldots\,\ldots}
=\sum^2_{i=1}\sum^2_{\bar i=1}\hat Y^{\ldots\,i\ldots\,\bar i
\ldots}_{\ \ldots\,\ldots\,\ldots}\ G^{\,q}_{i\kern 0.5pt\bar i}\,.
\mytag{7.12}
$$
Similarly, if we apply the transformation \mythetag{7.12} to an arbitrary
spin-tensorial field $\hat\bold Y$ with two upper spinor indices $i$ and
$\bar i$, then we can recover it from $\bold Y$ by means of the formula
\mythetag{7.11}. The Infeld-van der Waerden symbols satisfy the identities
$$
\xalignat 2
&\hskip -2em
\sum^2_{i=1}\sum^2_{\bar i=1}G^{i\kern 0.5pt\bar i}_p\ 
G^{\,q}_{i\kern 0.5pt\bar i}=\delta^q_p,
&&\sum^3_{q=0}G^{i\kern 0.5pt\bar i}_q\ 
G^{\,q}_{j\kern 0.5pt\bar j}=\delta^i_j
\,\delta^{\bar i}_{\bar j}.\qquad
\mytag{7.13}
\endxalignat
$$
The identities \mythetag{7.13} are compatible with \mythetag{7.10}.
They are derived by means of direct calculations from \mythetag{7.7} 
and \mythetag{7.9} and used in order to prove that the transformations
\mythetag{7.11} and \mythetag{7.12} are inverse to each other.\par
     As an example of their usage, let's apply the conversion formulas
\mythetag{7.11} and \mythetag{7.12} to the metric tensor $\bold g$. The
result is given by the following identities:
$$
\gather
\hskip -2em
\sum^3_{p=0}\sum^3_{q=0}g_{p\,q}\,G^{\,p}_{i\kern 0.5pt\bar i}
\ G^{\,q}_{j\kern 0.5pt\bar j}=\frac{d_{ij}\,
\bd_{\,\bar i\kern 0.5pt\bar j}}{2},
\mytag{7.14}\\
\hskip -2em
\sum^2_{i=1}\sum^2_{j=1}\sum^2_{\bar i=1}\sum^2_{\bar j=1}
\frac{d_{ij}\,\bd_{\,\bar i\kern 0.5pt\bar j}}{2}
\ G^{i\kern 0.5pt\bar i}_p\ G^{j\kern 0.5pt\bar j}_q=g_{p\,q}.
\mytag{7.15}
\endgather
$$
Applying the formulas \mythetag{7.11} and \mythetag{7.12} to the
dual metric tensor, we get
$$
\gather
\hskip -2em
\sum^3_{p=0}\sum^3_{q=0}g^{p\,q}\ G^{i\kern 0.5pt\bar i}_p
\ G^{j\kern 0.5pt\bar j}_q
=2\,d^{\,ij}\,\bd^{\,\bar i\kern 0.5pt \bar j},
\mytag{7.16}\\
\hskip -2em
\sum^2_{i=1}\sum^2_{j=1}\sum^2_{\bar i=1}\sum^2_{\bar j=1}
2\,d^{\,ij}\,\bd^{\,\bar i\kern 0.5pt \bar j}
\,G^{\,p}_{i\kern 0.5pt\bar i}
\ G^{\,q}_{j\kern 0.5pt\bar j}=g^{p\,q}.
\mytag{7.17}
\endgather
$$
The identities \mythetag{7.14}, \mythetag{7.15}, \mythetag{7.16}, and
\mythetag{7.17} are derived from \mythetag{7.7} and \mythetag{7.9} by
direct calculations.\par
     {\bf A remark}. The mapping \mythetag{7.3} can be extended to the
complexified tangent space $\Bbb CT_p(M)$. In this case the symmetry
condition \mythetag{7.4} is not valid and we have
$$
\hskip -2em
\CD
@>h>>\\
\vspace{-4ex}
\Bbb CT_p(M)@.S_p(M)\otimes S^{{\sssize\dagger}*}_p(M).\\
\vspace{-4.2ex}
@<<w< 
\endCD
\mytag{7.18}
$$
This formula \mythetag{7.18} is the replacement of the formula
\mythetag{7.5} for the case of the complexified tangent space 
$\Bbb CT_p(M)$.
\head
8. Dirac spinors.
\endhead
    The definition~\mythedefinition{5.2} introducing the concept of the
spinor bundle over the space-time manifold $M$ is a self-consistent and
self-sufficient definition. However, the transition matrices $S$ and $T$
in this definition are restricted to the special orthochronous Lorentz
group $\MatGrSO^+(1,3,\Bbb R)$ (see \mythetag{5.23}). In order to extend
$\MatGrSO^+(1,3,\Bbb R)$ up to the complete Lorentz group $\MatGrO(1,3,
\Bbb R)$ one should add the following two matrices to it:
$$
\xalignat 2
&\theta=\Vmatrix -1 & 0 & 0 & 0\\ 0 & 1 & 0 & 0\\ 0 & 0 & 1 & 0\\
0 & 0 & 0 & 1\endVmatrix,
&&P=\Vmatrix 1 & 0 & 0 & 0\\ 0 & -1 & 0 & 0\\ 0 & 0 & -1 & 0\\
0 & 0 & 0 & -1\endVmatrix.
\endxalignat
$$
For this purpose the other spinor bundle is introduced. It is 
constructed as the direct sum of $SM$ and its Hermitian conjugate 
bundle $S^{\sssize\dagger}\!M$:
$$
\hskip -2em
DM=SM\oplus S^{\sssize\dagger}\!M.
\mytag{8.1}
$$
The elements of this $4$-dimensional complex bundle \mythetag{8.1}
are called {\it Dirac spinors}. Though the Dirac spinors are very 
popular in physics, we shall not consider them in this paper since 
due to \mythetag{8.1} they are not primary objects --- they are 
reduced to the standard $2$-components spinors from $SM$.
\head
9. Composite spin-tensorial bundles.
\endhead
      Let's remember the section~4. There we have introduced the 
spin-tensorial bundle $S^\alpha_\beta\bar S^\nu_\gamma T^m_n\!M$.
It is composed by spin tensors of the type $(\alpha,\beta|\nu,\gamma|m,n)$.
Suppose that we have several bundles of that sort over the space-time $M$. 
Let's denote them 
$$
\hskip -2em
S^{\alpha_1}_{\beta_1}\!\bar S^{\nu_1}_{\gamma_1}
T^{m_1}_{n_1}\!M,\,\ldots,\,
S^{\alpha_J}_{\beta_J}\!\bar S^{\nu_J}_{\gamma_J}
T^{m_J}_{n_J}\!M.
\mytag{9.1}
$$
In addition to \mythetag{9.1} we consider several tensor bundles
$$
\hskip -2em
T^{r_1}_{s_1}\!M,\,\ldots,\,T^{r_Q}_{s_Q}\kern -1pt M.
\mytag{9.2}
$$
Then we construct the direct sum of the bundles \mythetag{9.1} and
\mythetag{9.2}:
$$
N=S^{\alpha_1}_{\beta_1}\!\bar S^{\nu_1}_{\gamma_1}
T^{m_1}_{n_1}\!M\oplus\ldots\oplus
S^{\alpha_J}_{\beta_J}\!\bar S^{\nu_J}_{\gamma_J}
T^{m_J}_{n_J}\!M\oplus T^{r_1}_{s_1}\!M\oplus
\ldots\oplus T^{r_Q}_{s_Q}\kern -1pt M\quad
\mytag{9.3}
$$
We shall call $N$ in \mythetag{9.3} the {\it composite spin-tensorial
bundle}. By definition a point $q$ of the composite spin-tensorial bundle
\mythetag{9.3} is a list
$$
\hskip -2em
q=(p,\,\bold S[1],\,\ldots,\,\bold S[J],\,\bold T[1],\,\ldots,\,\bold T[Q]),
\mytag{9.4}
$$
where $p$ is a point of the space-time $M$, $\bold S[1],\,\ldots,\,\bold
S[J]$ are spin-tensors of the types $(\alpha_1,\beta_1|\nu_1,\gamma_1
|m_1,n_1),\,\ldots,\,(\alpha_J,\beta_J|\nu_J,\gamma_J|m_J,n_J)$, and 
$\bold T[1],\,\ldots,\,\bold T[Q]$ are tensors of the types $(r_1,s_1),
\,\ldots,\,(r_Q,s_Q)$ at the point $p$. Comparing \mythetag{9.4} with
\mythetag{1.2} and \mythetag{1.3}, we  find that the composite tensor
bundle is a suitable geometric object for choosing it as the domain of 
the functions \mythetag{1.2} and \mythetag{1.3}. However, it is more
convenient to extend this domain. Indeed, note that the real tensor
bundles \mythetag{9.2} can be extended up to the complexified tensor 
bundles
$$
\hskip -2em
\Bbb CT^{r_1}_{s_1}\!M,\,\ldots,\,\Bbb CT^{r_Q}_{s_Q}\kern -1pt M.
\mytag{9.5}
$$
Then remember that the complexified tensor bundles \mythetag{9.5} 
coincide with the following spin-tensorial bundles (see formula
\mythetag{4.2}):
$$
\hskip -2em
S^0_0\bar S^0_0T^{r_1}_{s_1}\!M,\,\ldots,
\,S^0_0\bar S^0_0T^{r_Q}_{s_Q}\kern -1pt M
\mytag{9.6}
$$
Due to \mythetag{9.6} we can replace the direct sum \mythetag{9.3}
by the other direct sum 
$$
N=S^{\alpha_1}_{\beta_1}\!\bar S^{\nu_1}_{\gamma_1}
T^{m_1}_{n_1}\!M\oplus\ldots\oplus
S^{\alpha_{J+Q}}_{\beta_{J+Q}}\!\bar S^{\nu_{J+Q}}_{\gamma_{J+Q}}
T^{m_{J+Q}}_{n_{J+Q}}\!M,
\mytag{9.7}
$$
where $\alpha_P=\beta_P=\nu_P=\gamma_P=0$, $m_P=r_{P-J}$, $n_P=s_{P-J}$
for $J<P\leqslant J+Q$. Passing from \mythetag{9.3} to the composite tensor
bundle \mythetag{9.7}, we can treat in a more uniform way the components 
of the list \mythetag{9.4}. Indeed, now we can write
$$
\hskip -2em
q=(p,\,\bold S[1],\,\ldots,\,\bold S[J+Q]),
\mytag{9.8}
$$
where $\bold T[1]=\bold S[J+1],\,\ldots,\,\bold T[Q]=\bold S[J+Q]$.
In other words, upon complexification, we treat the tensors 
$\bold T[1],\,\ldots,\,\bold T[Q]$ as the spin-tensors $\bold S[J+1],
\,\ldots,\,\bold S[J+Q]$. The {\it canonical projection\/} of the
composite bundle \mythetag{9.7} is defined as the map 
$$
\hskip -2em
\pi\!:\,N\to M
\mytag{9.9}
$$ 
that takes a point $q$ of the form \mythetag{9.8} to the point $p\in M$.
Despite the complexification, the complex bundle $N$ is a real smooth
manifold. Its dimension 
$$
\hskip -2em
\dim_{\,\Bbb R}N=2\sum^{J+Q}_{i=1}2^{\,\alpha_i+\nu_i+\beta_i
+\gamma_i}\ 4^{m_i+n_i}+4.
\mytag{9.10}
$$
Below we shall use local charts of some special sort in $N$. In order to
construct a local chart of that sort one should choose some local chart 
$U$ of $M$ equipped with a positively polarized right oriented orthonormal
frame $\boldsymbol\Upsilon_0,\,\boldsymbol\Upsilon_1,\,
\boldsymbol\Upsilon_2,\,\boldsymbol\Upsilon_3$ and with its associated
spinor frame $\boldsymbol\Psi_1,\,\boldsymbol\Psi_2$. The first part of 
the local coordinates of a point $q$ in \mythetag{9.8} is the local
coordinates of its image $p=\pi(q)$ under the projection \mythetag{9.9}:
$$
\hskip -2em
x^0,\,x^1,\,x^2,\,x^3.
\mytag{9.11}
$$
The second part of the local coordinates of $q$ is composed by the components 
of the spin-tensors $\bold S[1],\,\ldots,\,\bold S[J+Q]$ 
in \mythetag{9.8}:
$$
\hskip -6em
\gathered
S^{1\kern 1pt\ldots\,1\,1\kern 1pt\ldots\,1\,0\kern 1pt\ldots\,
0}_{1\kern 1pt\ldots\,1\,1\kern 1pt\ldots\,1\,0\kern 1pt\ldots\,0}[1],
\,\ldots,\,
S^{2\kern 1pt\ldots\,2\,2\kern 1pt\ldots\,2\,3\kern 1pt\ldots\,
3}_{2\kern 1pt\ldots\,2\,2\kern 1pt\ldots\,2\,3\kern 1pt\ldots\,3}[1],
\,\ldots\\
\ldots,\,
S^{1\kern 1pt\ldots\,1\,1\kern 1pt\ldots\,1\,0\kern 1pt\ldots\,
0}_{1\kern 1pt\ldots\,1\,1\kern 1pt\ldots\,1\,0\kern 1pt\ldots\,0}[J+Q],
\,\ldots,\,
S^{2\kern 1pt\ldots\,2\,2\kern 1pt\ldots\,2\,3\kern 1pt\ldots\,
3}_{2\kern 1pt\ldots\,2\,2\kern 1pt\ldots\,2\,3\kern 1pt\ldots\,3}[J+Q].
\endgathered
\mytag{9.12}
$$
Each part of the variables yields its own contribution to the total
dimension of $N$ in \mythetag{9.10}. Note that \mythetag{9.12} are
complex variables. For this reason the contribution of \mythetag{9.12}
in \mythetag{9.10} is taken with the factor $2$.\par
     Under a change of local charts the coordinates \mythetag{9.11}
are transformed traditionally by means of the transition functions
$$
\xalignat 2
&\hskip -2em
\cases
\tx^0=\tx^1(x^0,\,\ldots,x^3),\\
. \ . \ . \ .\ . \ . \ . \ . \ . \ . \ 
. \ . \ . \ . \ .\\
\tx^3=\tx^n(x^0,\,\ldots,x^3).
\endcases
&&\cases
x^0=\tx^1(\tx^0,\,\ldots,\tx^3),\\
. \ . \ . \ .\ . \ . \ . \ . \ . \ . \ 
. \ . \ . \ . \ .\\
x^3=\tx^n(\tx^0,\,\ldots,\tx^3).
\endcases\quad
\mytag{9.13}
\endxalignat
$$
The coordinates \mythetag{9.12} are transformed as the components
of spin-tensors
$$
\align
&\hskip -4em
\left\{
\aligned
&\tilde S^{i_1\ldots\,i_\alpha\bar i_1\ldots\,\bar i_\nu
h_1\ldots\,h_m}_{j_1\ldots\,j_\beta\bar j_1\ldots\,\bar j_\gamma
k_1\ldots\,k_n}[P]
=\dsize\msum{2}\Sb a_1,\,\ldots,\,a_\alpha\\ b_1,\,\ldots,\,b_\beta\endSb
\dsize\msum{2}\Sb \bar a_1,\,\ldots,\,\bar a_\nu\\ 
\bar b_1,\,\ldots,\,\bar b_\gamma\endSb
\dsize\msum{3}
\Sb c_1,\,\ldots,\,c_m\\ d_1,\,\ldots,\,d_n\endSb
\goth T^{\,i_1}_{a_1}\ldots\,\goth T^{\,i_\alpha}_{a_\alpha}\,\times\\
&\kern 40pt
\times\,\goth S^{b_1}_{j_1}\ldots\,\goth S^{b_\beta}_{j_\beta}
\ \overline{\goth T^{\,\bar i_1}_{\bar a_1}}
\ldots\,\overline{\goth T^{\,\bar i_\nu}_{\bar a_\nu}}\ 
\ \overline{\goth S^{\,\bar b_1}_{\bar j_1}}
\ldots\,\overline{\goth S^{\,\bar b_\gamma}_{\bar j_\gamma}}\ 
T^{h_1}_{c_1}\ldots\,T^{h_m}_{c_m}\,\times\\
\vspace{1.5ex}
&\kern 60pt
\times\,S^{\,d_1}_{k_1}\ldots\,S^{\,d_n}_{k_n}\ 
S^{\,a_1\ldots\,a_\alpha\,\bar a_1\ldots\,\bar a_\nu\,
c_1\ldots\,c_m}_{\,b_1\ldots\,b_\beta\,\bar b_1\ldots\,\bar b_\gamma
\,d_1\ldots\,d_n}[P],
\endaligned
\right.
\mytag{9.14}\\
\displaybreak
&\hskip -4em
\left\{
\aligned
&S^{i_1\ldots\,i_\alpha\bar i_1\ldots\,\bar i_\nu
h_1\ldots\,h_m}_{j_1\ldots\,j_\beta\bar j_1\ldots\,\bar j_\gamma
k_1\ldots\,k_n}[P]
=\dsize\msum{2}\Sb a_1,\,\ldots,\,a_\alpha\\ b_1,\,\ldots,\,b_\beta\endSb
\dsize\msum{2}\Sb \bar a_1,\,\ldots,\,\bar a_\nu\\ 
\bar b_1,\,\ldots,\,\bar b_\gamma\endSb
\dsize\msum{3}
\Sb c_1,\,\ldots,\,c_m\\ d_1,\,\ldots,\,d_n\endSb
\goth S^{\,i_1}_{a_1}\ldots\,\goth S^{\,i_\alpha}_{a_\alpha}\,\times\\
&\kern 40pt
\times\,\goth T^{b_1}_{j_1}\ldots\,\goth T^{b_\beta}_{j_\beta}
\ \overline{\goth S^{\,\bar i_1}_{\bar a_1}}
\ldots\,\overline{\goth S^{\,\bar i_\nu}_{\bar a_\nu}}\ 
\ \overline{\goth T^{\,\bar b_1}_{\bar j_1}}
\ldots\,\overline{\goth T^{\,\bar b_\gamma}_{\bar j_\gamma}}\ 
S^{h_1}_{c_1}\ldots\,S^{h_m}_{c_m}\,\times\\
\vspace{1.5ex}
&\kern 60pt
\times\,T^{\,d_1}_{k_1}\ldots\,T^{\,d_n}_{k_n}\ 
\tilde S^{\,a_1\ldots\,a_\alpha\,\bar a_1\ldots\,\bar a_\nu\,
c_1\ldots\,c_m}_{\,b_1\ldots\,b_\beta\,\bar b_1\ldots\,\bar b_\gamma
\,d_1\ldots\,d_n}[P],
\endaligned
\right.
\mytag{9.15}
\endalign
$$
where $\alpha=\alpha_P$, $\beta=\beta_P$, $\nu=\nu_P$, $\gamma
=\gamma_P$, $m=m_P$, $n=n_P$, and the integer number $P$ runs from 
$1$ to $J+Q$. The coordinates \mythetag{9.12} are complex numbers. 
Actually, their real and imaginary parts form local coordinates 
of a point $q$ of the real manifold $N$. However, using the complex 
numbers \mythetag{9.12} is preferable at least because the transition
functions \mythetag{9.14} and \mythetag{9.15} look more simple in 
terms of these complex numbers.\par
    The components of the transition matrices $\goth S$, $\goth T$, $S$,
and $T$ in \mythetag{9.14} and \mythetag{9.15} are frame relative, not 
coordinate relative. For this reason they do not depend on the transition 
functions in \mythetag{9.13}. The components of the matrices $\goth S
\in\MatGrSL(2,\Bbb C)$ and $\goth T\in\MatGrSL(2,\Bbb C)$ are taken from
the frame relationships
$$
\xalignat 2
&\hskip -2em
\tilde{\boldsymbol\Psi}_i=\sum^2_{j=1}\goth S^j_i\,
\boldsymbol\Psi_j,
&&\boldsymbol\Psi_i=\sum^2_{j=1}\goth T^j_i\,
\tilde{\boldsymbol\Psi}_j.
\mytag{9.16}
\endxalignat
$$
The Lorentzian matrices $S=\varphi(\goth S)\in\MatGrSO^+(1,3,\Bbb R)$
and $T=\varphi(\goth T)\in\MatGrSO^+(1,3,\Bbb R)$ are obtained by applying
the homomorphism \mythetag{5.20}. Their components could be taken from the
corresponding frame relationships for the associated frames:
$$
\xalignat 2
&\hskip -2em
\tilde{\boldsymbol\Upsilon}_i=\sum^3_{j=0}S^j_i\,\boldsymbol\Upsilon_j,
&&\boldsymbol\Upsilon_i=\sum^3_{j=0}T^j_i\,\tilde{\boldsymbol\Upsilon}_j.
\mytag{9.17}
\endxalignat
$$
Compare \mythetag{9.16} and \mythetag{9.17} with \mythetag{5.24} and
\mythetag{5.25}. Note that the sign uncertainty in \mythetag{9.16} now 
is absent.\par
    Let's consider the tangent space $T_q(N)$ of the manifold 
\mythetag{9.7} at some of its points $q\in N$. Its complexification
$\Bbb C\otimes T_q(N)$ is the span of the vectors
$$
\align
&\hskip -4em
\bold X_i=\frac{\partial}{\partial x^i},
\mytag{9.18}\\
&\hskip -4em
\bold W^{j_1\ldots\,j_\beta\,\bar j_1\ldots\,\bar j_\gamma\,k_1\ldots
\,k_n}_{i_1\ldots\,i_\alpha\,\bar i_1\ldots\,\bar i_\nu\,h_1\ldots
\,h_m}[P]=\frac{\partial}
{\partial S^{\,i_1\ldots\,i_\alpha\,\bar i_1\ldots\,\bar i_\nu\,
h_1\ldots\,h_m}_{j_1\ldots\,j_\beta\,\bar j_1\ldots\,\bar j_\gamma\,
k_1\ldots\,k_n}[P]},
\mytag{9.19}\\
&\hskip -4em
\bar\bold W^{\bar j_1\ldots\,\bar j_\gamma\,j_1\ldots\,j_\beta\,k_1\ldots
\,k_n}_{\bar i_1\ldots\,\bar i_\nu\,i_1\ldots\,i_\alpha\,h_1\ldots
\,h_m}[P]=\frac{\partial}
{\partial\overline{S^{\,i_1\ldots\,i_\alpha\,\bar i_1\ldots\,\bar i_\nu\,
h_1\ldots\,h_m}_{j_1\ldots\,j_\beta\,\bar j_1\ldots\,\bar j_\gamma\,
k_1\ldots\,k_n}[P]}}.
\mytag{9.20}
\endalign
$$
Note that $\alpha=\alpha_P$, $\beta=\beta_P$, $\nu=\nu_P$,
$\gamma=\gamma_P$, $m=m_P$, $n=n_P$, and $P$ runs from $1$ to $J+Q$ 
in \mythetag{9.19} and in \mythetag{9.20}.
Under the change of local coordinates given by the formulas
\mythetag{9.13}, \mythetag{9.14}, and \mythetag{9.15} the 
tangent vectors \mythetag{9.18}, \mythetag{9.19}, and \mythetag{9.20} 
are transformed according to the next five formulas:
$$
\gather
\tilde\bold X_j=\sum^3_{i=0}\frac{\partial x^i}{\partial\tx^j}\,
\bold X_i+\sum^{J+Q}_{P=1}
\dsize\msum{2}\Sb i_1,\,\ldots,\,i_\alpha\\ j_1,\,\ldots,\,j_\beta\\
a_1,\,\ldots,\,a_\alpha\\ b_1,\,\ldots,\,b_\beta\endSb
\dsize\msum{2}\Sb \bar i_1,\,\ldots,\,\bar i_\nu\\ 
\bar j_1,\,\ldots,\,\bar j_\gamma\\
\bar a_1,\,\ldots,\,\bar a_\nu\\ 
\bar b_1,\,\ldots,\,\bar b_\gamma\endSb
\dsize\msum{3}
\Sb h_1,\,\ldots,\,h_m\\ k_1,\,\ldots,\,k_n\\
c_1,\,\ldots,\,c_m\\ d_1,\,\ldots,\,d_n\endSb
\goth S^{\,i_1}_{a_1}\ldots\,\goth S^{\,i_\alpha}_{a_\alpha}\,\times\\
\vspace{0pt plus 1pt minus 1pt}
\times\,\goth T^{b_1}_{j_1}\ldots\,\goth T^{b_\beta}_{j_\beta}
\ \overline{\goth S^{\,\bar i_1}_{\bar a_1}}
\ldots\,\overline{\goth S^{\,\bar i_\nu}_{\bar a_\nu}}\ 
\ \overline{\goth T^{\,\bar b_1}_{\bar j_1}}
\ldots\,\overline{\goth T^{\,\bar b_\gamma}_{\bar j_\gamma}}\ 
S^{h_1}_{c_1}\ldots\,S^{h_m}_{c_m}\,
T^{\,d_1}_{k_1}\ldots\,T^{\,d_n}_{k_n}\,\times\\
\vspace{0pt plus 1pt minus 1pt}
\times\!\left(\,\sum^\alpha_{\mu=1}\sum^2_{v_\mu=1}\!
\left(\,\shave{\sum^2_{a=1}}
\goth T^{a_\mu}_a\,\frac{\partial\goth S^a_{v_\mu}}{\partial\tx^j}\!\right)
\tilde S^{\,a_1\ldots\,v_\mu\ldots\,a_\alpha\,\bar a_1\ldots\,\bar a_\nu\,
c_1\ldots\,c_m}_{\,b_1\ldots\,\ldots\,\ldots\,b_\beta\,\bar b_1\ldots\,
\bar b_\gamma\,d_1\ldots\,d_n}[P]\,+\right.\\
\vspace{0pt plus 1pt minus 1pt}
+\sum^\beta_{\mu=1}\sum^2_{w_\mu=1}\!
\left(\,\shave{\sum^2_{a=1}}\frac{\partial\goth T^{w_\mu}_a}
{\partial\tx^j}\,\goth S^a_{b_\mu}\!\right)
\tilde S^{\,a_1\ldots\,\ldots\,\ldots\,a_\alpha\,\bar a_1\ldots\,
\bar a_\nu\,c_1\ldots\,c_m}_{\,b_1\ldots\,w_\mu\ldots\,b_\beta\,\bar b_1
\ldots\,\bar b_\gamma\,d_1\ldots\,d_n}[P]\,+\\
\vspace{0pt plus 1pt minus 1pt}
+\sum^\nu_{\mu=1}\sum^2_{v_\mu=1}\!\left(\,\shave{\sum^2_{\bar a=1}}
\overline{\goth T^{\bar a_\mu}_{\bar a}}\,
\frac{\partial\overline{\goth S^{\bar a}_{v_\mu}}}{\partial\tx^j}\!\right)
\tilde S^{\,a_1\ldots\,a_\alpha\,\bar a_1\ldots\,v_\mu\ldots\,\bar a_\nu\,
c_1\ldots\,c_m}_{\,b_1\ldots\,b_\beta\,\bar b_1\ldots\,\ldots\,\ldots\,
\bar b_\gamma\,d_1\ldots\,d_n}[P]\,+\\
\vspace{0pt plus 1pt minus 1pt}
+\sum^\gamma_{\mu=1}\sum^2_{w_\mu=1}\!
\left(\,\shave{\sum^2_{\bar a=1}}\frac{\partial\overline{\goth
T^{w_\mu}_{\bar a}}}{\partial\tx^j}\,\overline{\goth S^{\bar a}_{\bar
b_\mu}}\!\right) \tilde S^{\,a_1\ldots\,a_\alpha\,\bar a_1\ldots\,
\ldots\,\ldots\,\bar a_\nu\,c_1\ldots\,c_m}_{\,b_1\ldots\,b_\beta\,
\bar b_1\ldots\,w_\mu\ldots\,\bar b_\gamma\,d_1\ldots\,d_n}[P]\,+\\
\vspace{0pt plus 1pt minus 1pt}
+\sum^m_{\mu=1}\sum^3_{v_\mu=0}\!\left(\,\shave{\sum^3_{a=0}}
T^{\,c_\mu}_a\,\frac{\partial S^{\,a}_{v_m}}{\partial\tx^j}\!\right)
\tilde S^{\,a_1\ldots\,a_\alpha\,\bar a_1\ldots\,
\bar a_\nu\,c_1\ldots\,v_\mu\ldots
\,c_m}_{\,b_1\ldots\,b_\beta\,\bar b_1\ldots\,\bar b_\gamma\,d_1\ldots\,
\ldots\,\ldots\,d_n}[P]\,+\\
\vspace{0pt plus 1pt minus 1pt}
\left.+\sum^n_{\mu=1}\sum^3_{w_\mu=0}\!
\left(\,\shave{\sum^3_{a=0}}\frac{\partial T^{w_\mu}_a}
{\partial\tx^j}\,S^{\,a}_{d_\mu}\!\right)
\tilde S^{\,a_1\ldots\,a_\alpha\,\bar a_1\ldots\,
\bar a_\nu\,c_1\ldots\,\ldots\,\ldots\,c_m}_{\,b_1\ldots\,b_\beta\,
\bar b_1\ldots\,\bar b_\gamma\,d_1\ldots\,w_\mu\ldots
\,d_n}[P]\right)\times\\
\vspace{0pt plus 1pt minus 1pt}
\times\,\bold W^{j_1\ldots\,j_\beta\,\bar j_1\ldots\,\bar j_\gamma\,k_1\ldots
\,k_n}_{i_1\ldots\,i_\alpha\,\bar i_1\ldots\,\bar i_\nu\,h_1\ldots
\,h_m}[P]
+\sum^{J+Q}_{P=1}
\dsize\msum{2}\Sb i_1,\,\ldots,\,i_\alpha\\ j_1,\,\ldots,\,j_\beta\\
a_1,\,\ldots,\,a_\alpha\\ b_1,\,\ldots,\,b_\beta\endSb
\dsize\msum{2}\Sb \bar i_1,\,\ldots,\,\bar i_\nu\\ 
\bar j_1,\,\ldots,\,\bar j_\gamma\\
\bar a_1,\,\ldots,\,\bar a_\nu\\ 
\bar b_1,\,\ldots,\,\bar b_\gamma\endSb
\dsize\msum{3}
\Sb h_1,\,\ldots,\,h_m\\ k_1,\,\ldots,\,k_n\\
c_1,\,\ldots,\,c_m\\ d_1,\,\ldots,\,d_n\endSb
\overline{\goth S^{\,i_1}_{a_1}}\ldots\,
\overline{\goth S^{\,i_\alpha}_{a_\alpha}}\,\times\\
\vspace{0pt plus 1pt minus 1pt}
\times\,\overline{\goth T^{b_1}_{j_1}}\ldots\,
\overline{\goth T^{b_\beta}_{j_\beta}}
\ \goth S^{\,\bar i_1}_{\bar a_1}
\ldots\,\goth S^{\,\bar i_\nu}_{\bar a_\nu}\ 
\ \goth T^{\,\bar b_1}_{\bar j_1}
\ldots\,\goth T^{\,\bar b_\gamma}_{\bar j_\gamma}\ 
S^{h_1}_{c_1}\ldots\,S^{h_m}_{c_m}\,
T^{\,d_1}_{k_1}\ldots\,T^{\,d_n}_{k_n}\,\times\\
\vspace{0pt plus 1pt minus 1pt}
\times\!\left(\,\sum^\alpha_{\mu=1}\sum^2_{v_\mu=1}\!
\left(\,\shave{\sum^2_{a=1}}
\overline{\goth T^{a_\mu}_a}\,\frac{\partial
\overline{\goth S^a_{v_\mu}}}{\partial\tx^j}\!\right)
\overline{\tilde S^{\,a_1\ldots\,v_\mu\ldots\,a_\alpha\,
\bar a_1\ldots\,\bar a_\nu\,c_1\ldots\,c_m}_{\,b_1\ldots\,
\ldots\,\ldots\,b_\beta\,\bar b_1\ldots\,\bar b_\gamma
\,d_1\ldots\,d_n}[P]}\,+\right.\\
\vspace{0pt plus 1pt minus 1pt}
+\sum^\beta_{\mu=1}\sum^2_{w_\mu=1}\!
\left(\,\shave{\sum^2_{a=1}}\frac{\partial
\overline{\goth T^{w_\mu}_a}}{\partial\tx^j}\,
\overline{\goth S^a_{b_\mu}}\!\right)
\overline{\tilde S^{\,a_1\ldots\,\ldots\,\ldots\,a_\alpha\,\bar a_1\ldots\,
\bar a_\nu\,c_1\ldots\,c_m}_{\,b_1\ldots\,w_\mu\ldots\,b_\beta\,
\bar b_1\ldots\,\bar b_\gamma\,d_1\ldots\,d_n}[P]}\,+\\
\vspace{0pt plus 1pt minus 1pt}
+\sum^\nu_{\mu=1}\sum^2_{v_\mu=1}\!\left(\,\shave{\sum^2_{\bar a=1}}
\goth T^{\bar a_\mu}_{\bar a}\,
\frac{\partial\goth S^{\bar a}_{v_\mu}}{\partial\tx^j}\!\right)
\overline{\tilde S^{\,a_1\ldots\,a_\alpha\,\bar a_1\ldots\,v_\mu\ldots
\,\bar a_\nu\,
c_1\ldots\,c_m}_{\,b_1\ldots\,b_\beta\,\bar b_1\ldots\,\ldots\,\ldots\,
\bar b_\gamma\,d_1\ldots\,d_n}[P]}\,+\\
\vspace{0pt plus 1pt minus 1pt}
+\sum^\gamma_{\mu=1}\sum^2_{w_\mu=1}\!
\left(\,\shave{\sum^2_{\bar a=1}}\frac{\partial\goth
T^{w_\mu}_{\bar a}}{\partial\tx^j}\,\goth S^{\bar a}_{\bar
b_\mu}\!\right) \overline{\tilde S^{\,a_1\ldots\,a_\alpha\,\bar a_1\ldots\,
\ldots\,\ldots\,\bar a_\nu\,c_1\ldots\,c_m}_{\,b_1\ldots\,b_\beta\,
\bar b_1\ldots\,w_\mu\ldots\,\bar b_\gamma\,d_1\ldots\,d_n}[P]}\,+\\
\vspace{0pt plus 1pt minus 1pt}
+\sum^m_{\mu=1}\sum^3_{v_\mu=0}\!\left(\,\shave{\sum^3_{a=0}}
T^{\,c_\mu}_a\,\frac{\partial S^{\,a}_{v_m}}{\partial\tx^j}\!\right)
\overline{\tilde S^{\,a_1\ldots\,a_\alpha\,\bar a_1\ldots\,
\bar a_\nu\,c_1\ldots\,v_\mu\ldots
\,c_m}_{\,b_1\ldots\,b_\beta\,\bar b_1\ldots\,\bar b_\gamma\,d_1\ldots\,
\ldots\,\ldots\,d_n}[P]}\,+\\
\displaybreak
\left.+\sum^n_{\mu=1}\sum^3_{w_\mu=0}\!
\left(\,\shave{\sum^3_{a=0}}\frac{\partial T^{w_\mu}_a}
{\partial\tx^j}\,S^{\,a}_{d_\mu}\!\right)
\overline{\tilde S^{\,a_1\ldots\,a_\alpha\,\bar a_1\ldots\,
\bar a_\nu\,c_1\ldots\,\ldots\,\ldots\,c_m}_{\,b_1\ldots\,b_\beta\,
\bar b_1\ldots\,\bar b_\gamma\,d_1\ldots\,w_\mu\ldots
\,d_n}[P]}\right)\times\\
\times\,\bar\bold W^{\bar j_1\ldots\,\bar j_\gamma\,j_1\ldots\,j_\beta\,
k_1\ldots\,k_n}_{\bar i_1\ldots\,\bar i_\nu\,i_1\ldots\,i_\alpha\,
h_1\ldots\,h_m}[P].
\endgather
$$
As we see, the transformation rule for the vector \mythetag{9.18} is 
expressed by the formula bigger than the page size. For this reason 
we give it without its inverse counterpart.
$$
\align
&\hskip -4em
\left\{
\aligned
&\tilde\bold W^{j_1\ldots\,j_\beta\,\bar j_1\ldots\,\bar j_\gamma\,k_1\ldots
\,k_n}_{i_1\ldots\,i_\alpha\,\bar i_1\ldots\,\bar i_\nu\,h_1\ldots
\,h_m}[P]
=\dsize\msum{2}\Sb a_1,\,\ldots,\,a_\alpha\\ b_1,\,\ldots,\,b_\beta\endSb
\dsize\msum{2}\Sb \bar a_1,\,\ldots,\,\bar a_\nu\\ 
\bar b_1,\,\ldots,\,\bar b_\gamma\endSb
\dsize\msum{3}
\Sb c_1,\,\ldots,\,c_m\\ d_1,\,\ldots,\,d_n\endSb
\goth S^{a_1}_{i_1}\ldots\,\goth S^{a_\alpha}_{i_\alpha}\,\times\\
&\kern 40pt
\times\,\goth T^{j_1}_{b_1}\ldots\,\goth T^{j_\beta}_{b_\beta}
\ \overline{\goth S^{\,\bar a_1}_{\bar i_1}}
\ldots\,\overline{\goth S^{\,\bar a_\nu}_{\bar i_\nu}}\ 
\ \overline{\goth T^{\,\bar j_1}_{\bar b_1}}
\ldots\,\overline{\goth T^{\,\bar j_\gamma}_{\bar b_\gamma}}\ 
S^{c_1}_{h_1}\ldots\,S^{c_m}_{h_m}\,\times\\
\vspace{2ex}
&\kern 60pt
\times\,T^{\,k_1}_{d_1}\ldots\,T^{\,k_n}_{d_n}\ 
\bold W^{\,b_1\ldots\,b_\beta\,\bar b_1\ldots\,\bar b_\gamma\,
d_1\ldots\,d_n}_{\,a_1\ldots\,a_\alpha\,\bar a_1\ldots\,\bar a_\nu
\,c_1\ldots\,c_m}[P],
\endaligned
\right.
\mytag{9.21}\\
&\hskip -4em
\left\{
\aligned
&\bold W^{j_1\ldots\,j_\beta\,\bar j_1\ldots\,\bar j_\gamma\,k_1\ldots
\,k_n}_{i_1\ldots\,i_\alpha\,\bar i_1\ldots\,\bar i_\nu\,h_1\ldots
\,h_m}[P]
=\dsize\msum{2}\Sb a_1,\,\ldots,\,a_\alpha\\ b_1,\,\ldots,\,b_\beta\endSb
\dsize\msum{2}\Sb \bar a_1,\,\ldots,\,\bar a_\nu\\ 
\bar b_1,\,\ldots,\,\bar b_\gamma\endSb
\dsize\msum{3}
\Sb c_1,\,\ldots,\,c_m\\ d_1,\,\ldots,\,d_n\endSb
\goth T^{a_1}_{i_1}\ldots\,\goth T^{a_\alpha}_{i_\alpha}\,\times\\
&\kern 40pt
\times\,\goth S^{j_1}_{b_1}\ldots\,\goth S^{j_\beta}_{b_\beta}
\ \overline{\goth T^{\,\bar a_1}_{\bar i_1}}
\ldots\,\overline{\goth T^{\,\bar a_\nu}_{\bar i_\nu}}\ 
\ \overline{\goth S^{\,\bar j_1}_{\bar b_1}}
\ldots\,\overline{\goth S^{\,\bar j_\gamma}_{\bar b_\gamma}}\ 
T^{c_1}_{h_1}\ldots\,T^{c_m}_{h_m}\,\times\\
\vspace{2ex}
&\kern 60pt
\times\,S^{\,k_1}_{d_1}\ldots\,S^{\,k_n}_{d_n}\ 
\tilde\bold W^{\,b_1\ldots\,b_\beta\,\bar b_1\ldots\,\bar b_\gamma\,
d_1\ldots\,d_n}_{\,a_1\ldots\,a_\alpha\,\bar a_1\ldots\,\bar a_\nu
\,c_1\ldots\,c_m}[P],
\endaligned
\right.
\mytag{9.22}\\
\vspace{2ex}
&\hskip -4em
\left\{
\aligned
&\Tilde{\Bar{\bold W}}\vphantom{\bold W}^{\bar j_1\ldots\,\bar j_\gamma\,
j_1\ldots\,j_\beta\,k_1\ldots\,k_n}_{\bar i_1\ldots\,\bar i_\nu\,
i_1\ldots\,i_\alpha\,h_1\ldots\,h_m}[P]
=\dsize\msum{2}\Sb a_1,\,\ldots,\,a_\alpha\\ b_1,\,\ldots,\,b_\beta\endSb
\dsize\msum{2}\Sb \bar a_1,\,\ldots,\,\bar a_\nu\\ 
\bar b_1,\,\ldots,\,\bar b_\gamma\endSb
\dsize\msum{3}
\Sb c_1,\,\ldots,\,c_m\\ d_1,\,\ldots,\,d_n\endSb
\overline{\goth S^{a_1}_{i_1}}\ldots\,
\overline{\goth S^{a_\alpha}_{i_\alpha}}\,\times\\
&\kern 40pt
\times\,\overline{\goth T^{j_1}_{b_1}}\ldots\,
\overline{\goth T^{j_\beta}_{b_\beta}}
\ \goth S^{\,\bar a_1}_{\bar i_1}
\ldots\,\goth S^{\,\bar a_\nu}_{\bar i_\nu}\ 
\ \goth T^{\,\bar j_1}_{\bar b_1}
\ldots\,\goth T^{\,\bar j_\gamma}_{\bar b_\gamma}\ 
S^{c_1}_{h_1}\ldots\,S^{c_m}_{h_m}\,\times\\
\vspace{2ex}
&\kern 60pt
\times\,T^{\,k_1}_{d_1}\ldots\,T^{\,k_n}_{d_n}\ 
\bar\bold W^{\,\bar b_1\ldots\,\bar b_\gamma\,b_1\ldots\,b_\beta\,
d_1\ldots\,d_n}_{\,\bar a_1\ldots\,\bar a_\nu\,a_1\ldots\,a_\alpha\,
c_1\ldots\,c_m}[P],
\endaligned
\right.
\mytag{9.23}\\
&\hskip -4em
\left\{
\aligned
&\bar\bold W^{\bar j_1\ldots\,\bar j_\gamma\,j_1\ldots\,j_\beta\,k_1
\ldots\,k_n}_{\bar i_1\ldots\,\bar i_\nu\,i_1\ldots\,i_\alpha\,h_1
\ldots\,h_m}[P]
=\dsize\msum{2}\Sb a_1,\,\ldots,\,a_\alpha\\ b_1,\,\ldots,\,b_\beta\endSb
\dsize\msum{2}\Sb \bar a_1,\,\ldots,\,\bar a_\nu\\ 
\bar b_1,\,\ldots,\,\bar b_\gamma\endSb
\dsize\msum{3}
\Sb c_1,\,\ldots,\,c_m\\ d_1,\,\ldots,\,d_n\endSb
\overline{\goth T^{a_1}_{i_1}}\ldots\,
\overline{\goth T^{a_\alpha}_{i_\alpha}}\,\times\\
&\kern 40pt
\times\,\overline{\goth S^{j_1}_{b_1}}\ldots\,
\overline{\goth S^{j_\beta}_{b_\beta}}
\ \goth T^{\,\bar a_1}_{\bar i_1}
\ldots\,\goth T^{\,\bar a_\nu}_{\bar i_\nu}\ 
\ \goth S^{\,\bar j_1}_{\bar b_1}
\ldots\,\goth S^{\,\bar j_\gamma}_{\bar b_\gamma}\ 
T^{c_1}_{h_1}\ldots\,T^{c_m}_{h_m}\,\times\\
\vspace{2ex}
&\kern 60pt
\times\,S^{\,k_1}_{d_1}\ldots\,S^{\,k_n}_{d_n}\ 
\Tilde{\Bar{\bold W}}\vphantom{\bold W}^{\,\bar b_1\ldots\,\bar b_\gamma\,
b_1\ldots\,b_\beta\,d_1\ldots\,d_n}_{\,\bar a_1\ldots\,\bar a_\nu\,
a_1\ldots\,a_\alpha\,c_1\ldots\,c_m}[P].
\endaligned
\right.
\mytag{9.24}
\endalign
$$
The transformation rules \mythetag{9.21}, \mythetag{9.22}, 
\mythetag{9.23}, and \mythetag{9.24} are subdivided into 
pairs: each direct transformation rule is paired with an 
inverse one. Note also that the coefficients in \mythetag{9.21} 
and \mythetag{9.22} differ from those in \mythetag{9.23} and 
\mythetag{9.24} only by complex conjugation. Unlike the corresponding
formula in \mycite{5}, the transition matrices $S$ and $T$ here are 
not Jacobi matrices. In particular, we have
$$
\hskip -2em
\frac{\partial x^i}{\partial\tx^j}\neq S^i_j.
\mytag{9.25}
$$
In order to replace the partial derivatives \mythetag{9.25} by 
the components of the transition matrix $S$, we should pass from 
the holonomic frame given by the vectors \mythetag{9.18}, 
\mythetag{9.19}, and \mythetag{9.20} to some specially designed non-holonomic 
frame in $N$. For this purpose let's keep the vectors
\mythetag{9.19} and \mythetag{9.20} unchanged, but replace the tangent
vectors \mythetag{9.18} by the following ones:
$$
\hskip -2em
\bold U_i=\sum^3_{j=0}\Upsilon^j_i\,\bold X_j.
\mytag{9.26}
$$
The coefficients $\Upsilon^j_i$ in the formula \mythetag{9.26} are the
coordinates of the frame vectors $\boldsymbol\Upsilon_0,\,
\boldsymbol\Upsilon_1,\,\boldsymbol\Upsilon_2,\,\boldsymbol\Upsilon_3$,
i\.\,e\. they are taken from the expansion
$$
\hskip -2em
\boldsymbol\Upsilon_i=\sum^3_{j=0}\Upsilon^j_i\,\bold E_j,
\mytag{9.27}
$$
where $\bold E_0,\,\bold E_1,\,\bold E_2,\,\bold E_3$ are the coordinate
vectors \mythetag{5.3} forming a holonomic frame in $M$. Despite to the
similarity of the formulas \mythetag{9.26} and \mythetag{9.27} they
are different. The formula \mythetag{9.27} is the expansion in the
tangent space to the space-time manifold $M$, while the formula
\mythetag{9.26} is the expansion in the tangent space to the composite
spin-tensorial bundle \mythetag{9.7}.\par
     Upon passing from \mythetag{9.18} to the new vectors $\bold U_0,\,
\bold U_1,\,\bold U_2,\,\bold U_3$ in \mythetag{9.26} we introduce the
following $\theta$-parameters defined through the transition matrices:
$$
\align
&\hskip -6em
\tilde\theta^k_{ij}=\sum^3_{a=0}T^k_a\,L_{\tilde{\boldsymbol\Upsilon}_i}
\!(S^a_j)=\sum^3_{a=0}\sum^3_{v=0}T^k_a\,\tilde\Upsilon^v_i\,\frac{\partial
S^a_j}{\partial\tx^v}=-\sum^3_{a=0}L_{\tilde{\boldsymbol\Upsilon}_i}
\!(T^k_a)\,S^a_j,\hskip -2em
\mytag{9.28}\\
&\hskip -6em
\tilde\vartheta^k_{ij}=\sum^2_{a=1}\goth T^k_a\,
L_{\tilde{\boldsymbol\Upsilon}_i}(\goth S^a_j)
=\sum^2_{a=1}\sum^3_{v=0}\goth T^k_a\,\tilde\Upsilon^v_i\,
\frac{\partial\goth S^a_j}{\partial\tx^v}=
-\sum^2_{a=1}L_{\tilde{\boldsymbol\Upsilon}_i}
\!(\goth T^k_a)\,\goth S^a_j.\hskip -2em
\mytag{9.29}
\endalign
$$
The $\theta$-parameters without tilde are introduced in a similar way:
$$
\align
&\hskip -6em
\theta^k_{ij}=\sum^3_{a=0}S^k_a\,L_{\boldsymbol\Upsilon_i}
\!(T^a_j)=\sum^3_{a=0}\sum^3_{v=0}S^k_a\,\Upsilon^v_i\,\frac{\partial
T^a_j}{\partial x^v}=-\sum^3_{a=0}L_{\boldsymbol\Upsilon_i}
\!(S^k_a)\,T^a_j,\hskip -2em
\mytag{9.30}\\
&\hskip -6em
\vartheta^k_{ij}=\sum^2_{a=1}\goth S^k_a\,
L_{\boldsymbol\Upsilon_i}(\goth T^a_j)
=\sum^2_{a=1}\sum^3_{v=0}\goth S^k_a\,\Upsilon^v_i\,
\frac{\partial\goth T^a_j}{\partial\tx^v}=
-\sum^2_{a=1}L_{\boldsymbol\Upsilon_i}
\!(\goth S^k_a)\,\goth T^a_j.\hskip -2em
\mytag{9.31}
\endalign
$$
Here $L_{\boldsymbol\Upsilon_i}$ and $L_{\tilde{\boldsymbol\Upsilon}_i}$
are Lie derivatives (see \mycite{11}). Using \mythetag{9.28}, 
\mythetag{9.29}, \mythetag{9.30}, \mythetag{9.31} one easily derives the 
following identities similar to those in \mycite{5}:
$$
\align
&\hskip -2em
\theta^k_{ij}=-\sum^3_{a=0}\sum^3_{b=0}\sum^3_{c=0}
T^a_i\,\tilde\theta^{\,c}_{\!ab}\ S^k_c\,T^b_j,
\mytag{9.32}\\
&\hskip -2em
\tilde\theta^k_{ij}=-\sum^3_{a=0}\sum^3_{b=0}\sum^3_{c=0}
S^a_i\,\theta^{\,c}_{\!ab}\ T^k_c\,S^b_j,
\mytag{9.33}\\
&\hskip -2em
\vartheta^k_{ij}=-\sum^3_{a=0}\sum^2_{b=1}\sum^2_{c=1}
T^a_i\,\tilde\vartheta^{\,c}_{\!ab}\ \goth S^k_c\,\goth T^b_j,
\mytag{9.34}\\
&\hskip -2em
\tilde\vartheta^k_{ij}=-\sum^3_{a=0}\sum^2_{b=1}\sum^2_{c=1}
S^a_i\,\vartheta^{\,c}_{\!ab}\ \goth T^k_c\,\goth S^b_j.
\mytag{9.35}\\
\endalign
$$
The formulas \mythetag{9.32}, \mythetag{9.33}, \mythetag{9.34},
and \mythetag{9.35} relate the $\theta$-parameters with and without
tilde. Note that, unlike those introduced in \mycite{5}, the 
$\theta$-parameters introduced in \mythetag{9.28} and \mythetag{9.30}
are not symmetric with respect to $i$ and $j$. The extent of their
asymmetry is given by the following formulas:
$$
\xalignat 2
&\hskip -2em
\theta^k_{ij}-\theta^k_{j\,i}=c^k_{ij}
&&\tilde\theta^k_{ij}-\tilde\theta^k_{j\,i}=\tilde c^k_{ij}
\mytag{9.36}
\endxalignat
$$
The parameters $c^k_{ij}$ in \mythetag{9.36} are taken from
the following commutator relationships for the frame vectors
of the orthonormal frame:
$$
[\boldsymbol\Upsilon_i,\,\boldsymbol\Upsilon_j]=
\sum^3_{k=1}c^{\,k}_{ij}\,\boldsymbol\Upsilon_k.
$$
The parameters $\tilde c^k_{ij}$ are taken from the analogous formula for
$\tilde{\boldsymbol\Upsilon}_0,\,\tilde{\boldsymbol\Upsilon}_1,\,
\tilde{\boldsymbol\Upsilon}_2,\,\tilde{\boldsymbol\Upsilon}_3$.\par
    Now we apply \mythetag{9.26}, \mythetag{9.28}, \mythetag{9.29},
\mythetag{9.30}, \mythetag{9.31} to the above huge transformation formula.
Then it is written so that we can assign a number to this formula: 
$$
\left\{
\gathered
\tilde\bold U_j=\sum^3_{i=0}S^{\,i}_j\,\bold U_i+\\
+\sum^{J+Q}_{P=1}\dsize\msums{2}{3}
\Sb i_1,\,\ldots,\,i_\alpha\\ j_1,\,\ldots,\,j_\beta\\
\bar i_1,\,\ldots,\,\bar i_\nu\\ \bar j_1,\,\ldots,\,\bar j_\gamma\\
h_1,\,\ldots,\,h_m\\ k_1,\,\ldots,\,k_n\endSb
\left(\tilde w^{i_1\ldots\,i_\alpha\,\bar i_1\ldots\,\bar i_\nu\,h_1
\ldots\,h_m}_{j\,j_1\ldots\,j_\beta\,\bar j_1\ldots\,\bar j_\gamma\,k_1
\ldots\,k_n}[P]\ \bold W^{j_1\ldots\,j_\beta\,\bar j_1\ldots\,\bar j_\gamma
\,k_1\ldots\,k_n}_{i_1\ldots\,i_\alpha\,\bar i_1\ldots\,\bar i_\nu
\,h_1\ldots\,h_m}[P]\,+\right.\\
\vspace{3ex}
\left.\kern 70pt+\,\overline{\tilde w^{i_1\ldots\,i_\alpha\,
\bar i_1\ldots\,\bar i_\nu\,h_1\ldots\,h_m}_{j\,j_1\ldots\,j_\beta\,
\bar j_1\ldots\,\bar j_\gamma\,k_1\ldots\,k_n}[P]}
\ \bar\bold W^{\bar j_1\ldots\,\bar j_\gamma\,j_1\ldots\,j_\beta
\,k_1\ldots\,k_n}_{\bar i_1\ldots\,\bar i_\nu\,i_1\ldots\,i_\alpha
\,h_1\ldots\,h_m}[P]\right).
\endgathered\right.
\quad
\mytag{9.37}
$$
The inverse transformation formula for \mythetag{9.37} is written
similarly:
$$
\left\{
\gathered
\bold U_j=\sum^3_{i=0}S^{\,i}_j\,\tilde\bold U_i\,+\\
+\sum^{J+Q}_{P=1}\dsize\msums{2}{3}
\Sb i_1,\,\ldots,\,i_\alpha\\ j_1,\,\ldots,\,j_\beta\\
\bar i_1,\,\ldots,\,\bar i_\nu\\ \bar j_1,\,\ldots,\,\bar j_\gamma\\
h_1,\,\ldots,\,h_m\\ k_1,\,\ldots,\,k_n\endSb
\left(w^{i_1\ldots\,i_\alpha\,\bar i_1\ldots\,\bar i_\nu\,h_1
\ldots\,h_m}_{j\,j_1\ldots\,j_\beta\,\bar j_1\ldots\,\bar j_\gamma\,k_1
\ldots\,k_n}[P]\ \tilde \bold W^{j_1\ldots\,j_\beta\,\bar j_1\ldots\,
\bar j_\gamma\,k_1\ldots\,k_n}_{i_1\ldots\,i_\alpha\,\bar i_1\ldots\,
\bar i_\nu\,h_1\ldots\,h_m}[P]\,+\right.\\
\vspace{3ex}
\left.\kern 70pt+\,\overline{w^{i_1\ldots\,i_\alpha\,
\bar i_1\ldots\,\bar i_\nu\,h_1\ldots\,h_m}_{j\,j_1\ldots\,j_\beta\,
\bar j_1\ldots\,\bar j_\gamma\,k_1\ldots\,k_n}[P]}
\ \Tilde{\Bar{\bold W}}\vphantom{\bold W}^{\bar j_1\ldots\,
\bar j_\gamma\,j_1\ldots\,j_\beta\,k_1\ldots\,k_n}_{\bar i_1\ldots\,
\bar i_\nu\,i_1\ldots\,i_\alpha\,h_1\ldots\,h_m}[P]\right).
\endgathered\right.
\quad
\mytag{9.38}
$$
Here are the formulas for \pagebreak $\tilde w^{i_1\ldots\,i_\alpha\,
\bar i_1\ldots\,\bar i_\nu\,h_1\ldots\,h_m}_{j\,j_1\ldots\,j_\beta\,
\bar j_1\ldots\,\bar j_\gamma\,k_1\ldots\,k_n}[P]$ and $w^{i_1\ldots\,
i_\alpha\,\bar i_1\ldots\,\bar i_\nu\,h_1\ldots\,h_m}_{j\,j_1\ldots
\,j_\beta\,\bar j_1\ldots\,\bar j_\gamma\,k_1\ldots\,k_n}[P]$
in the above transformation formulas \mythetag{9.37} and \mythetag{9.38}:
$$
\allowdisplaybreaks
\gather
\tilde w^{i_1\ldots\,i_\alpha\,\bar i_1\ldots\,\bar i_\nu\,h_1
\ldots\,h_m}_{j\,j_1\ldots\,j_\beta\,\bar j_1\ldots\,\bar j_\gamma
\,k_1\ldots\,k_n}[P]=
\dsize\msum{2}\Sb a_1,\,\ldots,\,a_\alpha\\ b_1,\,\ldots,\,b_\beta\endSb
\dsize\msum{2}\Sb \bar a_1,\,\ldots,\,\bar a_\nu\\ 
\bar b_1,\,\ldots,\,\bar b_\gamma\endSb
\dsize\msum{3}
\Sb c_1,\,\ldots,\,c_m\\ d_1,\,\ldots,\,d_n\endSb
\goth S^{\,i_1}_{a_1}\ldots\,\goth S^{\,i_\alpha}_{a_\alpha}\,\times\\
\vspace{0pt plus 1pt minus 1pt}
\times\,\goth T^{b_1}_{j_1}\ldots\,\goth T^{b_\beta}_{j_\beta}
\ \overline{\goth S^{\,\bar i_1}_{\bar a_1}}
\ldots\,\overline{\goth S^{\,\bar i_\nu}_{\bar a_\nu}}\ 
\ \overline{\goth T^{\,\bar b_1}_{\bar j_1}}
\ldots\,\overline{\goth T^{\,\bar b_\gamma}_{\bar j_\gamma}}\ 
S^{h_1}_{c_1}\ldots\,S^{h_m}_{c_m}\,
T^{\,d_1}_{k_1}\ldots\,T^{\,d_n}_{k_n}\,\times\\
\times
\left(\shave{\,\sum^\alpha_{\mu=1}\sum^2_{v_\mu=1}}
\tilde\vartheta^{\,a_\mu}_{\!j\,v_\mu}\
\tilde S^{\,a_1\ldots\,v_\mu\ldots\,a_\alpha\,\bar a_1\ldots\,\bar a_\nu\,
c_1\ldots\,c_m}_{\,b_1\ldots\,\ldots\,\ldots\,b_\beta\,\bar b_1\ldots\,
\bar b_\gamma\,d_1\ldots\,d_n}[P]-\shave{\,\sum^\beta_{\mu=1}\sum^2_{w_\mu=1}}
\tilde\vartheta^{\,w_\mu}_{\!j\,b_\mu}\,\times\right.\\
\times\,\tilde S^{\,a_1\ldots\,\ldots\,\ldots\,a_\alpha\,\bar a_1\ldots\,
\bar a_\nu\,c_1\ldots\,c_m}_{\,b_1\ldots\,w_\mu\ldots\,b_\beta\,\bar b_1\ldots
\,\bar b_\gamma\,d_1\ldots\,d_n}[P]+
\shave{\,\sum^\nu_{\mu=1}\sum^2_{v_\mu=1}}
\overline{\tilde\vartheta^{\,\bar a_\mu}_{\!j\,v_\mu}}\
\tilde S^{\,a_1\ldots\,a_\alpha\,\bar a_1\ldots\,
\ldots\,\ldots\,\bar a_\nu\,c_1\ldots\,c_m}_{\,b_1\ldots\,b_\beta\,
\bar b_1\ldots\,w_\mu\ldots\,\bar b_\gamma\,d_1\ldots\,d_n}[P]\,-\\
-\shave{\,\sum^\gamma_{\mu=1}\sum^2_{w_\mu=1}}
\overline{\tilde\vartheta^{\,w_\mu}_{\!j\,\bar b_\mu}}\ 
\tilde S^{\,a_1\ldots\,a_\alpha\,\bar a_1\ldots\,
\ldots\,\ldots\,\bar a_\nu\,c_1\ldots\,c_m}_{\,b_1\ldots\,b_\beta\,
\bar b_1\ldots\,w_\mu\ldots\,\bar b_\gamma\,d_1\ldots\,d_n}[P]
+\shave{\,\sum^m_{\mu=1}\sum^3_{v_\mu=0}}
\tilde\theta^{\,c_\mu}_{\!j\,v_\mu}\,\times\\
\left.\times\,\tilde S^{\,a_1\ldots\,a_\alpha\,\bar a_1\ldots\,
\bar a_\nu\,c_1\ldots\,v_\mu\ldots
\,c_m}_{\,b_1\ldots\,b_\beta\,\bar b_1\ldots\,\bar b_\gamma\,d_1\ldots\,
\ldots\,\ldots\,d_n}[P]
-\shave{\sum^n_{\mu=1}\sum^3_{w_\mu=0}}
\tilde\theta^{\,w_\mu}_{\!j\,d_\mu}
\tilde S^{\,a_1\ldots\,a_\alpha\,\bar a_1\ldots\,
\bar a_\nu\,c_1\ldots\,\ldots\,\ldots\,c_m}_{\,b_1\ldots\,b_\beta\,
\bar b_1\ldots\,\bar b_\gamma\,d_1\ldots\,w_\mu\ldots
\,d_n}[P]\right),\\
\vspace{4ex}
w^{i_1\ldots\,i_\alpha\,\bar i_1\ldots\,\bar i_\nu\,h_1\ldots
\,h_m}_{j\,j_1\ldots\,j_\beta\,\bar j_1\ldots\,\bar j_\gamma\,k_1
\ldots\,k_n}[P]=
\dsize\msum{2}\Sb a_1,\,\ldots,\,a_\alpha\\ b_1,\,\ldots,\,b_\beta\endSb
\dsize\msum{2}\Sb \bar a_1,\,\ldots,\,\bar a_\nu\\ 
\bar b_1,\,\ldots,\,\bar b_\gamma\endSb
\dsize\msum{3}
\Sb c_1,\,\ldots,\,c_m\\ d_1,\,\ldots,\,d_n\endSb
\goth T^{\,i_1}_{a_1}\ldots\,\goth T^{\,i_\alpha}_{a_\alpha}\,\times\\
\vspace{0pt plus 1pt minus 1pt}
\times\,\goth S^{b_1}_{j_1}\ldots\,\goth S^{b_\beta}_{j_\beta}
\ \overline{\goth T^{\,\bar i_1}_{\bar a_1}}
\ldots\,\overline{\goth T^{\,\bar i_\nu}_{\bar a_\nu}}\ 
\ \overline{\goth S^{\,\bar b_1}_{\bar j_1}}
\ldots\,\overline{\goth S^{\,\bar b_\gamma}_{\bar j_\gamma}}\ 
T^{h_1}_{c_1}\ldots\,T^{h_m}_{c_m}\,
S^{\,d_1}_{k_1}\ldots\,S^{\,d_n}_{k_n}\,\times\\
\times
\left(\shave{\,\sum^\alpha_{\mu=1}\sum^2_{v_\mu=1}}
\vartheta^{\,a_\mu}_{\!j\,v_\mu}\
S^{\,a_1\ldots\,v_\mu\ldots\,a_\alpha\,\bar a_1\ldots\,\bar a_\nu\,
c_1\ldots\,c_m}_{\,b_1\ldots\,\ldots\,\ldots\,b_\beta\,\bar b_1\ldots\,
\bar b_\gamma\,d_1\ldots\,d_n}[P]-\shave{\,\sum^\beta_{\mu=1}\sum^2_{w_\mu=1}}
\vartheta^{\,w_\mu}_{\!j\,b_\mu}\,\times\right.\\
\times\,S^{\,a_1\ldots\,\ldots\,\ldots\,a_\alpha\,\bar a_1\ldots\,
\bar a_\nu\,c_1\ldots\,c_m}_{\,b_1\ldots\,w_\mu\ldots\,b_\beta\,\bar b_1
\ldots\,\bar b_\gamma\,d_1\ldots\,d_n}[P]+
\shave{\,\sum^\nu_{\mu=1}\sum^2_{v_\mu=1}}
\overline{\vartheta^{\,\bar a_\mu}_{\!j\,v_\mu}}\
S^{\,a_1\ldots\,a_\alpha\,\bar a_1\ldots\,
\ldots\,\ldots\,\bar a_\nu\,c_1\ldots\,c_m}_{\,b_1\ldots\,b_\beta\,
\bar b_1\ldots\,w_\mu\ldots\,\bar b_\gamma\,d_1\ldots\,d_n}[P]\,-\\
-\shave{\,\sum^\gamma_{\mu=1}\sum^2_{w_\mu=1}}
\overline{\vartheta^{\,w_\mu}_{\!j\,\bar b_\mu}}\ 
S^{\,a_1\ldots\,a_\alpha\,\bar a_1\ldots\,
\ldots\,\ldots\,\bar a_\nu\,c_1\ldots\,c_m}_{\,b_1\ldots\,b_\beta\,
\bar b_1\ldots\,w_\mu\ldots\,\bar b_\gamma\,d_1\ldots\,d_n}[P]
+\shave{\,\sum^m_{\mu=1}\sum^3_{v_\mu=0}}
\theta^{\,c_\mu}_{\!j\,v_\mu}\,\times\\
\left.\times\,S^{\,a_1\ldots\,a_\alpha\,\bar a_1\ldots\,
\bar a_\nu\,c_1\ldots\,v_\mu\ldots
\,c_m}_{\,b_1\ldots\,b_\beta\,\bar b_1\ldots\,
\bar b_\gamma\,d_1\ldots\,\ldots\,\ldots\,d_n}[P]
-\shave{\sum^n_{\mu=1}\sum^3_{w_\mu=0}}
\theta^{\,w_\mu}_{\!j\,d_\mu}
S^{\,a_1\ldots\,a_\alpha\,\bar a_1\ldots\,
\bar a_\nu\,c_1\ldots\,\ldots\,\ldots\,c_m}_{\,b_1\ldots
\,b_\beta\,\bar b_1\ldots\,\bar b_\gamma\,d_1\ldots\,w_\mu\ldots
\,d_n}[P]\right).
\endgather
$$
These formulas are again rather huge so that we do not mark them with 
a number.\par
\head
10. Extended tensorial and spin-tensorial fields.
\endhead
\mydefinition{10.1} Let $N$ be a composite spin-tensorial bundle over 
the space-time manifold $M$ (in the sense of the formula \mythetag{9.7}).
An extended tensor field $\bold X$ of the type $(e,f)$ is a
tensor-valued function in $N$ such that it takes each point $q\in N$ 
to some tensor $\bold X(q)\in T^e_f(p,M)$, where $p=\pi(q)$ is the projection
of $q$.
\enddefinition
     Informally speaking, an extended tensor field $\bold X$ is 
a tensorial function with one point argument $p\in M$ and several
spin-tensorial arguments \pagebreak $\bold S[1],\,\ldots,\,\bold S[J+Q]$ 
as shown in \mythetag{9.8}. An extended spin-tensorial field is defined
in a similar way.
\mydefinition{10.2} Let $N$ be a composite spin-tensorial bundle over 
the space-time manifold $M$ (in the sense of the formula \mythetag{9.7}).
An extended spin-tensorial field $\bold X$ of the type 
$(\varepsilon,\eta|\sigma,\zeta|e,f)$ is a spin-tensor-valued function 
in $N$ such that it takes each point $q\in N$ to some spin-tensor 
$\bold X(q)\in S^{\,\varepsilon}_\eta\bar S^\sigma_\zeta T^e_f(p,M)$, 
where $p=\pi(q)$ is the projection of $q$.
\enddefinition
     Remember that upon complexification the tensor bundle $T^e_fM$
can be identified with the spin-tensorial bundle $S^0_0\bar S^0_0T^e_fM$
(compare \mythetag{9.5} and \mythetag{9.6} above):
$$
\hskip -2em
\Bbb CT^e_fM=\Bbb C\otimes T^e_fM=S^0_0\bar S^0_0T^e_fM.
\mytag{10.1}
$$
Due to the formula \mythetag{10.1} extended tensorial fields introduced
in the definition~\mythedefinition{10.1} can be understood as some 
special cases of extended spin-tensorial fields introduced in the 
definition~\mythedefinition{10.2}. With this remark in mind, below we
do not study extended tensorial fields in a separate way. They are 
naturally included into the general framework of the study of extended 
spin-tensorial fields.\par
     In a local chart $U\subset M$ equipped with a positively polarized 
right orthonormal frame $\boldsymbol\Upsilon_0,\,\boldsymbol\Upsilon_1,
\,\boldsymbol\Upsilon_2,\,\boldsymbol\Upsilon_3$ and with its associated
spinor frame $\boldsymbol\Psi_1,\,\boldsymbol\Psi_2$ an extended
spin-tensorial field $\bold X$ is represented by its components
$X^{i_1\ldots\,i_\varepsilon\bar i_1\ldots\,\bar i_\sigma h_1\ldots
\,h_e}_{j_1\ldots\,j_\eta\bar j_1\ldots\,\bar j_\zeta k_1\ldots\,k_f}$ 
each of which is a function of the variables \mythetag{9.11} and 
\mythetag{9.12}. When passing from $U$ to an overlapping chart 
$\tilde U$ equipped with another positively polarized right orthonormal
frame $\tilde{\boldsymbol\Upsilon}_0,\,\tilde{\boldsymbol\Upsilon}_1,
\,\tilde{\boldsymbol\Upsilon}_2,\,\tilde{\boldsymbol\Upsilon}_3$ and 
associated spinor frame $\tilde{\boldsymbol\Psi}_1,\,
\tilde{\boldsymbol\Psi}_2$ the components of $\bold X$ are transformed 
according to the standard transformation formulas 
$$
\align
&\hskip -4em
\left\{
\aligned
&\tilde X^{i_1\ldots\,i_\varepsilon\bar i_1\ldots\,\bar i_\sigma
h_1\ldots\,h_e}_{j_1\ldots\,j_\eta\bar j_1\ldots\,\bar j_\zeta
k_1\ldots\,k_f}
=\dsize\msum{2}\Sb a_1,\,\ldots,\,a_\varepsilon\\ b_1,\,\ldots,\,b_\eta\endSb
\dsize\msum{2}\Sb \bar a_1,\,\ldots,\,\bar a_\sigma\\ 
\bar b_1,\,\ldots,\,\bar b_\zeta\endSb
\dsize\msum{3}
\Sb c_1,\,\ldots,\,c_e\\ d_1,\,\ldots,\,d_f\endSb
\goth T^{\,i_1}_{a_1}\ldots\,\goth T^{\,i_\varepsilon}_{a_\varepsilon}\,
\times\\
&\kern 40pt
\times\,\goth S^{b_1}_{j_1}\ldots\,\goth S^{b_\eta}_{j_\eta}
\ \overline{\goth T^{\,\bar i_1}_{\bar a_1}}
\ldots\,\overline{\goth T^{\,\bar i_\sigma}_{\bar a_\sigma}}\ 
\ \overline{\goth S^{\,\bar b_1}_{\bar j_1}}
\ldots\,\overline{\goth S^{\,\bar b_\zeta}_{\bar j_\zeta}}\ 
T^{h_1}_{c_1}\ldots\,T^{h_e}_{c_e}\,\times\\
\vspace{2ex}
&\kern 60pt
\times\,S^{\,d_1}_{k_1}\ldots\,S^{\,d_f}_{k_f}\ 
X^{\,a_1\ldots\,a_\varepsilon\,\bar a_1\ldots\,\bar a_\sigma\,
c_1\ldots\,c_e}_{\,b_1\ldots\,b_\eta\,\bar b_1\ldots\,\bar b_\zeta
\,d_1\ldots\,d_f},
\endaligned
\right.
\mytag{10.2}\\
&\hskip -4em
\left\{
\aligned
&X^{i_1\ldots\,i_\varepsilon\bar i_1\ldots\,\bar i_\sigma
h_1\ldots\,h_e}_{j_1\ldots\,j_\eta\bar j_1\ldots\,\bar j_\zeta
k_1\ldots\,k_f}
=\dsize\msum{2}\Sb a_1,\,\ldots,\,a_\varepsilon\\ b_1,\,\ldots,\,b_\eta\endSb
\dsize\msum{2}\Sb \bar a_1,\,\ldots,\,\bar a_\sigma\\ 
\bar b_1,\,\ldots,\,\bar b_\zeta\endSb
\dsize\msum{3}
\Sb c_1,\,\ldots,\,c_e\\ d_1,\,\ldots,\,d_f\endSb
\goth S^{\,i_1}_{a_1}\ldots\,\goth S^{\,i_\varepsilon}_{a_\varepsilon}\,
\times\\
&\kern 40pt
\times\,\goth T^{b_1}_{j_1}\ldots\,\goth T^{b_\eta}_{j_\eta}
\ \overline{\goth S^{\,\bar i_1}_{\bar a_1}}
\ldots\,\overline{\goth S^{\,\bar i_\sigma}_{\bar a_\sigma}}\ 
\ \overline{\goth T^{\,\bar b_1}_{\bar j_1}}
\ldots\,\overline{\goth T^{\,\bar b_\zeta}_{\bar j_\zeta}}\ 
S^{h_1}_{c_1}\ldots\,S^{h_e}_{c_e}\,\times\\
\vspace{2ex}
&\kern 60pt
\times\,T^{\,d_1}_{k_1}\ldots\,T^{\,d_f}_{k_f}\ 
\tilde X^{\,a_1\ldots\,a_\varepsilon\,\bar a_1\ldots\,\bar a_\sigma\,
c_1\ldots\,c_e}_{\,b_1\ldots\,b_\eta\,\bar b_1\ldots\,\bar b_\zeta
\,d_1\ldots\,d_f},
\endaligned
\right.
\mytag{10.3}
\endalign
$$
while their arguments are transformed according to the formulas
\mythetag{9.13}, \mythetag{9.14}, and \mythetag{9.15} as described 
above in section~9.
\mydefinition{10.3}  An extended spin-tensorial field $\bold X$ 
associated with a composite spin-tensorial bundle $N$ is called 
{\it smooth\/} if its components are smooth function of their 
arguments \mythetag{9.11} and \mythetag{9.12} in any local chart $U$ 
equipped with a positively polarized right orthonormal frame
$\boldsymbol\Upsilon_0,\,\boldsymbol\Upsilon_1,\,\boldsymbol\Upsilon_2,
\,\boldsymbol\Upsilon_3$ and associated frame $\boldsymbol\Psi_1,
\,\boldsymbol\Psi_2$.
\enddefinition
\head
11. The algebra of extended spin-tensorial fields.
\endhead
    Suppose that some composite spin-tensorial bundle $N$ given
by the formula \mythetag{9.7} is fixed. Let's denote by
$S^\varepsilon_{\!\eta}\bar S^\sigma_\zeta T^{\,e}_f(M)$ the set 
of all smooth extended spin-tensorial fields of the type
$(\varepsilon,\eta|\sigma,\zeta|e,f)$. Then the direct sum
$$
\hskip -2em
\bold S(M)=\bigoplus^\infty_{\varepsilon=0}\bigoplus^\infty_{\eta=0}
\bigoplus^\infty_{\sigma=0}\bigoplus^\infty_{\zeta=0}
\bigoplus^\infty_{e=0}\bigoplus^\infty_{f=0}S^\varepsilon_{\!\eta}
\bar S^\sigma_\zeta T^{\,e}_f(M)
\mytag{11.1}
$$
is called the {\it algebra of extended spin-tensorial fields}. It is an
algebra over the ring 
$$
\hskip -2em
S^0_0\bar S^0_0 T^0_0(M)=\goth F_{\Bbb C}(N)
\mytag{11.2}
$$
of smooth complex functions in $N$. The graded algebra \mythetag{11.1}
is equipped with the following four algebraic operations:
\roster
\item $S^\varepsilon_{\!\eta}\bar S^\sigma_\zeta T^{\,e}_f(M)+
S^\varepsilon_{\!\eta}\bar S^\sigma_\zeta T^{\,e}_f(M)\longrightarrow
S^\varepsilon_{\!\eta}\bar S^\sigma_\zeta T^{\,e}_f(M)$;
\item $S^0_0\bar S^0_0 T^0_0(M)\otimes 
S^\varepsilon_{\!\eta}\bar S^\sigma_\zeta T^{\,e}_f(M)
\longrightarrow
S^\varepsilon_{\!\eta}\bar S^\sigma_\zeta T^{\,e}_f(M)$;
\item $S^\varepsilon_{\!\eta}\bar S^\sigma_\zeta T^{\,e}_f(M)
\otimes S^\alpha_\beta\!\bar S^\nu_\gamma T^m_n(M)
\longrightarrow
S^{\varepsilon+\alpha}_{\!\eta+\beta}\bar S^{\sigma+\nu}_{\zeta+\gamma}
T^{\,e+m}_{f+n}(M)$;
\item $C\!:\,S^\varepsilon_{\!\eta}\bar S^\sigma_\zeta T^{\,e}_f(M)
\longrightarrow S^{\varepsilon-1}_{\!\eta-1}\bar S^\sigma_\zeta
T^{\,e}_f(M)$ for $\varepsilon\geqslant 1$ and $\eta\geqslant 1$,
\item" " $C\!:\,S^\varepsilon_{\!\eta}\bar S^\sigma_\zeta T^{\,e}_f(M)
\longrightarrow S^\varepsilon_{\!\eta}\bar S^{\sigma-1}_{\zeta-1}
T^{\,e}_f(M)$ for $\sigma\geqslant 1$ and $\zeta\geqslant 1$,
\item" " $C\!:\,S^\varepsilon_{\!\eta}\bar S^\sigma_\zeta T^{\,e}_f(M)
\longrightarrow S^\varepsilon_{\!\eta}\bar S^\sigma_\zeta
T^{\,e-1}_{f-1}(M)$ for $e\geqslant 1$ and $f\geqslant 1$.
\endroster
These operations are called {\it addition}, {\it multiplication by
scalars}, {\it tensor product}, and {\it contraction}. The last item
\therosteritem{4} is subdivided into three parts indicating that the
contraction operation is allowed only within certain groups of
indices, i\.\,e\. a spinor index can be contracted only with a spinor 
index, a barred spinor index only with another barred spinor index, and 
a tensorial index with another tensorial index.
\head
12. Differentiation of extended spin-tensorial fields.
\endhead
    Suppose that some composite spin-tensorial bundle $N$ over the 
space-time manifold $M$ is fixed (see \mythetag{9.7}). Then the algebra
$\bold S(M)$ is also fixed.
\mydefinition{12.1} A mapping $D\!:\,\bold S(M)\to\bold S(M)$ is called 
a {\it differentiation} of the algebra of extended spin-tensorial fields 
if the following conditions are fulfilled:
\roster
\rosteritemwd=10pt
\item $D$ is concordant with the grading: $D(S^\varepsilon_{\!\eta}
\bar S^\sigma_\zeta T^{\,e}_f(M))\subset S^\varepsilon_{\!\eta}
\bar S^\sigma_\zeta T^{\,e}_f(M)$;
\item $D$ is $\Bbb C$-linear: \hskip -0.95cm\vtop{\hsize 6.4cm
      \noindent $D(\bold X+\bold Y)=D(\bold X)+D(\bold Y)$
      \newline and $D(\lambda\bold X)=\lambda D(\bold X)$
      for $\lambda\in\Bbb C$;\vskip 1.3ex}
\item $D$ commutates with the contractions: $D(C(\bold X))=C(D(\bold X))$;
\item $D$ obeys the Leibniz rule: $D(\bold X\otimes\bold Y)=D(\bold X)
      \otimes\bold Y+\bold X\otimes D(\bold Y)$.
\endroster
\enddefinition
\noindent Let's consider the set of all differentiations of the extended
algebra $\bold S(M)$. We denote it $\goth D_{\bold S}(M)$. It is easy to
check up that 
\roster
\rosteritemwd=10pt
\item the sum of two differentiations is a differentiation of the algebra
$\bold S(M)$;
\item the product of a differentiation by a smooth complex function in $N$
is a differentiation of the algebra $\bold S(M)$.
\endroster
Now we see that $\goth D_{\bold S}(M)$ is equipped with the structure 
of a module over the ring of smooth complex functions $\goth F_{\Bbb C}(N)$
(see \mythetag{11.2}).\par
     In the module $\goth D_{\bold S}(M)$ the composition of  two
differentiations $D_1$ and $D_2$ is not a differentiation, but their
commutator 
$$
[D_1,\,D_2]=D_1\compos D_2-D_2\compos D_1
$$
is a differentiation. Therefore, $\goth D_{\bold S}(M)$ is a Lie algebra. 
Note, however, that $\goth D_{\bold S}(M)$ is not a Lie algebra over the 
ring of smooth complex functions $\goth F_{\Bbb C}(N)$. It is only a Lie
algebra over the field of complex numbers $\Bbb C$.\par
\head
13. Localization.
\endhead
    The results of this section are very similar to those in 
section~7 of the paper \mycite{5}. Nevertheless, for the reader's 
convenience we give them in full details.\par
    Smooth extended spin-tensorial fields are global objects related 
to the spin-tensorial bundle $N$ in whole. But they are functions and 
their values are local objects. This means that two different fields 
$\bold A\neq\bold B$ can take the same values at some particular points.
Whenever this happens, we write $\bold A_q=\bold B_q$, where $q\in N$ 
is a point of the composite spin-tensorial bundle \mythetag{9.7}.\par
    Differentiations of the algebra $\bold S(M)$, as they are introduced 
above in the definitions~\mythedefinition{12.1}, are global objects without
any explicit subdivision into parts related to separate points of the
bundle $N$. Below in this section we shall show that they also can be
represented as functions taking their values in some linear spaces
associated with separate points of the manifold $N$.\par
    Let $D\in\goth D_{\bold S}(M)$ be a differentiation of the algebra of 
extended spin-tensorial fields $\bold S(M)$. Let's denote by $\delta$ 
the restriction of the mapping $D\!:\bold S(M)\to\bold S(M)$ to the module
$S^0_0\bar S^0_0 T^0_0(M)$ of extended scalar fields in \mythetag{11.1}:
$$
\hskip -2em
\delta\!:\,S^0_0\bar S^0_0 T^0_0(M)\to S^0_0\bar S^0_0 T^0_0(M).
\mytag{13.1}
$$
Since $S^0_0\bar S^0_0 T^0_0(M)=\goth F_{\Bbb C}(N)$, the mapping $\delta$
in \mythetag{13.1} is a differentiation of the ring of smooth complex
functions in the smooth real manifold $N$. It is known (see \S\,1 in
Chapter~\uppercase\expandafter{\romannumeral 1} of \mycite{11}) 
that any differentiation of the ring of smooth functions 
of an arbitrary smooth manifold is determined by some vector field 
$\bold Z$ in this manifold. In our case $\bold Z$ is a complexified
vector field, i\.\,e\.
$$
\bold Z\in\Bbb CT^1_0(N)=\Bbb C\otimes T^1_0(N).
$$
Applying this fact, we get the following representation for the operator
$\delta$ in terms of the differential operators \mythetag{9.26}, 
\mythetag{9.19}, and \mythetag{9.20}:
$$
\pagebreak
\gathered
\delta=\bold Z=\sum^3_{i=0}Z^i\ \bold U_i\,+\\
+\sum^{J+Q}_{P=1}\dsize\msums{2}{3}\Sb i_1,\,\ldots,\,i_\alpha\\
j_1,\,\ldots,\,j_\beta\\ \bar i_1,\,\ldots,\,\bar i_\nu\\ 
\bar j_1,\,\ldots,\,\bar j_\gamma\\
h_1,\,\ldots,\,h_m\\ k_1,\,\ldots,\,k_n\endSb
\left(Z^{i_1\ldots\,i_\alpha\,\bar i_1\ldots\,\bar i_\nu\,h_1\ldots
\,h_m}_{j_1\ldots\,j_\beta\,\bar j_1\ldots\,\bar j_\gamma\,k_1\ldots
\,k_n}[P]\ \bold W^{j_1\ldots\,j_\beta\,\bar j_1\ldots\,
\bar j_\gamma\,k_1\ldots\,k_n}_{i_1\ldots\,i_\alpha\,\bar i_1\ldots\,
\bar i_\nu\,h_1\ldots\,h_m}[P]\,+\right.\\
\left.\kern 70pt+\,\bar Z^{\bar i_1\ldots\,\bar i_\nu\,i_1\ldots
\,i_\alpha\,h_1\ldots\,h_m}_{\bar j_1\ldots\,\bar j_\gamma\,j_1
\ldots\,j_\beta\,k_1\ldots\,k_n}[P]\ \bar{\bold W}^{\bar j_1\ldots
\,\bar j_\gamma\,j_1\ldots\,j_\beta\,k_1\ldots\,k_n}_{\bar i_1\ldots
\,\bar i_\nu\,i_1\ldots\,i_\alpha\,h_1\ldots\,h_m}[P]\right).
\endgathered\qquad
\mytag{13.2}
$$
Note that $Z^i$, $Z^{i_1\ldots\,i_\alpha\,\bar i_1\ldots\,\bar i_\nu
\,h_1\ldots\,h_m}_{j_1\ldots\,j_\beta\,\bar j_1\ldots\,\bar j_\gamma\,
k_1\ldots\,k_n}[P]$, and $\bar Z^{\bar i_1\ldots\,\bar i_\nu\,i_1\ldots
\,i_\alpha\,h_1\ldots\,h_m}_{\bar j_1\ldots\,\bar j_\gamma\,j_1\ldots
\,j_\beta\,k_1\ldots\,k_n}[P]$ in \mythetag{13.2} are arbitrary smooth
complex functions within the the local chart $U$ where the representation
\mythetag{13.2} is defined.
\mylemma{13.1} Let $\psi$ be an extended scalar field (a smooth function)
identically constant within some open subset $O\subset N$ and let 
$\varphi=D(\psi)$ for some differentiation $D$. Then $\varphi=0$ within 
the open set $O$.
\endproclaim
\demo{Proof} Since $D(\psi)=\delta(\psi)$, choosing some local 
chart and applying the differential operators \mythetag{9.26}, 
\mythetag{9.19}, and \mythetag{9.20} to a constant, from the 
formula \mythetag{13.2} we derive that $\delta(\psi)=\bold Z\psi=0$ 
at any point $q$ of the open set $O$.\qed\enddemo
\mylemma{13.2} Let\/ $\bold X$ be an extended spin-tensorial field. 
If\/ $\bold X\equiv 0$, then for any differentiation $D$ the field 
$D(\bold X)$ is also identically equal to zero.
\endproclaim
    The proof is trivial. Since $\bold X\equiv 0$, we can write $\bold
X=\lambda\,\bold X$ with $\lambda\neq 1$. Then, applying the item
\therosteritem{2} of the definition~\mythedefinition{12.1}, we get 
$D(\bold X)=\lambda\,D(\bold X)$. Since $\lambda\neq 1$, this yields 
the required equality $D(\bold X)\equiv 0$.
\mylemma{13.3} Let $\bold X$ be an extended spin-tensorial field
identically zero within some open set $O\subset N$. If\/ $\bold Y
=D(\bold X)$ for some differentiation $D$, then $\bold Y_{\!q}=0$ 
at any point $q\in O$.
\endproclaim
\demo{Proof} Let's choose some arbitrary point $q\in O$ and take some
smooth scalar function $\eta$ such that it is identically equal to the
unity in some open neighborhood $O{\,'}\subset O$ of the point $q$ and
identically equal to zero outside the open set $O$. The existence
of such a function $\eta$ is easily proved by choosing some local chart 
$U$ that covers the point $q$. The product $\eta\,\bold X$ is identically
equal to zero:
$$
\hskip -2em
\eta\otimes\bold X=\eta\,\bold X\equiv 0.
\mytag{13.3}
$$
Applying the differentiation $D$ to \mythetag{13.3}, then taking into 
account the lemma~\mythelemma{13.2} and the item \therosteritem{4} of 
the definition~\mythedefinition{12.1}, we obtain
$$
0=D(0)=D(\eta\otimes\bold X)=D(\eta)\otimes\bold X+\eta\otimes D(\bold X)
=D(\eta)\,\bold X+\eta\,D(\bold X).
$$
Note that $D(\eta)=0$ at the point $q$ due to the lemma~\mythelemma{13.1}.
Moreover, $\bold X_q=0$ and $\eta=1$ at the point $q$. Therefore, by
specifying the above equality to the point $q$ we get $D(\bold X)=0$ at
the point $q$. The lemma is proved.\qed\enddemo
\mylemma{13.4} If two extended spin-tensorial fields $\bold X$ and 
$\bold Y$ are equal within some open neighborhood $O$ of a point 
$q\in N$, then for any differentiation $D$ their images $D(\bold X)$ 
and $D(\bold Y)$ are equal at the point $q$.
\endproclaim
    The lemma~\mythelemma{13.4} follows immediately from the 
lemma~\mythelemma{13.3}. This lemma is a basic tool for our purposes
of localization in this section.\par
    Let $q$ be some point of $N$ and let $p=\pi(q)$ be its projection
in the base manifold $M$. Suppose that $U$ is a local chart of $M$
equipped with a positively polarized right oriented orthonormal
frame $\boldsymbol\Upsilon_0,\,\boldsymbol\Upsilon_1,\,
\boldsymbol\Upsilon_2,\,\boldsymbol\Upsilon_3$ and with its associated 
spinor frame $\boldsymbol\Psi_1,\,\boldsymbol\Psi_2$ and such that 
it covers the point $p$ in $M$. Then we can use its preimage 
$\pi^{-1}(U)$ as a local chart in $N$ covering the point $q$. The 
variables \mythetag{9.11} and \mythetag{9.12} form a complete set of
local coordinates in the chart $\pi^{-1}(U)$. Any extended spin-tensorial
field $\bold X$ of the type $(\varepsilon,\eta|\sigma,\zeta|e,f)$ is
represented  by the formula 
$$
\hskip -2em
\bold X=
\dsize\msums{2}{3}\Sb a_1,\,\ldots,\,a_\varepsilon\\
b_1,\,\ldots,\,b_\eta\\ \bar a_1,\,\ldots,\,\bar a_\sigma\\ 
\bar b_1,\,\ldots,\,\bar b_\zeta\\
c_1,\,\ldots,\,c_e\\ d_1,\,\ldots,\,d_f\endSb
X^{a_1\ldots\,a_\varepsilon\bar a_1\ldots\,\bar a_\sigma
h_1\ldots\,h_e}_{b_1\ldots\,b_\eta\bar b_1\ldots\,\bar b_\zeta
k_1\ldots\,k_f}\ 
\boldsymbol\Psi^{b_1\ldots\,b_\eta\bar b_1\ldots\,\bar b_\zeta
d_1\ldots\,d_f}_{a_1\ldots\,a_\varepsilon\bar a_1\ldots\,\bar a_\sigma
c_1\ldots\,c_e},
\mytag{13.4}
$$
where $\boldsymbol\Psi^{b_1\ldots\,b_\eta\bar b_1\ldots\,\bar b_\zeta
d_1\ldots\,d_f}_{a_1\ldots\,a_\varepsilon\bar a_1\ldots\,\bar a_\sigma
c_1\ldots\,c_e}$ are given by the tensor products \mythetag{5.33}. In 
the case of standard (traditional) spin-tensorial fields the coefficients
$X^{a_1\ldots\,a_\varepsilon\bar a_1\ldots\,\bar a_\sigma
c_1\ldots\,c_e}_{b_1\ldots\,b_\eta\bar b_1\ldots\,\bar b_\zeta
d_1\ldots\,d_f}$ in \mythetag{13.4} depend on the coordinates $x^0,\,x^1,
\,x^2,\,x^3$ of a point $p\in M$ only (see \mythetag{4.3} and compare 
\mythetag{5.34} with \mythetag{13.4}). In the case of extended
fields they depend on the whole set of variables \mythetag{9.11} and
\mythetag{9.12}.\par
     Taking some differentiations, we can apply them to the left hand 
side of \mythetag{13.4}, but we cannot apply them to each summand in the
right hand side of this formula. The matter is that the scalars 
$X^{a_1\ldots\,a_\varepsilon\bar a_1\ldots\,\bar a_\sigma c_1\ldots\,
c_e}_{b_1\ldots\,b_\eta\bar b_1\ldots\,\bar b_\zeta d_1\ldots\,d_f}$
and the spin-tensors $\boldsymbol\Psi^{b_1\ldots\,b_\eta\bar b_1\ldots\,
\bar b_\zeta d_1\ldots\,d_f}_{a_1\ldots\,a_\varepsilon\bar a_1\ldots\,
\bar a_\sigma c_1\ldots\,c_e}$ are defined locally only within the chart
$\pi^{-1}(U)$. Therefore, they do not fit the
definition~\mythedefinition{10.2}. In order to convert them to global 
fields we choose some smooth real function $\eta\in\goth F_{\Bbb R}(N)
\subset\goth F_{\Bbb C}(N)$ such that it is identically equal to the unity
within some open neighborhood of the point $q$ and is identically zero
outside the chart $\pi^{-1}(U)$. Then we define the following global 
extended fields:
$$
\align
\hskip -4em
\hat X^{a_1\ldots\,a_\varepsilon\bar a_1\ldots\,\bar a_\sigma
c_1\ldots\,c_e}_{b_1\ldots\,b_\eta\bar b_1\ldots\,\bar b_\zeta
d_1\ldots\,d_f}
&=\cases
\eta\ X^{a_1\ldots\,a_\varepsilon\bar a_1\ldots\,\bar a_\sigma
c_1\ldots\,c_e}_{b_1\ldots\,b_\eta\bar b_1\ldots\,\bar b_\zeta
d_1\ldots\,d_f} &\text{within \ }\pi^{-1}(U),\\
\ \ 0&\text{outside \ }\pi^{-1}(U),\endcases
\mytag{13.5}\\
\vspace{2ex}
\hskip -2em
\hat{\boldsymbol\Upsilon}_i&=\cases\eta\ \boldsymbol\Upsilon_i 
&\text{within \ }\pi^{-1}(U),\\
\ \ 0&\text{outside \ }\pi^{-1}(U),\endcases
\mytag{13.6}\\
\vspace{2ex}
\hskip -2em
\hat{\boldsymbol\eta}^i&=\cases\eta\
\boldsymbol\eta^i &\text{within \ }\pi^{-1}(U),\\
\ \ 0&\text{outside \ }\pi^{-1}(U),\endcases
\mytag{13.7}\\
\vspace{2ex}
\hskip -2em
\hat{\boldsymbol\Psi}_i&=\cases\eta\
\boldsymbol\Psi_i &\text{within \ }\pi^{-1}(U),\\
\ \ 0&\text{outside \ }\pi^{-1}(U),\endcases
\mytag{13.8}\\
\vspace{2ex}
\hskip -2em
\hat{\boldsymbol\vartheta}^i&=\cases\eta\
\boldsymbol\vartheta^i &\text{within \ }\pi^{-1}(U),\\
\ \ 0&\text{outside \ }\pi^{-1}(U).\endcases
\mytag{13.9}
\endalign
$$
Then by analogy to \mythetag{5.28} we write
$$
\hskip -2em
\aligned
&\hat{\bPsi}_i=\tau(\hat{\boldsymbol\Psi}_i),\\
&\hat{\overline{\boldsymbol\vartheta}}\vphantom{\vartheta}^{\,i}
=\tau(\hat{\boldsymbol\vartheta}^{\,i})
\endaligned
\mytag{13.10}
$$
and by analogy to \mythetag{5.30}, \mythetag{5.31}, \mythetag{5.32}, 
and \mythetag{5.33} \pagebreak we introduce the following tensor products, 
which are global extended fields:
$$
\align
&\hskip -2em
\hat{\boldsymbol\Upsilon}^{\,d_1\ldots\,d_f}_{c_1\ldots\,c_e}
=\hat{\boldsymbol\Upsilon}_{c_1}\otimes\ldots\otimes
\hat{\boldsymbol\Upsilon}_{c_e}
\otimes\hat{\boldsymbol\eta}^{d_1}\otimes\ldots\otimes
\hat{\boldsymbol\eta}^{d_f},
\mytag{13.11}\\
\vspace{1ex}
&\hskip -2em
\hat{\boldsymbol\Psi}^{b_1\ldots\,b_\eta}_{a_1\ldots\,a_\varepsilon}
=\hat{\boldsymbol\Psi}_{a_1}\otimes\ldots\otimes
\hat{\boldsymbol\Psi}_{a_\varepsilon}\otimes
\hat{\boldsymbol\vartheta}^{\,b_1}
\otimes\ldots\otimes
\hat{\boldsymbol\vartheta}^{\,b_\eta},
\mytag{13.12}\\
\vspace{1ex}
&\hskip -2em
\hat{\bPsi}\vphantom{\bPsi}^{\bar b_1\ldots\,\bar b_\zeta}_{\bar a_1
\ldots\,\bar a_\sigma}
=\hat{\bPsi}_{\bar a_1}\otimes\ldots\otimes
\hat{\bPsi}_{\bar a_\sigma}\otimes
\hat{\overline{\boldsymbol\vartheta}}\vphantom{\vartheta}^{\,b_1}
\otimes\ldots\otimes
\hat{\overline{\boldsymbol\vartheta}}\vphantom{\vartheta}^{\,b_\zeta}.
\mytag{13.13}\\
\vspace{1ex}
&\hskip -3em
\hat{\boldsymbol\Psi}^{b_1\ldots\,b_\eta\bar b_1\ldots\,\bar b_\zeta
d_1\ldots\,d_f}_{a_1\ldots\,a_\varepsilon\bar a_1\ldots\,\bar a_\sigma
c_1\ldots\,c_e}
=\hat{\boldsymbol\Psi}^{b_1\ldots\,b_\eta}_{a_1\ldots\,a_\varepsilon}
\otimes\hat{\bPsi}\vphantom{\bPsi}^{\bar b_1\ldots\,
\bar b_\zeta}_{\bar a_1\ldots\,\bar a_\sigma}\otimes
\hat{\boldsymbol\Upsilon}^{\,d_1\ldots\,d_f}_{c_1\ldots\,c_e}.
\mytag{13.14}
\endalign
$$
Taking into account \mythetag{13.5} and \mythetag{13.14}, from 
\mythetag{13.4} we derive
$$
\eta^{\,\omega}\ \bold X=
\dsize\msums{2}{3}\Sb a_1,\,\ldots,\,a_\varepsilon\\
b_1,\,\ldots,\,b_\eta\\ \bar a_1,\,\ldots,\,\bar a_\sigma\\ 
\bar b_1,\,\ldots,\,\bar b_\zeta\\
c_1,\,\ldots,\,c_e\\ d_1,\,\ldots,\,d_f\endSb
\hat X^{a_1\ldots\,a_\varepsilon\bar a_1\ldots\,\bar a_\sigma
c_1\ldots\,c_e}_{b_1\ldots\,b_\eta\bar b_1\ldots\,\bar b_\zeta
d_1\ldots\,d_f}\ 
\hat{\boldsymbol\Psi}^{b_1\ldots\,b_\eta\bar b_1\ldots\,\bar b_\zeta
d_1\ldots\,d_f}_{a_1\ldots\,a_\varepsilon\bar a_1\ldots\,\bar a_\sigma
c_1\ldots\,c_e},\quad
\mytag{13.15}
$$
where $\omega=\varepsilon+\eta+\sigma+\zeta+e+f+1$.
Now we can apply a differentiation $D$ not only to left hand side
of \mythetag{13.15}, but to each summand in the right hand side of 
this equality. Using the item \therosteritem{4} of the
definitions~\mythedefinition{12.1}, from \mythetag{13.15} we derive
$$
\gathered
D(\eta^\omega\ \bold X)=
\dsize\msums{2}{3}\Sb a_1,\,\ldots,\,a_\varepsilon\\
b_1,\,\ldots,\,b_\eta\\ \bar a_1,\,\ldots,\,\bar a_\sigma\\ 
\bar b_1,\,\ldots,\,\bar b_\zeta\\
c_1,\,\ldots,\,c_e\\ d_1,\,\ldots,\,d_f\endSb
D(\hat X^{a_1\ldots\,a_\varepsilon\bar a_1\ldots\,\bar a_\sigma
c_1\ldots\,c_e}_{b_1\ldots\,b_\eta\bar b_1\ldots\,\bar b_\zeta
d_1\ldots\,d_f})\ 
\hat{\boldsymbol\Psi}^{b_1\ldots\,b_\eta\bar b_1\ldots\,\bar b_\zeta
d_1\ldots\,d_f}_{a_1\ldots\,a_\varepsilon\bar a_1\ldots\,\bar a_\sigma
c_1\ldots\,c_e}\,+\\
\vspace{2ex}
+\dsize\msums{2}{3}\Sb a_1,\,\ldots,\,a_\varepsilon\\
b_1,\,\ldots,\,b_\eta\\ \bar a_1,\,\ldots,\,\bar a_\sigma\\ 
\bar b_1,\,\ldots,\,\bar b_\zeta\\
c_1,\,\ldots,\,c_e\\ d_1,\,\ldots,\,d_f\endSb
\hat X^{a_1\ldots\,a_\varepsilon\bar a_1\ldots\,\bar a_\sigma
c_1\ldots\,c_e}_{b_1\ldots\,b_\eta\bar b_1\ldots\,\bar b_\zeta
d_1\ldots\,d_f}\ 
D(\hat{\boldsymbol\Psi}^{b_1\ldots\,b_\eta\bar b_1\ldots\,\bar b_\zeta
d_1\ldots\,d_f}_{a_1\ldots\,a_\varepsilon\bar a_1\ldots\,\bar a_\sigma
c_1\ldots\,c_e}).
\endgathered\quad
\mytag{13.16}
$$
Due to the lemma~\mythelemma{13.4} we have $D(\eta^{\,\omega}\ \bold X)
=D(\bold X)$ at the point $q$. Moreover, due to \mythetag{13.5},
\mythetag{13.6}, \mythetag{13.7}, \mythetag{13.8}, \mythetag{13.9},
\mythetag{13.10}, \mythetag{13.11}, \mythetag{13.12}, \mythetag{13.13},
and \mythetag{13.14} and since $\eta(q)=1$, we have
$\hat{\boldsymbol\Psi}^{b_1\ldots\,b_\eta\bar b_1\ldots\,\bar b_\zeta
d_1\ldots\,d_f}_{a_1\ldots\,a_\varepsilon\bar a_1\ldots\,\bar a_\sigma
c_1\ldots\,c_e}=\boldsymbol\Psi^{b_1\ldots\,b_\eta\bar b_1\ldots\,
\bar b_\zeta d_1\ldots\,d_f}_{a_1\ldots\,a_\varepsilon\bar a_1\ldots\,
\bar a_\sigma c_1\ldots\,c_e}$ at this point. As for the fields 
$D(\hat X^{a_1\ldots\,a_\varepsilon\bar a_1\ldots\,\bar a_\sigma
c_1\ldots\,c_e}_{b_1\ldots\,b_\eta\bar b_1\ldots\,\bar b_\zeta d_1
\ldots\,d_f})$ and $D(\hat{\boldsymbol\Psi}^{b_1\ldots\,b_\eta\bar b_1
\ldots\,\bar b_\zeta d_1\ldots\,d_f}_{a_1\ldots\,a_\varepsilon\bar a_1
\ldots\,\bar a_\sigma c_1\ldots\,c_e})$ in \mythetag{13.16}, again due 
to the lemma~\mythelemma{13.4} their values at the point $q$ do not 
depend on a particular choice of the function $\eta$.\par
     Since $\hat X^{a_1\ldots\,a_\varepsilon\bar a_1\ldots\,\bar a_\sigma
c_1\ldots\,c_e}_{b_1\ldots\,b_\eta\bar b_1\ldots\,\bar b_\zeta d_1\ldots\,
d_f}$ is a function, i\.\,e\. it is interpreted as an extended scalar
field, we have $D(\hat X^{a_1\ldots\,a_\varepsilon\bar a_1
\ldots\,\bar a_\sigma c_1\ldots\,c_e}_{b_1\ldots\,b_\eta\bar b_1\ldots\,
\bar b_\zeta d_1\ldots\,d_f})=\delta(\hat X^{a_1\ldots\,a_\varepsilon
\bar a_1\ldots\,\bar a_\sigma c_1\ldots\,c_e}_{b_1\ldots\,b_\eta\bar b_1
\ldots\,\bar b_\zeta d_1\ldots\,d_f})$. Then from the formulas 
\mythetag{13.2}, \mythetag{9.26}, \mythetag{9.18}, \mythetag{9.19}, and
\mythetag{9.20} for the value of the function $D(\hat X^{a_1\ldots\,
a_\varepsilon\bar a_1\ldots\,\bar a_\sigma c_1\ldots\,c_e}_{b_1\ldots\,
b_\eta\bar b_1\ldots\,\bar b_\zeta d_1\ldots\,d_f})$ at the point $q$ we
derive the following formula:
$$
\gathered
D(\hat X^{a_1\ldots\,a_\varepsilon\bar a_1\ldots\,\bar a_\sigma c_1\ldots\,
c_e}_{b_1\ldots\,b_\eta\bar b_1\ldots\,\bar b_\zeta d_1\ldots\,d_f})
=\sum^3_{i=0}\sum^3_{j=0}Z^i\,\Upsilon^j_i\,
\frac{\partial X^{a_1\ldots\,a_\varepsilon\bar a_1\ldots\,\bar a_\sigma
c_1\ldots\,c_e}_{b_1\ldots\,b_\eta\bar b_1\ldots\,\bar b_\zeta d_1\ldots\,
d_f}}{\partial x^j}\,+\\
+\sum^{J+Q}_{P=1}\dsize\msums{2}{3}\Sb i_1,\,\ldots,\,i_\alpha\\
j_1,\,\ldots,\,j_\beta\\ \bar i_1,\,\ldots,\,\bar i_\nu\\ 
\bar j_1,\,\ldots,\,\bar j_\gamma\\
h_1,\,\ldots,\,h_m\\ k_1,\,\ldots,\,k_n\endSb
Z^{i_1\ldots\,i_\alpha\,\bar i_1\ldots\,\bar i_\nu\,h_1\ldots
\,h_m}_{j_1\ldots\,j_\beta\,\bar j_1\ldots\,\bar j_\gamma\,k_1\ldots
\,k_n}[P]\ 
\frac{\partial X^{a_1\ldots\,a_\varepsilon\bar a_1\ldots\,\bar a_\sigma
c_1\ldots\,c_e}_{b_1\ldots\,b_\eta\bar b_1\ldots\,\bar b_\zeta d_1\ldots\,
d_f}}{\partial S^{\,i_1\ldots\,i_\alpha\,\bar i_1\ldots\,\bar i_\nu\,
h_1\ldots\,h_m}_{j_1\ldots\,j_\beta\,\bar j_1\ldots\,\bar j_\gamma\,
k_1\ldots\,k_n}[P]}\,+\\
+\sum^{J+Q}_{P=1}\dsize\msums{2}{3}\Sb i_1,\,\ldots,\,i_\alpha\\
j_1,\,\ldots,\,j_\beta\\ \bar i_1,\,\ldots,\,\bar i_\nu\\ 
\bar j_1,\,\ldots,\,\bar j_\gamma\\
h_1,\,\ldots,\,h_m\\ k_1,\,\ldots,\,k_n\endSb
\bar Z^{\bar i_1\ldots\,\bar i_\nu\,i_1\ldots\,i_\alpha\,h_1\ldots
\,h_m}_{\bar j_1\ldots\,\bar j_\gamma\,j_1\ldots\,j_\beta\,k_1\ldots
\,k_n}[P]\ 
\frac{\partial X^{a_1\ldots\,a_\varepsilon\bar a_1\ldots\,\bar a_\sigma
c_1\ldots\,c_e}_{b_1\ldots\,b_\eta\bar b_1\ldots\,\bar b_\zeta d_1\ldots\,
d_f}}{\partial\overline{S^{\,i_1\ldots\,i_\alpha\,\bar i_1\ldots\,
\bar i_\nu\,h_1\ldots\,h_m}_{j_1\ldots\,j_\beta\,\bar j_1\ldots\,
\bar j_\gamma\,k_1\ldots\,k_n}[P]}\,}.
\endgathered\qquad
\mytag{13.17}
$$
In order to evaluate $D(\hat{\boldsymbol\Psi}^{b_1\ldots\,b_\eta\bar b_1
\ldots\,\bar b_\zeta d_1\ldots\,d_f}_{a_1\ldots\,a_\varepsilon\bar a_1
\ldots\,\bar a_\sigma c_1\ldots\,c_e})$ at the point $q$ let's apply $D$ 
to \mythetag{13.14}. Using the item \therosteritem{4} of the 
definition~\mythedefinition{12.1}, we obtain the following equality:
$$
\hskip -0.1em
\gathered
D(\hat{\boldsymbol\Psi}^{b_1\ldots\,b_\eta\bar b_1\ldots\,\bar b_\zeta
d_1\ldots\,d_f}_{a_1\ldots\,a_\varepsilon\bar a_1\ldots\,\bar a_\sigma
c_1\ldots\,c_e})
=D(\hat{\boldsymbol\Psi}^{b_1\ldots\,b_\eta}_{a_1\ldots\,a_\varepsilon})
\otimes\hat{\bPsi}\vphantom{\bPsi}^{\bar b_1\ldots\,
\bar b_\zeta}_{\bar a_1\ldots\,\bar a_\sigma}\otimes
\hat{\boldsymbol\Upsilon}^{d_1\ldots\,d_f}_{c_1\ldots\,c_e}\,+\\
\vspace{2ex}
+\,\hat{\boldsymbol\Psi}^{b_1\ldots\,b_\eta}_{a_1\ldots\,a_\varepsilon}
\otimes D(\hat{\bPsi}\vphantom{\bPsi}^{\bar b_1\ldots\,
\bar b_\zeta}_{\bar a_1\ldots\,\bar a_\sigma})\otimes
\hat{\boldsymbol\Upsilon}^{d_1\ldots\,d_f}_{c_1\ldots\,c_e}
+\hat{\boldsymbol\Psi}^{b_1\ldots\,b_\eta}_{a_1\ldots\,a_\varepsilon}
\otimes\hat{\bPsi}\vphantom{\bPsi}^{\bar b_1\ldots\,
\bar b_\zeta}_{\bar a_1\ldots\,\bar a_\sigma}\otimes
D(\hat{\boldsymbol\Upsilon}^{d_1\ldots\,d_f}_{c_1\ldots\,c_e}).
\endgathered
\mytag{13.18}
$$
In order to evaluate $D(\hat{\boldsymbol\Psi}^{b_1\ldots\,b_\eta}_{a_1
\ldots\,a_\varepsilon})$ in \mythetag{13.18} we apply $D$ to the equality
\mythetag{13.12}:
$$
\hskip -2em
\gathered
D(\hat{\boldsymbol\Psi}^{b_1\ldots\,b_\eta}_{a_1\ldots\,a_\varepsilon})
=\sum^\varepsilon_{v=1}\boldsymbol\Psi_{a_1}\otimes\ldots\otimes\,
D(\hat{\boldsymbol\Psi}_{a_v})\otimes\ldots\otimes
\boldsymbol\Psi_{a_\varepsilon}\otimes\\
\otimes\ \boldsymbol\vartheta^{\,b_1}\otimes\ldots\otimes\,
\boldsymbol\vartheta^{\,b_\eta}+\sum^\eta_{w=1}
\boldsymbol\Psi_{a_1}\otimes\ldots\otimes\,
\boldsymbol\Psi_{a_\varepsilon}\otimes\\
\otimes\ \boldsymbol\vartheta^{\,b_1}\otimes\ldots\otimes 
D(\hat{\boldsymbol\vartheta}^{\,b_w})\otimes\ldots
\otimes\boldsymbol\vartheta^{\,b_\eta}.\vphantom{\sum^\eta_{w=1}}
\endgathered
\mytag{13.19}
$$
Similarly, in order to evaluate $D(\hat{\bPsi}\vphantom{\bPsi}^{\bar
b_1\ldots\,\bar b_\zeta}_{\bar a_1\ldots\,\bar a_\sigma})$ in 
\mythetag{13.18} we apply $D$ to \mythetag{13.13}:
$$
\hskip -2em
\gathered
D(\hat{\bPsi}\vphantom{\bPsi}^{\bar b_1\ldots\,\bar b_\zeta}_{\bar a_1
\ldots\,\bar a_\sigma})
=\sum^\sigma_{v=1}\hat{\bPsi}_{\bar a_1}\otimes\ldots\otimes\,
D(\hat{\bPsi}_{\bar a_v})\otimes\ldots\otimes
\hat{\bPsi}_{\bar a_\sigma}\otimes\\
\otimes\ 
\hat{\overline{\boldsymbol\vartheta}}\vphantom{\vartheta}^{\,b_1}
\otimes\ldots\otimes\,
\hat{\overline{\boldsymbol\vartheta}}\vphantom{\vartheta}^{\,b_\zeta}
+\sum^\zeta_{w=1}\hat{\bPsi}_{\bar a_1}\otimes\ldots\otimes\,
\hat{\bPsi}_{\bar a_\sigma}\otimes\\
\otimes\ \hat{\overline{\boldsymbol\vartheta}}\vphantom{\vartheta}^{\,b_1}
\otimes\ldots\otimes 
D(\hat{\overline{\boldsymbol\vartheta}}\vphantom{\vartheta}^{\,b_w})
\otimes\ldots\otimes
\hat{\overline{\boldsymbol\vartheta}}\vphantom{\vartheta}^{\,b_\zeta}.
\vphantom{\sum^\zeta_{w=1}}
\endgathered
\mytag{13.20}
$$
And finally, in order to evaluate $D(\hat{\boldsymbol\Upsilon}^{d_1
\ldots\,d_f}_{c_1\ldots\,c_e})$ in \mythetag{13.18} we apply $D$ to
\mythetag{13.11}. As a result we get the following equality:
$$
\hskip -2em
\gathered
D(\hat{\boldsymbol\Upsilon}^{d_1\ldots\,
d_f}_{c_1\ldots\,c_e})=\sum^e_{v=1}
\boldsymbol\Upsilon_{c_1}\otimes\ldots\otimes\,
D(\hat{\boldsymbol\Upsilon}_{c_v})\otimes
\ldots\otimes\boldsymbol\Upsilon_{c_e}\otimes\\
\otimes\ \boldsymbol\eta^{d_1}\otimes\ldots\otimes\,\boldsymbol\eta^{d_f}
+\sum^f_{w=1}\boldsymbol\Upsilon_{c_1}\otimes\ldots\otimes
\boldsymbol\Upsilon_{c_e}\otimes\\
\otimes\ \boldsymbol\eta^{d_1}\otimes\ldots\otimes 
D(\hat{\boldsymbol\eta}^{d_w})\otimes\ldots
\otimes\boldsymbol\eta^{d_f}.\vphantom{\sum^f_{w=1}}
\endgathered
\mytag{13.21}
$$
Applying \mythetag{13.19}, \mythetag{13.20}, \mythetag{13.21} 
to \mythetag{13.18} and then substituting \mythetag{13.18} and
\mythetag{13.17} into the equality \mythetag{13.16}, we derive 
the following lemma.
\mylemma{13.5} Any differentiation $D$ of the algebra of extended 
spin-tensorial fields is uniquely fixed by its restrictions to the 
modules $S^0_0\bar S^0_0 T^0_0(M)$, $S^1_0\bar S^0_0 T^0_0(M)$,
$S^0_1\bar S^0_0 T^0_0(M)$, $S^0_0\bar S^1_0 T^0_0(M)$, 
$S^0_0\bar S^0_1 T^0_0(M)$, $S^0_0\bar S^0_0 T^1_0(M)$, and
$S^0_0\bar S^0_0 T^0_1(M)$ in \mythetag{11.1}.
\endproclaim
    In \mythetag{13.21} the value of the extended vector field 
$D(\hat{\boldsymbol\Upsilon}_{h})$ at the point $q$ is a vector 
of $\Bbb CT_{\pi(q)}(M)$. We can write the following expansion 
for this vector:
$$
\hskip -2em
D(\hat{\boldsymbol\Upsilon}_h)
=\sum^3_{k=0}\Gamma^k_h\,\boldsymbol\Upsilon_k.
\mytag{13.22}
$$
Due to the lemma~\mythelemma{13.4} the left hand side of the equality
\mythetag{13.22} does not depend on a particular choice of the function
$\eta$ in \mythetag{13.6}. Therefore, the coefficients $\Gamma^k_h$ in
\mythetag{13.22} represent the differentiation $D$ at the point $q$ for
a given local chart $U$ in $M$. The same is true for $Z^i$,
$Z^{i_1\ldots\,i_r}_{j_1\ldots\,j_s}[P]$, $Z^{i_1\ldots\,i_\alpha\,
\bar i_1\ldots\,\bar i_\nu\,h_1\ldots\,h_m}_{j_1\ldots\,j_\beta\,
\bar j_1\ldots\,\bar j_\gamma\,k_1\ldots\,k_n}[P]$, and $\bar Z^{i_1\ldots
\,i_\alpha\,\bar i_1\ldots\,\bar i_\nu\,h_1\ldots\,h_m}_{j_1\ldots\,
j_\beta\,\bar j_1\ldots\,\bar j_\gamma\,k_1\ldots\,k_n}[P]$ in 
\mythetag{13.17}. Being dependent on $q$, all these coefficients are 
some smooth complex-valued functions of the variables \mythetag{9.11} 
and \mythetag{9.12}. However, if $q$ is fixed, they all are complex
constants.\par
     For the fields $D(\hat{\boldsymbol\Psi}_i)$ in \mythetag{13.19} 
and for $D(\hat{\bPsi}_i)$ in \mythetag{13.20} we have the expansions
similar to the above expansion \mythetag{13.22}:
$$
\aligned
&\hskip -2em
D(\hat{\boldsymbol\Psi}_i)=\sum^2_{k=1}\Alpha^k_i\,\boldsymbol\Psi_k,\\
&\hskip -2em
D(\hat{\bPsi}_i)
=\sum^2_{k=1}\bar{\Alpha}\vphantom{\Alpha}^k_i\,\bPsi_k.
\endaligned
\mytag{13.23}
$$
In general case $\Alpha^j_i$ and $\bar{\Alpha}\vphantom{\Alpha}^j_i$ are
arbitrary complex-valued functions of the variables \mythetag{9.11} and
\mythetag{9.12}. They are complex constants if $q\in N$ is fixed.\par
     Let's return back to the formulas \mythetag{13.19}, \mythetag{13.20}, 
and \mythetag{13.21}. For the values of the fields
$D(\hat{\boldsymbol\vartheta}^{\,j})$,
$D(\hat{\overline{\boldsymbol\vartheta}}\vphantom{\vartheta}^{\,j})$,
$D(\hat{\boldsymbol\eta}^k)$ at the point $q$ we write
$$
\allowdisplaybreaks
\align
&\hskip -2em
D(\hat{\boldsymbol\eta}^k)
=-\sum^3_{h=0}\Gamma^k_h\ \boldsymbol\eta^{\,h},
\mytag{13.24}\\
&\hskip -2em
D(\hat{\boldsymbol\vartheta}^{\,j})
=-\sum^3_{h=0}\Alpha^j_h\ \boldsymbol\vartheta^{\,h},
\mytag{13.25}\\
&\hskip -2em
D(\hat{\overline{\boldsymbol\vartheta}}\vphantom{\vartheta}^{\,j})
=-\sum^3_{h=0}\bar{\Alpha}\vphantom{\Alpha}^j_h\ 
\overline{\boldsymbol\vartheta}\vphantom{\vartheta}^{\,h}.
\mytag{13.26}\\
\endalign
$$\par
     Let $\bold X$ be an extended spin-tensorial field of the type 
$(\varepsilon,\eta|\sigma,\zeta|e,f)$ given by the expansion \mythetag{13.4}. 
Then $D(\bold X)$ is given by the analogous
expansion:
$$
\gather
D(\bold X)=
\dsize\msums{2}{3}\Sb a_1,\,\ldots,\,a_\varepsilon\\
b_1,\,\ldots,\,b_\eta\\ \bar a_1,\,\ldots,\,\bar a_\sigma\\ 
\bar b_1,\,\ldots,\,\bar b_\zeta\\
c_1,\,\ldots,\,c_e\\ d_1,\,\ldots,\,d_f\endSb
DX^{a_1\ldots\,a_\varepsilon\bar a_1\ldots\,\bar a_\sigma
h_1\ldots\,h_e}_{b_1\ldots\,b_\eta\bar b_1\ldots\,\bar b_\zeta
k_1\ldots\,k_f}\ 
\boldsymbol\Psi^{b_1\ldots\,b_\eta\bar b_1\ldots\,\bar b_\zeta
d_1\ldots\,d_f}_{a_1\ldots\,a_\varepsilon\bar a_1\ldots\,\bar a_\sigma
c_1\ldots\,c_e}.\qquad\quad
\mytag{13.27}\\
\vspace{-0.9ex}
\intertext{where}\\
\vspace{-4ex}
\gathered
DX^{a_1\ldots\,a_\varepsilon\bar a_1\ldots\,\bar a_\sigma c_1
\ldots\,c_e}_{b_1\ldots\,b_\eta\bar b_1\ldots\,\bar b_\zeta
d_1\ldots\,d_f}
=\sum^3_{i=0}\sum^3_{j=0}Z^i\,\Upsilon^j_i\,
\frac{\partial X^{a_1\ldots\,a_\varepsilon\bar a_1\ldots\,
\bar a_\sigma c_1\ldots\,c_e}_{b_1\ldots\,b_\eta\bar b_1
\ldots\,\bar b_\zeta d_1\ldots\,
d_f}}{\partial x^j}\,+\\
+\sum^{J+Q}_{P=1}\dsize\msums{2}{3}\Sb i_1,\,\ldots,\,i_\alpha\\
j_1,\,\ldots,\,j_\beta\\ \bar i_1,\,\ldots,\,\bar i_\nu\\ 
\bar j_1,\,\ldots,\,\bar j_\gamma\\
h_1,\,\ldots,\,h_m\\ k_1,\,\ldots,\,k_n\endSb
Z^{i_1\ldots\,i_\alpha\,\bar i_1\ldots\,\bar i_\nu\,h_1\ldots
\,h_m}_{j_1\ldots\,j_\beta\,\bar j_1\ldots\,\bar j_\gamma\,k_1\ldots
\,k_n}[P]\ 
\frac{\partial X^{a_1\ldots\,a_\varepsilon\bar a_1\ldots\,
\bar a_\sigma c_1\ldots\,c_e}_{b_1\ldots\,b_\eta\bar b_1
\ldots\,\bar b_\zeta d_1\ldots\,
d_f}}{\partial S^{\,i_1\ldots\,i_\alpha\,\bar i_1\ldots\,\bar i_\nu\,
h_1\ldots\,h_m}_{j_1\ldots\,j_\beta\,\bar j_1\ldots\,\bar j_\gamma\,
k_1\ldots\,k_n}[P]}\,+\\
+\sum^{J+Q}_{P=1}\dsize\msums{2}{3}\Sb i_1,\,\ldots,\,i_\alpha\\
j_1,\,\ldots,\,j_\beta\\ \bar i_1,\,\ldots,\,\bar i_\nu\\ 
\bar j_1,\,\ldots,\,\bar j_\gamma\\
h_1,\,\ldots,\,h_m\\ k_1,\,\ldots,\,k_n\endSb
\bar Z^{\bar i_1\ldots\,\bar i_\nu\,i_1\ldots\,i_\alpha\,h_1\ldots
\,h_m}_{\bar j_1\ldots\,\bar j_\gamma\,j_1\ldots\,j_\beta\,k_1\ldots
\,k_n}[P]\ 
\frac{\partial X^{a_1\ldots\,a_\varepsilon\bar a_1\ldots\,
\bar a_\sigma c_1\ldots\,c_e}_{b_1\ldots\,b_\eta\bar b_1\ldots\,
\bar b_\zeta d_1\ldots\,
d_f}}{\partial\overline{S^{\,i_1\ldots\,i_\alpha\,\bar i_1\ldots\,
\bar i_\nu\,h_1\ldots\,h_m}_{j_1\ldots\,j_\beta\,\bar j_1\ldots\,
\bar j_\gamma\,k_1\ldots\,k_n}[P]}\,}\,+\\
\kern -9em
+\sum^\varepsilon_{\mu=1}\sum^2_{v_\mu=1}\Alpha^{a_\mu}_{v_\mu}\ 
X^{a_1\ldots\,v_\mu\,\ldots\,a_\varepsilon\bar a_1\ldots\,\bar a_\sigma
c_1\ldots\,c_e}_{b_1\ldots\,\ldots\,\ldots\,b_\eta\bar b_1\ldots\,
\bar b_\zeta d_1\ldots\,d_f}\,-\\
\kern 9em-\sum^\eta_{\mu=1}\sum^2_{w_\mu=1}\Alpha^{w_\mu}_{b_\mu}\
X^{a_1\ldots\,\ldots\,\ldots\,a_\varepsilon\bar a_1\ldots\,\bar a_\sigma
c_1\ldots\,c_e}_{b_1\ldots\,w_\mu\,\ldots\,b_\eta\bar b_1\ldots\,
\bar b_\zeta d_1\ldots\,d_f}\,+\\
\kern -9em
+\sum^\sigma_{\mu=1}\sum^2_{v_\mu=1}
\bar{\Alpha}\vphantom{\Alpha}^{\bar a_\mu}_{v_\mu}\ 
X^{a_1\ldots\,a_\varepsilon\bar a_1\ldots\,v_\mu\,\ldots\,\bar a_\sigma
c_1\ldots\,c_e}_{b_1\ldots\,b_\eta\bar b_1\ldots\,\ldots\,\ldots\,
\bar b_\zeta d_1\ldots\,d_f}\,-\\
\kern 9em-\sum^\zeta_{\mu=1}\sum^2_{w_\mu=1}
\bar{\Alpha}\vphantom{\Alpha}^{w_\mu}_{\bar b_\mu}\
X^{a_1\ldots\,a_\varepsilon\bar a_1\ldots\,\ldots\,\ldots\,\bar a_\sigma
c_1\ldots\,c_e}_{b_1\ldots\,b_\eta\bar b_1\ldots\,w_\mu\,\ldots\,
\bar b_\zeta d_1\ldots\,d_f}\,+\\
\kern -9em+\sum^e_{\mu=1}\sum^3_{v_\mu=0}\Gamma^{c_\mu}_{v_\mu}\ 
X^{a_1\ldots\,a_\varepsilon\bar a_1\ldots\,\bar a_\sigma
c_1\ldots\,v_\mu\,\ldots\,c_e}_{b_1\ldots\,b_\eta\bar b_1\ldots\,
\bar b_\zeta d_1\ldots\,\ldots\,\ldots\,d_f}\,-\\
\kern 9em-\sum^f_{\mu=1}\sum^3_{w_\mu=0}\Gamma^{w_\mu}_{b_\mu}\
X^{a_1\ldots\,a_\varepsilon\bar a_1\ldots\,\bar a_\sigma
c_1\ldots\,\ldots\,\ldots\,c_e}_{b_1\ldots\,b_\eta\bar b_1\ldots\,
\bar b_\zeta d_1\ldots\,w_\mu\,\ldots\,d_f}.
\endgathered\qquad\quad
\mytag{13.28}
\endgather
$$
\mylemma{13.6} Any differentiation $D$ of the algebra of extended 
spin-tensorial fields $\bold S(M)$ is uniquely fixed by its restrictions 
to the modules $S^0_0\bar S^0_0 T^0_0(M)$, $S^1_0\bar S^0_0 T^0_0(M)$,
$S^0_0\bar S^1_0 T^0_0(M)$, $S^0_0\bar S^0_0 T^1_0(M)$ in the direct
sum \mythetag{11.1}.
\endproclaim
    This lemma~\mythelemma{13.6} strengthens the previous
lemma~\mythelemma{13.5}. It is proved by direct calculations when we 
apply \mythetag{13.24}, \mythetag{13.25}, \mythetag{13.26} to
\mythetag{13.19}, \mythetag{13.20}, and \mythetag{13.21}. Another
proof can be based on \mythetag{13.27} and on the general explicit formula
\mythetag{13.28} for the components of the field $D(\bold X)$. Indeed, 
we see that any differentiation $D$ is completely determined by the 
parameters $Z^i$, $Z^{i_1\ldots\,i_\alpha\,\bar i_1\ldots\,\bar i_\nu
\,h_1\ldots\,h_m}_{j_1\ldots\,j_\beta\,\bar j_1\ldots\,\bar j_\gamma\,k_1
\ldots\,k_n}[P]$, $\bar Z^{\bar i_1\ldots\,\bar i_\nu\,i_1\ldots\,i_\alpha
\,h_1\ldots\,h_m}_{\bar j_1\ldots\,\bar j_\gamma\,j_1\ldots\,j_\beta\,k_1
\ldots\,k_n}[P]$ characterizing the restriction of $D$ to 
$S^0_0\bar S^0_0 T^0_0(M)$ and by the parameters $\Alpha^k_i$,
$\bar{\Alpha}\vphantom{\Alpha}^k_i$, $\Gamma^k_h$ describing the
restrictions of $D$ to $S^1_0\bar S^0_0 T^0_0(M)$, $S^0_0\bar S^1_0
T^0_0(M)$, and $S^0_0\bar S^0_0 T^1_0(M)$ respectively.\par 
     $Z$-parameters, $\Alpha$-parameters and $\Gamma$-parameters
determining a differentiation $D$ in \mythetag{13.28} obey some 
definite transformation rules under a change of a local chart. 
The transformation rule for $Z^i$ is the most simple one:
$$
\hskip -2em
Z^i=\dsize\sum^3_{j=0}S^i_j\,\tilde Z^j.
\mytag{13.29}
$$
The transformation rule for $Z^{i_1\ldots\,i_\alpha\,\bar i_1\ldots\,
\bar i_\nu\,h_1\ldots\,h_m}_{j_1\ldots\,j_\beta\,\bar j_1\ldots\,
\bar j_\gamma\,k_1\ldots\,k_n}[P]$ is much more huge:
$$
\gathered
Z^{i_1\ldots\,i_\alpha\,\bar i_1\ldots\,\bar i_\nu\,h_1
\ldots\,h_m}_{j_1\ldots\,j_\beta\,\bar j_1\ldots\,\bar j_\gamma
\,k_1\ldots\,k_n}[P]=
\dsize\msum{2}\Sb a_1,\,\ldots,\,a_\alpha\\ b_1,\,\ldots,\,b_\beta\\
\bar a_1,\,\ldots,\,\bar a_\nu\\ \bar b_1,\,\ldots,\,\bar b_\gamma\endSb
\dsize\msum{3}\Sb c_1,\,\ldots,\,c_m\\ d_1,\,\ldots,\,d_n\endSb
\goth S^{\,i_1}_{a_1}\ldots\,\goth S^{\,i_\alpha}_{a_\alpha}\times\\
\times\,\goth T^{b_1}_{j_1}\ldots\,\goth T^{b_\beta}_{j_\beta}
\ \overline{\goth S^{\,\bar i_1}_{\bar a_1}}
\ldots\,\overline{\goth S^{\,\bar i_\nu}_{\bar a_\nu}}\ 
\ \overline{\goth T^{\,\bar b_1}_{\bar j_1}}
\ldots\,\overline{\goth T^{\,\bar b_\gamma}_{\bar j_\gamma}}\ 
S^{h_1}_{c_1}\ldots\,S^{h_m}_{c_m}\,
T^{\,d_1}_{k_1}\ldots\,T^{\,d_n}_{k_n}\,\times\\
\times\,\tilde Z^{a_1\ldots\,a_\alpha\,\bar a_1\ldots\,\bar a_\nu\,c_1
\ldots\,c_m}_{b_1\ldots\,b_\beta\,\bar b_1\ldots\,\bar b_\gamma
\,d_1\ldots\,d_n}[P]-\sum^\alpha_{\mu=1}\sum^3_{i=0}\sum^3_{j=0}
\sum^2_{v_\mu=1}\vartheta^{\,i_\mu}_{i\,v_\mu}\,\times\\
\times\ S^{\,i_1\ldots\,v_\mu\ldots\,i_\alpha\,\bar i_1\ldots\,\bar i_\nu\,
h_1\ldots\,h_m}_{\,j_1\ldots\,\ldots\,\ldots\,j_\beta\,\bar j_1\ldots\,
\bar j_\gamma\,k_1\ldots\,k_n}[P]\ S^i_j\ \tilde Z^j
+\sum^\beta_{\mu=1}\sum^3_{i=0}\sum^3_{j=0}\sum^2_{w_\mu=1}
\vartheta^{\,w_\mu}_{i\,j_\mu}\,\times\\
\times\ S^{\,i_1\ldots\,\ldots\,\ldots\,i_\alpha\,\bar i_1\ldots\,
\bar i_\nu\,h_1\ldots\,h_m}_{\,j_1\ldots\,w_\mu\,\ldots\,j_\beta\,
\bar j_1\ldots\,\bar j_\gamma\,k_1\ldots\,k_n}[P]\ S^i_j\ \tilde Z^j
-\sum^\nu_{\mu=1}\sum^3_{i=0}\sum^3_{j=0}\sum^2_{v_\mu=1}
\overline{\vartheta^{\,\bar i_\mu}_{i\,v_\mu}}\,\times\\
\times\ S^{\,i_1\ldots\,i_\alpha\,\bar i_1\ldots\,v_\mu\ldots\,\bar i_\nu\,
h_1\ldots\,h_m}_{\,j_1\ldots\,j_\beta\,\bar j_1\ldots\,\ldots\,\ldots\,
\bar j_\gamma\,k_1\ldots\,k_n}[P]\ S^i_j\ \tilde Z^j
+\sum^\gamma_{\mu=1}\sum^3_{i=0}\sum^3_{j=0}\sum^2_{w_\mu=1}
\overline{\vartheta^{\,w_\mu}_{i\,\bar j_\mu}}\,\times\\
\times\ S^{\,i_1\ldots\,i_\alpha\,\bar i_1\ldots\,\ldots\,\ldots\,
\bar i_\nu\,h_1\ldots\,h_m}_{\,j_1\ldots\,j_\beta\,\bar j_1\ldots\,
w_\mu\,\ldots\,\bar j_\gamma\,k_1\ldots\,k_n}[P]\ S^i_j\ \tilde Z^j
-\sum^m_{\mu=1}\sum^3_{i=0}\sum^3_{j=0}\sum^3_{v_\mu=0}
\theta^{\,i_\mu}_{i\,v_\mu}\,\times\\
\times\ S^{\,i_1\ldots\,i_\alpha\,\bar i_1\ldots\,\bar i_\nu\,
h_1\ldots\,v_\mu\ldots\,h_m}_{\,j_1\ldots\,j_\beta\,\bar j_1\ldots\,
\bar j_\gamma\,k_1\ldots\,\ldots\,\ldots\,k_n}[P]\ S^i_j\ \tilde Z^j
+\sum^n_{\mu=1}\sum^3_{i=0}\sum^3_{j=0}\sum^3_{w_\mu=0}
\theta^{\,w_\mu}_{i\,j_\mu}\,\times\\
\times\ S^{\,i_1\ldots\,i_\alpha\,\bar i_1\ldots\,\bar i_\nu\,
h_1\ldots\,\ldots\,\ldots\,h_m}_{\,j_1\ldots\,j_\beta\,\bar j_1\ldots\,
\bar j_\gamma\,k_1\ldots\,w_\mu\,\ldots\,k_n}[P]\ S^i_j\ \tilde Z^j.
\vphantom{\sum^n_{\mu=1}}
\endgathered\quad
\mytag{13.30}
$$
The transformation rule for $\bar Z^{\bar i_1\ldots\,\bar i_\nu\,
i_1\ldots\,i_\alpha\,h_1\ldots\,h_m}_{\bar j_1\ldots\,\bar j_\gamma\,
j_1\ldots\,j_\beta\,k_1\ldots\,k_n}[P]$ is equally huge as the previous
transformation rule for $Z^{i_1\ldots\,i_\alpha\,\bar i_1\ldots\,
\bar i_\nu\,h_1\ldots\,h_m}_{j_1\ldots\,j_\beta\,\bar j_1\ldots\,
\bar j_\gamma\,k_1\ldots\,k_n}[P]$:
$$
\gathered
\bar Z^{\bar i_1\ldots\,\bar i_\nu\,
i_1\ldots\,i_\alpha\,h_1\ldots\,h_m}_{\bar j_1\ldots\,\bar j_\gamma\,
j_1\ldots\,j_\beta\,k_1\ldots\,k_n}[P]=
\dsize\msum{2}\Sb a_1,\,\ldots,\,a_\alpha\\ b_1,\,\ldots,\,b_\beta\\
\bar a_1,\,\ldots,\,\bar a_\nu\\ \bar b_1,\,\ldots,\,\bar b_\gamma\endSb
\dsize\msum{3}\Sb c_1,\,\ldots,\,c_m\\ d_1,\,\ldots,\,d_n\endSb
\goth S^{\,\bar i_1}_{\bar a_1}\ldots\,
\goth S^{\,\bar i_\nu}_{\bar a_\nu}\times\\
\times\,\goth T^{\bar b_1}_{\bar j_1}\ldots\,
\goth T^{\bar b_\gamma}_{\bar j_\gamma}
\ \overline{\goth S^{\,i_1}_{a_1}}\ldots\,
\overline{\goth S^{\,i_\alpha}_{a_\alpha}}\ 
\ \overline{\goth T^{\,b_1}_{j_1}}
\ldots\,\overline{\goth T^{\,b_\beta}_{j_\beta}}\ 
S^{h_1}_{c_1}\ldots\,S^{h_m}_{c_m}\,
T^{\,d_1}_{k_1}\ldots\,T^{\,d_n}_{k_n}\,\times\\
\times\,\tbZ^{\bar a_1\ldots\,\bar a_\nu\,a_1\ldots\,a_\alpha\,c_1
\ldots\,c_m}_{\bar b_1\ldots\,\bar b_\gamma\,b_1\ldots\,b_\beta\,
d_1\ldots\,d_n}[P]-\sum^\alpha_{\mu=1}\sum^3_{i=0}\sum^3_{j=0}
\sum^2_{v_\mu=1}\overline{\vartheta^{\,i_\mu}_{i\,v_\mu}}\,\times\\
\times\ \overline{S^{\,i_1\ldots\,v_\mu\ldots\,i_\alpha\,\bar i_1
\ldots\,\bar i_\nu\,h_1\ldots\,h_m}_{\,j_1\ldots\,\ldots\,\ldots\,
j_\beta\,\bar j_1\ldots\,\bar j_\gamma\,k_1\ldots\,k_n}[P]}\ S^i_j
\ \tilde Z^j+\sum^\beta_{\mu=1}\sum^3_{i=0}\sum^3_{j=0}\sum^2_{w_\mu=1}
\overline{\vartheta^{\,w_\mu}_{i\,j_\mu}}\,\times\\
\times\ \overline{S^{\,i_1\ldots\,\ldots\,\ldots\,i_\alpha\,\bar i_1
\ldots\,\bar i_\nu\,h_1\ldots\,h_m}_{\,j_1\ldots\,w_\mu\,\ldots\,j_\beta
\,\bar j_1\ldots\,\bar j_\gamma\,k_1\ldots\,k_n}[P]}\ S^i_j\ \tilde Z^j
-\sum^\nu_{\mu=1}\sum^3_{i=0}\sum^3_{j=0}\sum^2_{v_\mu=1}
\vartheta^{\,\bar i_\mu}_{i\,v_\mu}\,\times\\
\times\ \overline{S^{\,i_1\ldots\,i_\alpha\,\bar i_1\ldots\,v_\mu
\ldots\,\bar i_\nu\,h_1\ldots\,h_m}_{\,j_1\ldots\,j_\beta\,\bar j_1
\ldots\,\ldots\,\ldots\,\bar j_\gamma\,k_1\ldots\,k_n}[P]}\ S^i_j\ 
\tilde Z^j+\sum^\gamma_{\mu=1}\sum^3_{i=0}\sum^3_{j=0}\sum^2_{w_\mu=1}
\vartheta^{\,w_\mu}_{i\,\bar j_\mu}\,\times\\
\times\ \overline{S^{\,i_1\ldots\,i_\alpha\,\bar i_1\ldots\,\ldots\,
\ldots\,\bar i_\nu\,h_1\ldots\,h_m}_{\,j_1\ldots\,j_\beta\,\bar j_1
\ldots\,w_\mu\,\ldots\,\bar j_\gamma\,k_1\ldots\,k_n}[P]}\ S^i_j
\ \tilde Z^j-\sum^m_{\mu=1}\sum^3_{i=0}\sum^3_{j=0}\sum^3_{v_\mu=0}
\theta^{\,i_\mu}_{i\,v_\mu}\,\times\\
\times\ \overline{S^{\,i_1\ldots\,i_\alpha\,\bar i_1\ldots\,\bar i_\nu
\,h_1\ldots\,v_\mu\ldots\,h_m}_{\,j_1\ldots\,j_\beta\,\bar j_1\ldots\,
\bar j_\gamma\,k_1\ldots\,\ldots\,\ldots\,k_n}[P]}\ S^i_j\ \tilde Z^j
+\sum^n_{\mu=1}\sum^3_{i=0}\sum^3_{j=0}\sum^3_{w_\mu=0}
\theta^{\,w_\mu}_{i\,j_\mu}\,\times\\
\times\ \overline{S^{\,i_1\ldots\,i_\alpha\,\bar i_1\ldots\,\bar i_\nu
\,h_1\ldots\,\ldots\,\ldots\,h_m}_{\,j_1\ldots\,j_\beta\,\bar j_1\ldots
\,\bar j_\gamma\,k_1\ldots\,w_\mu\,\ldots\,k_n}[P]}\ S^i_j\ \tilde Z^j.
\vphantom{\sum^n_{\mu=1}}
\endgathered\quad
\mytag{13.31}
$$
The transformation rules \mythetag{13.29}, \mythetag{13.30}, 
\mythetag{13.31} are completed with the transformation rules 
for $\Alpha^k_i$, $\bar{\Alpha}\vphantom{\Alpha}^k_i$ and
$\Gamma^k_h$. They are less huge:
$$
\gather
\hskip -2em
\Alpha^k_i=\sum^2_{b=1}\sum^2_{a=1}\goth S^k_a\,\goth T^b_i\ 
\tilde{\Alpha}\vphantom{\Alpha}^a_b
+\sum^3_{a=0}Z^a\,\vartheta^k_{ai},
\mytag{13.32}\\
\vspace{1ex}
\hskip -2em
\bar{\Alpha}\vphantom{\Alpha}^k_i=\sum^2_{b=1}\sum^2_{a=1}
\overline{\goth S^k_a}\ \overline{\goth T^b_i}\ 
\tilde{\bar{\Alpha}}\vphantom{\Alpha}^a_b
+\sum^3_{a=0}Z^a\,\overline{\vartheta^k_{ai}},
\mytag{13.33}\\
\vspace{1ex}
\hskip -2em
\Gamma^k_h=\sum^3_{b=0}\sum^3_{a=0}S^k_a\,T^b_h\ \tilde\Gamma^a_b
+\sum^3_{a=0}Z^a\,\theta^k_{ah}.
\mytag{13.34}
\endgather
$$
These transformation rules \mythetag{13.29}, \mythetag{13.30},
\mythetag{13.31} are derived from \mythetag{13.2} using \mythetag{9.21},
\mythetag{9.23}, and \mythetag{9.37}. The transformation rules 
\mythetag{13.32}, \mythetag{13.33}, and \mythetag{13.34} are derived from
\mythetag{13.23} and \mythetag{13.22} using the formulas \mythetag{9.17},
\mythetag{9.16}, \mythetag{9.29}, \mythetag{9.28}, \mythetag{9.33}, 
and \mythetag{9.35}.\par 
     Note that a differentiation $D$ acts as a first order linear
differential operator upon the components of an extended spin-tensorial
field $\bold X$ in the formula \mythetag{13.28}. The coefficients $Z^i$, 
$Z^{h_1\ldots\,h_r}_{k_1\ldots\,k_s}[P]$, $Z^{i_1\ldots\,i_\alpha\,
\bar i_1\ldots\,\bar i_\nu\,h_1\ldots\,h_m}_{j_1\ldots\,j_\beta\,
\bar j_1\ldots\,\bar j_\gamma\,k_1\ldots\,k_n}[P]$, $\bar Z^{\bar i_1
\ldots\,\bar i_\nu\,i_1\ldots\,i_\alpha\,h_1\ldots\,h_m}_{\bar j_1\ldots
\,\bar j_\gamma\,j_1\ldots\,j_\beta\,k_1\ldots\,k_n}[P]$, $\Alpha^k_i$,
$\bar{\Alpha}\vphantom{\alpha}^k_i$, $\Gamma^k_i$ of the linear operator 
in this formula are not differentiated. Therefore, fixing some point 
$q\in N$ and taking the values of these coefficients at the point
$q$, we can say that we know this linear operator at that particular point
$q$ even if we don't know the values of theese coefficients at other points.
\mydefinition{13.1} Let $N$ be a composite spin-tensorial bundle over the
space-time manifold $M$ (in the sense of the formula \mythetag{9.7}).
A {\it spin-tensorial first order differential operator\/} $D_q$ at a 
point $q\in N$ is a geometric object associated with the point $q$ and
represented by the set of constants $Z^i$, $Z^{i_1\ldots\,i_\alpha\,
\bar i_1\ldots\,\bar i_\nu
\,h_1\ldots\,h_m}_{j_1\ldots\,j_\beta\,\bar j_1\ldots\,\bar j_\gamma
\,k_1\ldots\,k_n}[P]$, $\bar Z^{\bar i_1\ldots\,\bar i_\nu\,i_1\ldots
\,i_\alpha\,h_1\ldots\,h_m}_{\bar j_1\ldots\,\bar j_\gamma\,j_1\ldots
\,j_\beta\,k_1\ldots\,k_n}[P]$, $\Gamma^k_i$, $\Alpha^k_i$,
$\bar{\Alpha}\vphantom{\alpha}^k_i$ in a proper local chart such that 
they obey the transformation rules \mythetag{13.29}, \mythetag{13.30},
\mythetag{13.31}, \mythetag{13.32}, \mythetag{13.33}, \mythetag{13.34}.
\enddefinition
     Spin-tensorial first order differential operators at a fixed
point $q$ constitute a finite-dimensional linear space over the field of
complex numbers $\Bbb C$, we denote it with the symbol
$\goth D(q,M)$. Its dimension is given by the formula
$$
\hskip -2em
\dim_{\,\Bbb C}\goth D(q,M)=\dim_{\,\Bbb R}N+4^2+2^2+2^2.
\mytag{13.35}
$$
The dimension of the composite spin-tensorial bundle $N$ in \mythetag{13.35}
is given by the formula \mythetag{9.10}. The linear spaces $\goth D(q,M)$ 
with $q$ running over $N$ are glued into a complex vector bundle over $N$. 
\mytheorem{13.1} Any differentiation $D$ of the algebra of extended 
spin-tensorial fields $\bold S(M)$ is represented as a field of
differential operators $D_q\in\goth D(q,M)$, one per each 
point $q\in N$. Conversely, each smooth field of spin-tensorial first 
order differential operators is a differentiation of the algebra 
$\bold S(M)$.
\endproclaim
The theorem~\mythetheorem{13.1} solves the problem of localization 
announced in the very beginning of this section. 
\head
14. Degenerate differentiations.
\endhead
\mydefinition{14.1} A differentiation $D$ of the algebra of extended 
spin-tensorial fields $\bold S(M)$ is called a {\it degenerate
differentiation\/} if its restriction \mythetag{13.1} to the module
$S^0_0\bar S^0_0 T^0_0(M)$ is identically zero.
\enddefinition
    For a degenerate differentiation $D$ from the formula \mythetag{13.2}
we derive that all its $Z$-components are identically equal to zero:
$$
\gather
\hskip -2em
Z^i=0,
\mytag{14.1}\\
\vspace{1.5ex}
\hskip -2em
Z^{i_1\ldots\,i_\alpha\,\bar i_1\ldots\,\bar i_\nu\,h_1\ldots
\,h_m}_{j_1\ldots\,j_\beta\,\bar j_1\ldots\,\bar j_\gamma\,k_1\ldots
\,k_n}[P]=0,
\mytag{14.2}\\
\vspace{1.5ex}
\hskip -2em
\bar Z^{\bar i_1\ldots\,\bar i_\nu\,i_1\ldots\,i_\alpha\,h_1\ldots
\,h_m}_{\bar j_1\ldots\,\bar j_\gamma\,j_1\ldots\,j_\beta\,k_1\ldots
\,k_n}[P]=0.
\mytag{14.3}
\endgather
$$
Unlike $Z$-components, $\Alpha$-components and $\Gamma$-components of 
such a differentiation $D$ in general case are not zero. Substituting
\mythetag{14.1} into the formulas \mythetag{13.32}, \mythetag{13.33}, 
and \mythetag{13.33}, we find that the transformation rules for 
$\Alpha$-components and $\Gamma$-components of a degenerate differentiation
$D$ reduce to the following ones:
$$
\gather
\hskip -2em
\Alpha^k_i=\sum^2_{b=1}\sum^2_{a=1}\goth S^k_a\,\goth T^b_i\ 
\tilde{\Alpha}\vphantom{\Alpha}^a_b,
\mytag{14.4}\\
\vspace{1ex}
\hskip -2em
\bar{\Alpha}\vphantom{\Alpha}^k_i=\sum^2_{b=1}\sum^2_{a=1}
\overline{\goth S^k_a}\ \overline{\goth T^b_i}\ 
\tilde{\bar{\Alpha}}\vphantom{\Alpha}^a_b.
\mytag{14.5}\\
\hskip -2em
\Gamma^k_h=\sum^3_{b=0}\sum^3_{a=0}S^k_a\,T^b_h\ \tilde\Gamma^a_b.
\mytag{14.6}
\endgather
$$
Note that these transformation rules \mythetag{14.4}, \mythetag{14.5},
\mythetag{14.6} are special cases of the transformation rule \mythetag{10.3} 
for the components of an extended spin-tensorial field,
while \mythetag{10.2} is inverse to \mythetag{10.3}. This observation
proves the following theorem.
\mytheorem{14.1} Defining a degenerate differentiation $D$ of the algebra
$\bold S(M)$ is equivalent to defining three extended spin-tensorial fields
$\eufb S$, $\bar{\eufb S}$, and $\bold S$ of the types $(1,1|0,0|0,0)$,
$(0,0|1,1|0,0)$, and $(0,0|0,0|1,1)$ respectively.
\endproclaim
     By tradition we use the symbols $\eufb S$, $\bar{\eufb S}$, and 
$\bold S$ for the extended spin-tensorial fields determined by the
$\Alpha$ and the $\Gamma$-components of a degenerate differentiation $D$.
One should be careful for not to confuse the components of the field 
$\bold S$ and the components of the direct transition matrix $S$ (e\.\,g\.
in the formula \mythetag{14.6}).\par
\head
15. Covariant differentiations.
\endhead
     The set of differentiations of the algebra of extended spin-tensorial
fields $\bold S(M)$ possesses the structure of a module over the ring of
smooth complex functions $\goth F_{\Bbb C}(N)$. The set of complexified 
extended vector fields $\Bbb CT^1_0(M)=S^0_0\bar S^0_0T^1_0(M)$ (see 
\mythetag{10.1}) also is a module over the same ring $\goth F_{\Bbb C}(N)$.
Therefore, the following definition is consistent.
\mydefinition{15.1} Say that in the space-time manifold $M$ a covariant
differentiation of the algebra of extended spin-tensorial fields $\bold
S(M)$ is given if some homomorphism of $\goth F_{\Bbb C}(N)$-modules
$\nabla\!:\,\Bbb CT^1_0(M)\to\goth D_{\bold S}(M)$ is given. The image 
of a vector field $\bold Y$ under such homomorphism is denoted by
$\nabla_{\bold Y}$. The differentiation $D=\nabla_{\bold Y}$
is called the {\it covariant differentiation along the vector field\/}
$\bold Y$.
\enddefinition
     Let's note that the $\goth F_{\Bbb C}(N)$-module $\goth D_{\bold
S}(M)$ admits a localization in the sense of the following definition
(compare with the definition~8.2 in \mycite{5}). 
\mydefinition{15.2} Let $A$ be a module over the ring of smooth 
complex functions $\goth F_{\Bbb C}(M)$ in some smooth real 
manifold $M$. We say that the module $A$ admits a {\it localization\/} 
if it is isomorphic to a functional module so that each element $\bold 
a\in A$ is represented as some function $\bold a(q)=\bold a_q$ in $M$
taking its values in some $\Bbb C$-linear spaces $A_q$ associated with 
each point $q$ of the manifold $M$.
\enddefinition
\noindent Indeed, according to the theorem~\mythetheorem{13.1} each 
differentiation $D$ is a field of differential operators. As for the
differential operators themselves, they form finite-dimensional 
$\Bbb C$-linear spaces $\goth D(q,M)$ of the dimension \mythetag{13.35}, 
\pagebreak one per each point $q\in N$.\par
     The definition~\mythedefinition{15.1} is a complexified version of
the definition~8.2 from \mycite{5}. One can easily formulate complexified
versions for the definition~8.3 and for the theorems~8.2 and 8.3 from
\mycite{5}.
\mydefinition{15.3} Let $A$ be a module that admits a localization 
in the sense of the definition~\mythedefinition{15.2}. We say that
the localization of $A$ is a {\it complete localization} if the 
following two conditions are fulfilled:
\roster
\rosteritemwd=5pt
\item for any point $q\in M$ and for any vector $\bold v\in A_q$ 
      there exists an element $\bold a\in A$ such that $\bold a_q
      =\bold v$;
\item if $\bold a_q=0$ at some point $q\in M$, then there exist 
some finite set of elements $\bold E_0,\,\ldots,\,\bold E_n$ in
$A$ and some smooth complex functions $\alpha_0\,\,\dots,\,\alpha_n$ 
vanishing at the point $q$ such that $\bold a=\alpha_0\,\bold E_0
+\ldots+\alpha_n\,\bold E_n$.
\endroster      
\enddefinition
\mytheorem{15.1} Let $A$ and $B$ be two $\goth F_{\Bbb C}(M)$-modules 
that admit localizations. If the localization of $A$ is a complete
localization, then each homomorphism $f\!:\,A\to B$ is represented 
by a family of\/ $\Bbb C$-linear mappings 
$$
\hskip -2em
F_q\!:\,A_q\to B_q
\mytag{15.1}
$$
so that if\/ $\bold a\in A$ and $\bold b=f(\bold a)$, then $\bold b_q
=F_q(\bold a_q)$ for each point $q\in M$.
\endproclaim
\mytheorem{15.2} Let $\pi\!:\,VM\to M$ be a smooth $n$-dimensional complex
vector bundle over some smooth real base manifold $M$ and let $A$ be the
$\goth F_{\Bbb C}(M)$-module of all global smooth sections\footnotemark\ of
this bundle. Then $A$ admits a complete localization in the sense of the 
definition~\mythedefinition{15.3}.
\endproclaim
\footnotetext{ \ See the definition of {\it sections\/} and {\it smooth
sections\/} in \mycite{10}.}
\adjustfootnotemark{-1}
     The proof of the theorems~\mythetheorem{15.1} and 
\mythetheorem{15.2} is analogous to the proof of the theorems~8.1 
and 8.2 in \mycite{5}. We leave this proof to the reader.\par
     Now let's return back to the spaces $\goth D(q,M)$ with $q$ running
over $N$. Remember that they are glued into a complex vector bundle for
which $N$ is a base manifold. Similarly, the complexified tangent spaces
$\Bbb CT_{\pi(q)}(M)=\Bbb CT^1_0(\pi(q),M)$ are also glued into a complex
vector bundle over $N$. This bundle is called an {\it induced bundle}. It
is induced by the projection map \mythetag{9.9} from the complexified
tangent bundle $\Bbb CT\!M$.\par
    Note that the $\goth F_{\Bbb C}(N)$ modules $\Bbb CT^1_0(M)$ and 
$\goth D_{\bold S}(M)$ in the definition~\mythedefinition{15.1} are
represented as the sets of all smooth global sections for the above
two bundles. For this reason we can apply the theorem~\mythetheorem{15.2}
to $\Bbb CT^1_0(M)$ and $\goth D_{\bold S}(M)$. As a result we conclude 
that each covariant differentiation $\nabla$ of the algebra of extended 
spin-tensorial fields $\bold S(M)$ is composed by $\Bbb C$-linear maps 
$\Bbb CT_{\pi(q)}(M)\to\goth D(q,M)$ specific to each point $q\in N$. 
This fact is expressed by the following formula:
$$
\hskip -2em
\nabla_{\bold Y}\bold X=C(\bold Y\otimes\nabla\bold X).
\mytag{15.2}
$$
If $\bold Y$ in \mythetag{15.2} is given by the expansion 
$$
\bold Y=\sum^3_{j=0}Y^j\ \boldsymbol\Upsilon_j 
$$
in some positively polarized right orthonormal frame $\boldsymbol\Upsilon_0,
\,\boldsymbol\Upsilon_1,\,\boldsymbol\Upsilon_2,\,\boldsymbol\Upsilon_3$
and if $\bold X$ is an extended spin-tensorial field of the type
$(\varepsilon,\eta|\sigma,\zeta|e,f)$ given by the expansion
\mythetag{13.4}, then the formula \mythetag{15.2} is specified to the 
following one:
$$
\nabla_{\bold Y}\bold X
=\sum^3_{j=0}\dsize\msums{2}{3}\Sb a_1,\,\ldots,\,a_\varepsilon\\
b_1,\,\ldots,\,b_\eta\\ \bar a_1,\,\ldots,\,\bar a_\sigma\\ 
\bar b_1,\,\ldots,\,\bar b_\zeta\\
c_1,\,\ldots,\,c_e\\ d_1,\,\ldots,\,d_f\endSb
Y^j\ \nabla_{\!\!j}X^{a_1\ldots\,a_\varepsilon\bar a_1\ldots\,
\bar a_\sigma c_1\ldots\,c_e}_{b_1\ldots\,b_\eta\bar b_1\ldots\,
\bar b_\zeta d_1\ldots\,d_f}\ 
\boldsymbol\Psi^{b_1\ldots\,b_\eta\bar b_1\ldots\,\bar b_\zeta
d_1\ldots\,d_f}_{a_1\ldots\,a_\varepsilon\bar a_1\ldots\,
\bar a_\sigma c_1\ldots\,c_e}.\quad
\mytag{15.3}
$$
Looking at the formula \mythetag{15.2}, we see that each covariant
differentiation $\nabla$ can be treated as an operator producing 
the extended spin-tensorial field $\nabla\bold X$ of the type
$(\varepsilon,\eta|\sigma,\zeta|e,f+1)$ from any given extended 
spin-tensorial field $\bold X$ of the type
$(\varepsilon,\eta|\sigma,\zeta|e,f)$. This operator increases by 
one the number of covariant tensorial indices of a spin-tensorial
field $\bold X$. It is called the operator of {\it covariant
differential\/} associated with the covariant differentiation 
$\nabla$. By $\nabla_{\!\!j}X^{a_1\ldots\,a_\varepsilon\bar a_1\ldots
\,\bar a_\sigma c_1\ldots\,c_e}_{b_1\ldots\,b_\eta\bar b_1\ldots\,
\bar b_\zeta d_1\ldots\,d_f}$ in \mythetag{15.3} we denote the 
components of the field $\nabla\bold X$. Comparing \mythetag{13.27}
and \mythetag{15.3}, we can write the following equality
for the differentiation $D=\nabla_{\bold Y}$:
$$
\hskip -2em
DX^{a_1\ldots\,a_\varepsilon\bar a_1\ldots\,\bar a_\sigma h_1
\ldots\,h_e}_{b_1\ldots\,b_\eta\bar b_1\ldots\,\bar b_\zeta k_1
\ldots\,k_f}
=\sum^3_{j=0}Y^j\ \nabla_{\!\!j}X^{a_1\ldots\,a_\varepsilon\bar a_1
\ldots\,\bar a_\sigma c_1\ldots\,c_e}_{b_1\ldots\,b_\eta\bar b_1\ldots
\,\bar b_\zeta d_1\ldots\,d_f}.
\mytag{15.4}
$$
Note that $\nabla_{\!\!j}$ in \mythetag{15.4} is the symbol of a covariant
derivative. Unlike $\nabla$ and $\nabla_{\bold Y}$, covariant derivatives
are applied not to spin-tensorial fields, but to its components in 
expansions like \mythetag{13.4}.\par
     Let's consider the linear map $\Bbb CT_{\pi(q)}(M)\to
\goth D(q,M)$ produced by some covariant differentiation $\nabla$ 
at some particular point $q\in N$. In a local chart this map 
is given by some linear functions expressing the components of the
differential operator $D_q$, where $D=\nabla_{\bold Y}$, through 
the components of the vector $\bold Y_{\!q}$:
$$
\allowdisplaybreaks
\align
&\hskip -2em
Z^{\,i}=\sum^3_{j=0}Z^{\,i}_j\ Y^j,
\mytag{15.5}\\
&\hskip -2em
Z^{\,i_1\ldots\,i_\alpha\,\bar i_1\ldots\,\bar i_\nu\,h_1\ldots
\,h_m}_{j_1\ldots\,j_\beta\,\bar j_1\ldots\,\bar j_\gamma\,k_1
\ldots\,k_n}[P]=\sum^3_{j=0}Z^{\,i_1\ldots\,i_\alpha\,\bar i_1\ldots
\,\bar i_\nu\,h_1\ldots\,h_m}_{j\,j_1\ldots\,j_\beta\,\bar j_1\ldots
\,\bar j_\gamma\,k_1\ldots\,k_n}[P]\ Y^j,
\mytag{15.6}\\
&\hskip -2em
\bar Z^{\,\bar i_1\ldots\,\bar i_\nu\,i_1\ldots\,i_\alpha\,h_1\ldots
\,h_m}_{\bar j_1\ldots\,\bar j_\gamma\,j_1\ldots\,j_\beta\,k_1\ldots
\,k_n}[P]=\sum^3_{j=0}\bar Z^{\,\bar i_1\ldots\,\bar i_\nu\,i_1\ldots
\,i_\alpha\,h_1\ldots\,h_m}_{j\,\bar j_1\ldots\,\bar j_\gamma\,j_1\ldots
\,j_\beta\,k_1\ldots\,k_n}[P]\ Y^j
\mytag{15.7}\\
&\hskip -2em
\Gamma^k_i=\sum^3_{j=0}\Gamma^k_{j\,i}\ Y^j,
\mytag{15.8}\\
&\hskip -2em
\Alpha^k_i=\sum^3_{j=0}\Alpha^k_{j\,i}\ Y^j,
\mytag{15.9}\\
&\hskip -2em
\bar{\Alpha}\vphantom{\Alpha}^k_i=\sum^3_{j=0}
\bar{\Alpha}\vphantom{\Alpha}^k_{j\,i}\ Y^j.
\mytag{15.10}
\endalign
$$
Substituting \mythetag{15.5}, \mythetag{15.6}, \mythetag{15.7},
\mythetag{15.8}, \mythetag{15.9}, and \mythetag{15.10} into the 
formula \mythetag{13.28}, then taking into account 
\mythetag{15.4}, we derive the following formula:
$$
\gathered
\nabla_{\!\!j}X^{a_1\ldots\,a_\varepsilon\bar a_1\ldots\,
\bar a_\sigma c_1\ldots\,c_e}_{b_1\ldots\,b_\eta\bar b_1\ldots\,\bar b_\zeta
d_1\ldots\,d_f}
=\sum^3_{i=0}\sum^3_{k=0}Z^{\,i}_j\,\Upsilon^k_i\,
\frac{\partial X^{a_1\ldots\,a_\varepsilon\bar a_1\ldots\,
\bar a_\sigma c_1\ldots\,c_e}_{b_1\ldots\,b_\eta\bar b_1
\ldots\,\bar b_\zeta d_1\ldots\,
d_f}}{\partial x^k}\,+\\
\vspace{1.5ex}
+\sum^{J+Q}_{P=1}\dsize\msums{2}{3}\Sb i_1,\,\ldots,\,i_\alpha\\
j_1,\,\ldots,\,j_\beta\\ \bar i_1,\,\ldots,\,\bar i_\nu\\ 
\bar j_1,\,\ldots,\,\bar j_\gamma\\
h_1,\,\ldots,\,h_m\\ k_1,\,\ldots,\,k_n\endSb
Z^{\,i_1\ldots\,i_\alpha\,\bar i_1\ldots\,\bar i_\nu\,h_1\ldots
\,h_m}_{j\,j_1\ldots\,j_\beta\,\bar j_1\ldots\,\bar j_\gamma\,k_1\ldots
\,k_n}[P]\ 
\frac{\partial X^{a_1\ldots\,a_\varepsilon\bar a_1\ldots\,
\bar a_\sigma c_1\ldots\,c_e}_{b_1\ldots\,b_\eta\bar b_1\ldots\,
\bar b_\zeta d_1\ldots\,
d_f}}{\partial S^{\,i_1\ldots\,i_\alpha\,\bar i_1\ldots\,\bar i_\nu\,
h_1\ldots\,h_m}_{j_1\ldots\,j_\beta\,\bar j_1\ldots\,\bar j_\gamma\,
k_1\ldots\,k_n}[P]}\,+\\
\vspace{1.5ex}
+\sum^{J+Q}_{P=1}\dsize\msums{2}{3}\Sb i_1,\,\ldots,\,i_\alpha\\
j_1,\,\ldots,\,j_\beta\\ \bar i_1,\,\ldots,\,\bar i_\nu\\ 
\bar j_1,\,\ldots,\,\bar j_\gamma\\
h_1,\,\ldots,\,h_m\\ k_1,\,\ldots,\,k_n\endSb
\bar Z^{\,\bar i_1\ldots\,\bar i_\nu\,i_1\ldots\,i_\alpha\,h_1\ldots
\,h_m}_{j\,\bar j_1\ldots\,\bar j_\gamma\,j_1\ldots\,j_\beta\,k_1\ldots
\,k_n}[P]\ 
\frac{\partial X^{a_1\ldots\,a_\varepsilon\bar a_1\ldots\,\bar a_\sigma 
c_1\ldots\,c_e}_{b_1\ldots\,b_\eta\bar b_1\ldots\,\bar b_\zeta d_1\ldots\,
d_f}}{\partial\overline{S^{\,i_1\ldots\,i_\alpha\,\bar i_1\ldots\,
\bar i_\nu\,h_1\ldots\,h_m}_{j_1\ldots\,j_\beta\,\bar j_1\ldots\,
\bar j_\gamma\,k_1\ldots\,k_n}[P]}\,}\,+\\
\vspace{1.5ex}
\kern -9em
+\sum^\varepsilon_{\mu=1}\sum^2_{v_\mu=1}\Alpha^{a_\mu}_{j\,v_\mu}\ 
X^{a_1\ldots\,v_\mu\,\ldots\,a_\varepsilon\bar a_1\ldots\,\bar a_\sigma
c_1\ldots\,c_e}_{b_1\ldots\,\ldots\,\ldots\,b_\eta\bar b_1\ldots\,
\bar b_\zeta d_1\ldots\,d_f}\,-\\
\kern 9em-\sum^\eta_{\mu=1}\sum^2_{w_\mu=1}\Alpha^{w_\mu}_{j\,b_\mu}\
X^{a_1\ldots\,\ldots\,\ldots\,a_\varepsilon\bar a_1\ldots\,\bar a_\sigma
c_1\ldots\,c_e}_{b_1\ldots\,w_\mu\,\ldots\,b_\eta\bar b_1\ldots\,
\bar b_\zeta d_1\ldots\,d_f}\,+\\
\kern -9em
+\sum^\sigma_{\mu=1}\sum^2_{v_\mu=1}
\bar{\Alpha}\vphantom{\Alpha}^{\bar a_\mu}_{j\,v_\mu}\ 
X^{a_1\ldots\,a_\varepsilon\bar a_1\ldots\,v_\mu\,\ldots\,\bar a_\sigma
c_1\ldots\,c_e}_{b_1\ldots\,b_\eta\bar b_1\ldots\,\ldots\,\ldots\,
\bar b_\zeta d_1\ldots\,d_f}\,-\\
\kern 9em-\sum^\zeta_{\mu=1}\sum^2_{w_\mu=1}
\bar{\Alpha}\vphantom{\Alpha}^{w_\mu}_{j\,\bar b_\mu}\
X^{a_1\ldots\,a_\varepsilon\bar a_1\ldots\,\ldots\,\ldots\,\bar a_\sigma
c_1\ldots\,c_e}_{b_1\ldots\,b_\eta\bar b_1\ldots\,w_\mu\,\ldots\,
\bar b_\zeta d_1\ldots\,d_f}\,+\\
\kern -9em+\sum^e_{\mu=1}\sum^3_{v_\mu=0}\Gamma^{c_\mu}_{j\,v_\mu}\ 
X^{a_1\ldots\,a_\varepsilon\bar a_1\ldots\,\bar a_\sigma
c_1\ldots\,v_\mu\,\ldots\,c_e}_{b_1\ldots\,b_\eta\bar b_1\ldots\,
\bar b_\zeta d_1\ldots\,\ldots\,\ldots\,d_f}\,-\\
\kern 9em-\sum^f_{\mu=1}\sum^3_{w_\mu=0}\Gamma^{w_\mu}_{j\,b_\mu}\
X^{a_1\ldots\,a_\varepsilon\bar a_1\ldots\,\bar a_\sigma
c_1\ldots\,\ldots\,\ldots\,c_e}_{b_1\ldots\,b_\eta\bar b_1\ldots\,
\bar b_\zeta d_1\ldots\,w_\mu\,\ldots\,d_f}.
\endgathered\qquad
\mytag{15.11}
$$
The formula \mythetag{15.11} is an explicit formula for the
covariant derivative $\nabla_{\!\!j}$ associated with the covariant
differentiation $\nabla$. It is also called a {\it coordinate
representation\/} of the covariant differentiation $\nabla$. The 
quantities $Z^{\,i}_j$, $Z^{\,i_1\ldots\,i_\alpha\,\bar i_1\ldots\,
\bar i_\nu\,h_1\ldots\,h_m}_{j\,j_1\ldots\,j_\beta\,\bar j_1\ldots\,
\bar j_\gamma\,k_1\ldots\,k_n}[P]$, $\bar Z^{\,\bar i_1\ldots\,
\bar i_\nu\,i_1\ldots\,i_\alpha\,h_1\ldots\,h_m}_{j\,\bar j_1\ldots
\,\bar j_\gamma\,j_1\ldots\,j_\beta\,k_1\ldots\,k_n}[P]$, $\Alpha^k_{j\,i}$,
$\bar{\Alpha}\vphantom{\Alpha}k_{j\,i}$, and $\Gamma^k_{j\,i}$ are the
components of the covariant differentiation $\nabla$ in this coordinate
representation. They obey some definite transformation rules under a change
of a local chart. From \mythetag{13.29}, \mythetag{13.30}, \mythetag{13.31},
\mythetag{13.32}, \mythetag{13.33}, and \mythetag{13.34} we derive 
a series of transformation formulas for these quantities 
$Z^{\,i}_j$, $Z^{\,i_1\ldots\,i_\alpha\,\bar i_1\ldots\,\bar i_\nu\,h_1\ldots
\,h_m}_{j\,j_1\ldots\,j_\beta\,\bar j_1\ldots\,\bar j_\gamma\,k_1\ldots
\,k_n}[P]$, $\bar Z^{\,\bar i_1\ldots\,\bar i_\nu\,i_1\ldots\,i_\alpha
\,h_1\ldots\,h_m}_{j\,\bar j_1\ldots\,\bar j_\gamma\,j_1\ldots\,j_\beta
\,k_1\ldots\,k_n}[P]$, $\Alpha^k_{j\,i}$,
$\bar{\Alpha}\vphantom{\Alpha}k_{j\,i}$, and $\Gamma^k_{j\,i}$. Here is the
transformation rule for $Z^{\,i}_j$. It is the most simple:
$$
\hskip -2em
Z^{\,i}_j=\dsize\sum^3_{h=0}\sum^3_{k=0}S^i_h\ T^k_j\ \tilde Z^{\,h}_k.
\mytag{15.12}
$$
It is followed by the transformation rule for $Z^{\,i_1\ldots\,i_\alpha\,
\bar i_1\ldots\,\bar i_\nu\,h_1\ldots\,h_m}_{j\,j_1\ldots\,j_\beta\,
\bar j_1\ldots\,\bar j_\gamma\,k_1\ldots\,k_n}[P]$:
$$
\hskip 0.1em
\gathered
Z^{\,i_1\ldots\,i_\alpha\,\bar i_1\ldots\,\bar i_\nu\,h_1
\ldots\,h_m}_{j\,j_1\ldots\,j_\beta\,\bar j_1\ldots\,\bar j_\gamma
\,k_1\ldots\,k_n}[P]=\sum^3_{d=0}
\dsize\msum{2}\Sb a_1,\,\ldots,\,a_\alpha\\ b_1,\,\ldots,\,b_\beta\\
\bar a_1,\,\ldots,\,\bar a_\nu\\ \bar b_1,\,\ldots,\,\bar b_\gamma\endSb
\dsize\msum{3}\Sb c_1,\,\ldots,\,c_m\\ d_1,\,\ldots,\,d_n\endSb
\goth S^{\,i_1}_{a_1}\ldots\,\goth S^{\,i_\alpha}_{a_\alpha}\times\\
\times\,\goth T^{b_1}_{j_1}\ldots\,\goth T^{b_\beta}_{j_\beta}
\ \overline{\goth S^{\,\bar i_1}_{\bar a_1}}
\ldots\,\overline{\goth S^{\,\bar i_\nu}_{\bar a_\nu}}\ 
\ \overline{\goth T^{\,\bar b_1}_{\bar j_1}}
\ldots\,\overline{\goth T^{\,\bar b_\gamma}_{\bar j_\gamma}}\ 
S^{h_1}_{c_1}\ldots\,S^{h_m}_{c_m}\,
T^{\,d_1}_{k_1}\ldots\,T^{\,d_n}_{k_n}\,\times\\
\times\,T^d_j\ \tilde Z^{a_1\ldots\,a_\alpha\,\bar a_1\ldots\,
\bar a_\nu\,c_1\ldots\,c_m}_{d\,b_1\ldots\,b_\beta\,\bar b_1\ldots\,
\bar b_\gamma\,d_1\ldots\,d_n}[P]-\sum^\alpha_{\mu=1}\sum^3_{i=0}
\sum^3_{c=0}\sum^3_{d=0}\sum^2_{v_\mu=1}\vartheta^{\,i_\mu}_{i\,v_\mu}
\,\times\\
\times\ S^{\,i_1\ldots\,v_\mu\ldots\,i_\alpha\,\bar i_1\ldots\,\bar i_\nu\,
h_1\ldots\,h_m}_{\,j_1\ldots\,\ldots\,\ldots\,j_\beta\,\bar j_1\ldots\,
\bar j_\gamma\,k_1\ldots\,k_n}[P]\ S^i_c\ T^d_j\ \tilde Z^c_d
+\sum^\beta_{\mu=1}\sum^3_{i=0}\sum^3_{c=0}\sum^3_{d=0}\sum^2_{w_\mu=1}
\vartheta^{\,w_\mu}_{i\,j_\mu}\,\times\\
\times\ S^{\,i_1\ldots\,\ldots\,\ldots\,i_\alpha\,\bar i_1\ldots\,
\bar i_\nu\,h_1\ldots\,h_m}_{\,j_1\ldots\,w_\mu\,\ldots\,j_\beta\,
\bar j_1\ldots\,\bar j_\gamma\,k_1\ldots\,k_n}[P]\ S^i_c\ T^d_j\ 
\tilde Z^c_d-\sum^\nu_{\mu=1}\sum^3_{i=0}\sum^3_{c=0}\sum^3_{d=0}
\sum^2_{v_\mu=1}\overline{\vartheta^{\,\bar i_\mu}_{i\,v_\mu}}\,\times\\
\times\ S^{\,i_1\ldots\,i_\alpha\,\bar i_1\ldots\,v_\mu\ldots\,\bar i_\nu\,
h_1\ldots\,h_m}_{\,j_1\ldots\,j_\beta\,\bar j_1\ldots\,\ldots\,\ldots\,
\bar j_\gamma\,k_1\ldots\,k_n}[P]\ S^i_c\ T^d_j\ \tilde Z^c_d
+\sum^\gamma_{\mu=1}\sum^3_{i=0}\sum^3_{c=0}\sum^3_{d=0}\sum^2_{w_\mu=1}
\overline{\vartheta^{\,w_\mu}_{i\,\bar j_\mu}}\,\times\\
\times\ S^{\,i_1\ldots\,i_\alpha\,\bar i_1\ldots\,\ldots\,\ldots\,
\bar i_\nu\,h_1\ldots\,h_m}_{\,j_1\ldots\,j_\beta\,\bar j_1\ldots\,
w_\mu\,\ldots\,\bar j_\gamma\,k_1\ldots\,k_n}[P]\ S^i_c\ T^d_j\ 
\tilde Z^c_d-\sum^m_{\mu=1}\sum^3_{i=0}\sum^3_{c=0}\sum^3_{d=0}
\sum^3_{v_\mu=0}\theta^{\,i_\mu}_{i\,v_\mu}\,\times\\
\times\ S^{\,i_1\ldots\,i_\alpha\,\bar i_1\ldots\,\bar i_\nu\,
h_1\ldots\,v_\mu\ldots\,h_m}_{\,j_1\ldots\,j_\beta\,\bar j_1\ldots\,
\bar j_\gamma\,k_1\ldots\,\ldots\,\ldots\,k_n}[P]\ S^i_c\ T^d_j\ 
\tilde Z^c_d+\sum^n_{\mu=1}\sum^3_{i=0}\sum^3_{c=0}\sum^3_{d=0}
\sum^3_{w_\mu=0}\theta^{\,w_\mu}_{i\,j_\mu}\,\times\\
\times\ S^{\,i_1\ldots\,i_\alpha\,\bar i_1\ldots\,\bar i_\nu\,
h_1\ldots\,\ldots\,\ldots\,h_m}_{\,j_1\ldots\,j_\beta\,\bar j_1\ldots\,
\bar j_\gamma\,k_1\ldots\,w_\mu\,\ldots\,k_n}[P]\ S^i_c\ T^d_j\ 
\tilde Z^c_d.\vphantom{\sum^n_{\mu=1}}
\endgathered
\mytag{15.13}
$$
The transformation rules for the quantities $\Alpha^k_{j\,i}$, 
$\bar{\Alpha}\vphantom{\Alpha}^k_{j\,i}$, and $\Gamma^k_{j\,i}$
are presented by the following three formulas, which are not so
huge as \mythetag{15.13}:
$$
\align
&\hskip -2em
\Gamma^k_{j\,i}=\dsize\sum^3_{b=0}\sum^3_{a=0}\sum^3_{c=0}
S^k_a\,T^b_i\,T^c_j\ \tilde\Gamma^a_{c\,b}+
\sum^3_{a=0}Z^a_j\,\theta^k_{ai},
\mytag{15.14}\\
\vspace{2ex}
&\hskip -2em
\Alpha^k_{j\,i}=\dsize\sum^2_{b=1}\sum^2_{a=1}\sum^3_{c=0}
\goth S^k_a\,\goth T^b_i\,T^c_j\ \tilde{\Alpha}\vphantom{\Alpha}^a_{c\,b}
+\sum^3_{a=0}Z^a_j\,\vartheta^k_{ai},
\mytag{15.15}\\
\vspace{2ex}
&\hskip -2em
\bar{\Alpha}\vphantom{\Alpha}^k_{j\,i}=\sum^2_{b=1}\sum^2_{a=1}
\sum^3_{c=0}\overline{\goth S^k_a}\ \overline{\goth T^b_i}\,T^c_j\  
\tilde{\bar{\Alpha}}\vphantom{\Alpha}^a_b
+\sum^3_{a=0}Z^a_j\,\overline{\vartheta^k_{ai}}.
\mytag{15.16}
\endalign
$$
And finally, we write the transformation rule for $\bar Z^{\,\bar i_1
\ldots\,\bar i_\nu\,i_1\ldots\,i_\alpha\,h_1\ldots\,h_m}_{j\,\bar j_1
\ldots\,\bar j_\gamma\,j_1\ldots\,j_\beta\,k_1\ldots\,k_n}[P]$:
$$
\hskip 0.1em
\gathered
\bar Z^{\,\bar i_1\ldots\,\bar i_\nu\,i_1\ldots\,i_\alpha\,h_1\ldots
\,h_m}_{j\,\bar j_1\ldots\,\bar j_\gamma\,j_1\ldots\,j_\beta\,k_1
\ldots\,k_n}[P]=\sum^3_{d=0}
\dsize\msum{2}\Sb a_1,\,\ldots,\,a_\alpha\\ b_1,\,\ldots,\,b_\beta\\
\bar a_1,\,\ldots,\,\bar a_\nu\\ \bar b_1,\,\ldots,\,\bar b_\gamma\endSb
\dsize\msum{3}\Sb c_1,\,\ldots,\,c_m\\ d_1,\,\ldots,\,d_n\endSb
\goth S^{\,\bar i_1}_{\bar a_1}\ldots\,
\goth S^{\,\bar i_\nu}_{\bar a_\nu}\times\\
\times\,\goth T^{\bar b_1}_{\bar j_1}\ldots\,
\goth T^{\bar b_\gamma}_{\bar j_\gamma}
\ \overline{\goth S^{\,i_1}_{a_1}}\ldots\,
\overline{\goth S^{\,i_\alpha}_{a_\alpha}}\ 
\ \overline{\goth T^{\,b_1}_{j_1}}
\ldots\,\overline{\goth T^{\,b_\beta}_{j_\beta}}\ 
S^{h_1}_{c_1}\ldots\,S^{h_m}_{c_m}\,
T^{\,d_1}_{k_1}\ldots\,T^{\,d_n}_{k_n}\,\times\\
\times\,T^d_j\ \tbZ^{\bar a_1\ldots\,\bar a_\nu\,a_1\ldots\,a_\alpha
\,c_1\ldots\,c_m}_{d\,\bar b_1\ldots\,\bar b_\gamma\,b_1\ldots\,b_\beta
\,d_1\ldots\,d_n}[P]-\sum^\alpha_{\mu=1}\sum^3_{i=0}\sum^3_{c=0}
\sum^3_{d=0}\sum^2_{v_\mu=1}\overline{\vartheta^{\,i_\mu}_{i\,v_\mu}}
\,\times\\
\times\ \overline{S^{\,i_1\ldots\,v_\mu\ldots\,i_\alpha\,\bar i_1
\ldots\,\bar i_\nu\,h_1\ldots\,h_m}_{\,j_1\ldots\,\ldots\,\ldots\,
j_\beta\,\bar j_1\ldots\,\bar j_\gamma\,k_1\ldots\,k_n}[P]}\ S^i_c
\ T^d_j\ \tilde Z^c_d+\sum^\beta_{\mu=1}\sum^3_{i=0}\sum^3_{c=0}
\sum^3_{d=0}\sum^2_{w_\mu=1}
\overline{\vartheta^{\,w_\mu}_{i\,j_\mu}}\,\times\\
\times\ \overline{S^{\,i_1\ldots\,\ldots\,\ldots\,i_\alpha\,\bar i_1
\ldots\,\bar i_\nu\,h_1\ldots\,h_m}_{\,j_1\ldots\,w_\mu\,\ldots\,j_\beta
\,\bar j_1\ldots\,\bar j_\gamma\,k_1\ldots\,k_n}[P]}\ S^i_c\ T^d_j\ 
\tilde Z^c_d-\sum^\nu_{\mu=1}\sum^3_{i=0}\sum^3_{c=0}\sum^3_{d=0}
\sum^2_{v_\mu=1}\vartheta^{\,\bar i_\mu}_{i\,v_\mu}\,\times\\
\times\ \overline{S^{\,i_1\ldots\,i_\alpha\,\bar i_1\ldots\,v_\mu
\ldots\,\bar i_\nu\,h_1\ldots\,h_m}_{\,j_1\ldots\,j_\beta\,\bar j_1
\ldots\,\ldots\,\ldots\,\bar j_\gamma\,k_1\ldots\,k_n}[P]}\ S^i_c\ 
T^d_j\ \tilde Z^c_d+\sum^\gamma_{\mu=1}\sum^3_{i=0}\sum^3_{c=0}
\sum^3_{d=0}\sum^2_{w_\mu=1}\vartheta^{\,w_\mu}_{i\,\bar j_\mu}\,
\times\\
\times\ \overline{S^{\,i_1\ldots\,i_\alpha\,\bar i_1\ldots\,\ldots\,
\ldots\,\bar i_\nu\,h_1\ldots\,h_m}_{\,j_1\ldots\,j_\beta\,\bar j_1
\ldots\,w_\mu\,\ldots\,\bar j_\gamma\,k_1\ldots\,k_n}[P]}\ S^i_c
\ T^d_j\ \tilde Z^c_d-\sum^m_{\mu=1}\sum^3_{i=0}\sum^3_{c=0}
\sum^3_{d=0}\sum^3_{v_\mu=0}\theta^{\,i_\mu}_{i\,v_\mu}\,\times\\
\times\ \overline{S^{\,i_1\ldots\,i_\alpha\,\bar i_1\ldots\,\bar i_\nu
\,h_1\ldots\,v_\mu\ldots\,h_m}_{\,j_1\ldots\,j_\beta\,\bar j_1\ldots\,
\bar j_\gamma\,k_1\ldots\,\ldots\,\ldots\,k_n}[P]}\ S^i_c\ T^d_j\ 
\tilde Z^c_d+\sum^n_{\mu=1}\sum^3_{i=0}\sum^3_{c=0}\sum^3_{d=0}
\sum^3_{w_\mu=0}\theta^{\,w_\mu}_{i\,j_\mu}\,\times\\
\times\ \overline{S^{\,i_1\ldots\,i_\alpha\,\bar i_1\ldots\,\bar i_\nu
\,h_1\ldots\,\ldots\,\ldots\,h_m}_{\,j_1\ldots\,j_\beta\,\bar j_1\ldots
\,\bar j_\gamma\,k_1\ldots\,w_\mu\,\ldots\,k_n}[P]}\ S^i_c\ T^d_j\ 
\tilde Z^c_d.\vphantom{\sum^n_{\mu=1}}
\endgathered
\mytag{15.17}
$$
The formulas \mythetag{15.3} and \mythetag{15.11} yield the explicit
expressions for arbitrary covariant derivatives in general case. 
However, below we consider some specializations of these formulas 
which appear to be more valuable than the initial formulas \mythetag{15.3}
and \mythetag{15.11} themselves.
\head
16. Degenerate covariant differentiations.
\endhead
\mydefinition{16.1} A covariant differentiation $\nabla$ is said to 
be {\it degenerate\/} if \ $\nabla_{\bold Y}\psi=0$ for any extended
scalar field $\psi$ and for any extended vector field $\bold Y$.
\enddefinition
     This definition is concordant with the
definition~\mythedefinition{14.1}. For degenerate covariant differentiations 
we have a theorem which is analogous to the theorems~\mythetheorem{14.1}.
\mytheorem{16.1} Defining a degenerate covariant differentiation
$\nabla$ of the algebra of extended spin-tensorial fields $\bold S(M)$ 
is equivalent to defining three extended spin-tensorial fields $\eufb S$,
$\bar{\eufb S}$, and $\bold S$ of the types $(1,1|0,0|0,1)$,
$(0,0|1,1|0,1)$, and $(0,0|0,0|1,2)$ respectively.
\endproclaim
\demo{Proof} Let $\nabla$ be an arbitrary degenerate covariant
differentiation.  For this differentiation from \mythetag{15.5}, 
\mythetag{15.6}, and \mythetag{15.7} we derive 
the following equalities:
$$
\gather
\hskip -2em
Z^{\,i}_j=0,
\mytag{16.1}\\
\vspace{1ex}
\hskip -2em
Z^{\,i_1\ldots\,i_\alpha\,\bar i_1\ldots\,\bar i_\nu\,h_1\ldots
\,h_m}_{j\,j_1\ldots\,j_\beta\,\bar j_1\ldots\,\bar j_\gamma\,k_1
\ldots\,k_n}[P]=0,
\mytag{16.2}\\
\vspace{1.5ex}
\hskip -2em
\bar Z^{\,\bar i_1\ldots\,\bar i_\nu\,i_1\ldots
\,i_\alpha\,h_1\ldots\,h_m}_{j\,\bar j_1\ldots\,\bar j_\gamma\,j_1\ldots
\,j_\beta\,k_1\ldots\,k_n}[P]=0.
\mytag{16.3}
\endgather
$$
The equalities \mythetag{16.1}, \mythetag{16.2}, and \mythetag{16.3}
are analogs of the equalities \mythetag{14.1}, \mythetag{14.2}, and 
\mythetag{14.3} respectively. Applying \mythetag{16.1} to \mythetag{15.14},
\mythetag{15.15}, and \mythetag{15.16}, we obtain the following 
transformation rules for the $\Alpha$-components and the
$\Gamma$-components of the degenerate covariant differentiation
$\nabla$:
$$
\align
&\hskip -2em
\Alpha^k_{j\,i}=\dsize\sum^2_{b=1}\sum^2_{a=1}\sum^3_{c=0}
\goth S^k_a\,\goth T^b_i\,T^c_j\ \tilde{\Alpha}\vphantom{\Alpha}^a_{c\,b},
\mytag{16.4}\\
\vspace{2ex}
&\hskip -2em
\bar{\Alpha}\vphantom{\Alpha}^k_{j\,i}=\sum^2_{b=1}\sum^2_{a=1}
\sum^3_{c=0}\overline{\goth S^k_a}\ \overline{\goth T^b_i}\,T^c_j\  
\tilde{\bar{\Alpha}}\vphantom{\Alpha}^a_b,
\mytag{16.5}\\
\vspace{2ex}
&\hskip -2em
\Gamma^k_{j\,i}=\dsize\sum^3_{b=0}\sum^3_{a=0}\sum^3_{c=0}
S^k_a\,T^b_i\,T^c_j\ \tilde\Gamma^a_{c\,b}.
\mytag{16.6}
\endalign
$$
From \mythetag{16.4}, \mythetag{16.5}, and \mythetag{16.6}, we see that
$\Alpha$-components of the degenerate covariant differentiation $\nabla$
define two extended spin-tensorial fields of the types $(1,1|0,0|0,1)$ 
and $(0,0|1,1|0,1)$, while its  $\Gamma$-components define an extended
spin-tensorial field of the type $(0,0|0,0|1,2)$. The theorem is proved.
\qed\enddemo
\head
17. Horizontal and vertical covariant differentiations.
\endhead
     Suppose again that $N$ is a composite spin-tensorial bundle over the
space-time manifold $M$ (in the sense of \mythetag{9.7}). 
Let $\nabla$ be a covariant differentiation of the algebra of extended
spin-tensorial fields $\bold S(M)$. Then $D=\nabla_{\bold Y}$ is a
differentiation of $\bold S(M)$, its restriction to the set of scalar
fields is given by some vector field $\bold Z=\bold Z(\bold Y)$ in $N$ 
(see formula \mythetag{13.2} above). In other words, we have a homomorphism
$$
\hskip -2em
\Bbb CT^1_0(M)\to\Bbb CT^1_0(N)
\mytag{17.1}
$$
that maps an extended vector field $\bold Y$ of $M$ to some regular
vector field of $N$. Applying the localization theorem~\mythetheorem{15.1}
to the homomorphism \mythetag{17.1}, we come to the following definition
and to the theorem after it.
\mydefinition{17.1} Suppose that for each point $q$ of the composite 
spin-tensorial bundle $N$ over the space-time $M$ some $\Bbb C$-linear 
map of the vector spaces 
$$
\hskip -2em
f_q\!:\ \Bbb CT_{\pi(q)}(M)\to \Bbb CT_q(N)
\mytag{17.2}
$$ 
is given. Then we say that a {\it lift\/} of vectors 
from $M$ to the bundle $N$ is defined.
\enddefinition
\mytheorem{17.1} Any homomorphism of $\goth F_{\Bbb C}(N)$-modules
\mythetag{17.1} is uniquely associated with some smooth lift of vectors
from $M$ to $N$. It is represented by this lift as a collection of 
$\Bbb C$-linear maps \mythetag{17.2} specific to each point $q\in N$.
\endproclaim
    The $\Bbb C$-linear maps \mythetag{17.2} here are special cases 
of general $\Bbb C$-linear maps \mythetag{15.1} declared in the 
localization theorem~\mythetheorem{15.1}.\par
    Now let's consider the canonical projection $\pi\!:\,N\to M$. The
differential of this map (upon complexification) acts in the 
direction opposite to the lift of vectors \mythetag{17.2}  introduced in the
definition~\mythedefinition{17.1}. Indeed, we have $\pi_*\!:\,
\Bbb CT_q(N)\to\Bbb CT_{\pi(q)}(M)$ at each point $q\in N$. Therefore, 
the composition $f\compos\pi_*$ acts from $\Bbb CT_{\pi(q)}(M)$ to 
$\Bbb CT_{\pi(q)}(M)$. This composite map determines an extended operator
field (a spin-tensorial field of the type $(0,0|0,0|1,1)$).\par
\mydefinition{17.2} A lift of vectors $f$ from $M$ to $N$ is called
{\it vertical\/} if $\pi_*\compos f=0.$
\enddefinition
\mydefinition{17.3} A lift of vectors $f$ from $M$ to $N$ is called
{\it horizontal\/} if $\pi_*\compos f=\idop$, i\.\,e\. if the composition
$\pi_*\compos f$ coincides with the field of identical operators.
\enddefinition
     Like any other bundle, the composite spin-tensorial bundle $N$
naturally subdivides into fibers over the points of the base manifold $M$.
The set of vectors tangent to the fiber at a point $q$ (upon
complexification) is a $\Bbb C$-linear subspace within the complexified
tangent space $\Bbb CT_q(N)$. This subspace coincides with the kernel of
the complexified mapping $\pi_*$. We denote this subspace 
$$
\hskip -2em
V_q(N)=\Ker\pi_*
\mytag{17.3}
$$ 
and call it the {\it vertical subspace}. Any vertical lift of vectors 
determines a set linear mappings from $\Bbb CT_{\pi(q)}$  to the vertical 
subspace \mythetag{17.3} at each point $q\in N$.
\mylemma{17.1} The difference of two horizontal lifts is a vertical lift 
of vectors from the base manifold $M$ to the bundle $N$.
\endproclaim
Indeed, if one takes two horizontal lifts of vectors $f_1$ and $f_2$, 
then $\pi_*\compos(f_1-f_2)=\ =\pi_*\compos f_1-\pi_*\compos f_2=\idop
-\idop=0$. This means that the difference $f_1-f_2$ is a vertical lift
according to the definition~\mythedefinition{17.2}.\par
    Each covariant differentiation $\nabla$ is associated with some
lift of vectors (see the definition~\mythedefinition{17.1}, and the 
theorem~\mythetheorem{17.1} above).
\mydefinition{17.4} A covariant differentiation $\nabla$ is called 
a {\it horizontal covariant differentiation\/} (or a {\it vertical 
covariant differentiation\/}) if the corresponding lift of vectors 
is {\it horizontal\/} (or {\it vertical}).
\enddefinition
\mylemma{17.2} The difference of two horizontal covariant differentiations
is a vertical covariant differentiation.
\endproclaim
    The lemma~\mythelemma{17.2} is an immediate consequence of the
lemma~\mythelemma{17.1}.
\head
18. Native extended spin-tensorial fields\\
and vertical multivariate differentiations.
\endhead
     Let $N$ be a composite spin-tensorial bundle over the space-time
manifold $M$. Then each its point $q$ is represented by a list 
$q=(p,\,\bold S[1],\,\ldots,\,\bold S[J+Q])$, where $p\in M$ and 
$\bold S[1],\,\ldots,\,\bold S[J+Q]$ are some spin-tensors at the point 
$p$ (see formula \mythetag{9.8} above). Let's consider the map that takes 
$q$ to the $P$-th spin-tensor $\bold S[P]$ in this list. According to the 
definition~\mythedefinition{10.2}, this map is an extended spin-tensorial
field of the type $(\alpha_P,\beta_P|\nu_P,\gamma_P|m_P,n_P)$. It is 
canonically associated with the bundle $N$. Therefore, it is called a 
{\it native extended spin-tensorial field}. Totally, we have $J+Q$ native
extended spin-tensorial fields associated with the composite spin-tensorial
bundle $N$, we denote them $\bold S[1],\,\ldots,\,\bold S[J+Q]$.\par
\mydefinition{18.1} A {\it multivariate differentiation\/} of the type
$(\beta,\alpha|\gamma,\nu|n,m)$ in the algebra $\bold S(M)$ is a 
homomorphism of $\goth F_{\Bbb C}(N)$-modules
$$
\hskip -2em
\nabla\!:\,S^\alpha_\beta\bar S^\nu_\gamma T^m_n(M)\to\goth D_{\bold S}(M).
\mytag{18.1}
$$
If\/ $\bold Y$ is an extended spin-tensorial field of the type
$(\alpha,\beta|\nu,\gamma|m,n)$, then we can apply the homomorphism
\mythetag{18.1} to it. As a result we get the differentiation 
$D=\nabla_{\bold Y}$ of the algebra of extended spin-tensorial fields
$\bold S(M)$. It is called the {\it multivariate differentiation along
the spin-tensorial field\/} $\bold Y$.
\enddefinition
     Note that the type of a multivariate differentiation 
$(\beta,\alpha|\gamma,\nu|n,m)$ in the above
definition~\mythedefinition{18.1} is dual to the type of the
module $S^\alpha_\beta\bar S^\nu_\gamma T^m_n(M)$ in the formula 
\mythetag{18.1}. If $\alpha=0$, $\beta=0$, $\nu=0$, $\gamma=0$,
$m=1$, and $n=0$, then the definition~\mythedefinition{18.1} reduces 
to the definition~\mythedefinition{15.1}. This means that a covariant
differentiation is a special multivariate differentiation whose type
is $(0,0|0,0|0,1)$. Similarly, a multivariate differentiation of the 
type $(0,0|0,0|1,0)$ is called a {\it contravariant differentiation}.
\par
    Let $\nabla$ be some multivariate differentiation of the algebra 
of extended spin-tensorial fields $\bold S(M)$. Then, applying the
localization theorem~\mythetheorem{15.1} to the homomorphism
$\nabla\!:\,S^\alpha_\beta\bar S^\nu_\gamma T^m_n(M)\to
\goth D_{\bold S}(M)$, we find that this homomorphism is composed by 
$\Bbb C$-linear maps $S^{\,\alpha}_\beta\bar S^{\nu}_\gamma 
T^m_n(\pi(q),M)\to\goth D(q,M)$ specific to each point $q\in N$. This 
fact is expressed by the formula coinciding with \mythetag{15.2}:
$$
\hskip -2em
\nabla_{\bold Y}\bold X=C(\bold Y\otimes\nabla\bold X).
\mytag{18.2}
$$
However, instead of \mythetag{15.3} in this case we have
$$
\gathered
\nabla_{\bold Y}\bold X
=\dsize\msums{2}{3}\Sb a_1,\,\ldots,\,a_\varepsilon\\
b_1,\,\ldots,\,b_\eta\\ \bar a_1,\,\ldots,\,\bar a_\sigma\\ 
\bar b_1,\,\ldots,\,\bar b_\zeta\\
c_1,\,\ldots,\,c_e\\ d_1,\,\ldots,\,d_f\endSb
\dsize\msums{2}{3}\Sb i_1,\,\ldots,\,i_\alpha\\
j_1,\,\ldots,\,j_\beta\\ \bar i_1,\,\ldots,\,\bar i_\nu\\ 
\bar j_1,\,\ldots,\,\bar j_\gamma\\
h_1,\,\ldots,\,h_m\\ k_1,\,\ldots,\,k_n\endSb
Y^{i_1\ldots\,i_\alpha\bar i_1\ldots\,
\bar i_\nu h_1\ldots\,h_m}_{j_1\ldots\,j_\beta\bar j_1\ldots\,
\bar j_\gamma k_1\ldots\,k_n}\times\\
\times\, 
\nabla^{j_1\ldots\,j_\beta\bar j_1\ldots\,\bar j_\gamma k_1\ldots
\,k_n}_{\!i_1\ldots\,i_\alpha\bar i_1\ldots\,\bar i_\nu h_1\ldots\,h_m}
X^{a_1\ldots\,a_\varepsilon\bar a_1\ldots\,\bar a_\sigma c_1\ldots
\,c_e}_{b_1\ldots\,b_\eta\bar b_1\ldots\,\bar b_\zeta d_1\ldots\,d_f}
\ \boldsymbol\Psi^{b_1\ldots\,b_\eta\bar b_1\ldots\,\bar b_\zeta d_1
\ldots\,d_f}_{a_1\ldots\,a_\varepsilon\bar a_1\ldots\,\bar a_\sigma
c_1\ldots\,c_e}.
\endgathered
\mytag{18.3}
$$
Looking at \mythetag{18.3}, we see that each multivariate differentiation
$\nabla$ of the type $(\beta,\alpha|\gamma,\nu|n,m)$ can be treated as an
operator producing the extended spin-tensorial field $\nabla\bold X$ of the
type $(\varepsilon+\beta,\eta+\alpha|\sigma+\gamma,\zeta+\nu|r+n,s+m)$ from
any given extended spin-tensorial field $\bold X$ of the type
$(\varepsilon,\eta|\sigma,\zeta|r,s)$. This operator is called the operator
of {\it multivariate differential\/} of the type
$(\beta,\alpha|\gamma,\nu|n,m)$.\par
     Let $P$ be an integer number such that $1\leqslant P\leqslant J+Q$ 
and let $\bold Y$ be an extended spin-tensorial field of the type 
$(\alpha_P,\beta_P|\nu_P,\gamma_P|m_P,n_P)$. Remember that each point $q$
of the composite spin-tensorial bundle $N$ is a list of the form
\mythetag{9.8}:
$$
\hskip -2em
q=(p,\,\bold S[1],\,\ldots,\,\bold S[J+Q]),
\mytag{18.4}
$$
Note that the $P$-th spin-tensor $\bold S[P]$ in the list \mythetag{18.4} 
has the same type as the spin-tensor $\bold Y=\bold Y_{\!q}$ (the value of 
the extended spin-tensorial field $\bold Y$ at the point $q$). They both 
belong to the same spin-tensorial space $S^{\alpha_P}_{\beta_P}
\bar S^{\nu_P}_{\gamma_P} T^{m_P}_{n_P}(p,M)$, therefore we can add them.
This means that we can treat the list
$$
\hskip -2em
q(t)=(p,\,\bold S[1],\,\ldots,\,\bold S[P]+t\,\bold Y_{\!q},
\,\ldots,\,\bold S[J+Q])
\mytag{18.5}
$$
as a one-parametric set of points in $N$, the scalar variable $t$ being 
its parameter. Thus in \mythetag{18.5} we have a line (a 
straight line) passing through the initial point $q\in N$ and lying 
completely within the fiber over the point $p=\pi(q)\in M$. Suppose that
$\bold X$ is some extended spin-tensorial field of the type 
$(\varepsilon,\eta|\sigma,\zeta|e,f)$. Denote by $\bold X(t)$ the values 
of this field at the points of the above parametric line \mythetag{18.5}:
$$
\hskip -2em
\bold X(t)=\bold X_{q(t)}.
\mytag{18.6}
$$
Since $\pi(q(t))=p=\const$ for any $t$, the values of the spin-tensor-valued
function \mythetag{18.6} all belong to the same space
$S^\varepsilon_\eta\bar S^\sigma_\zeta T^e_f(p,M)$. Hence, we can add and
subtract them, and, since $\bold X$ is smooth, we can take the following 
limit of the ratio:
$$
\hskip -2em
\dot X(t)=\lim_{\tau\to\,0}\frac{\bold X(t+\tau)-\bold X(t)}{\tau}.
\mytag{18.7}
$$
Let's denote by $\bold Z_q$ the value of the derivative
\mythetag{18.7} for $t=0$:
$$
\hskip -2em
\bold Z_q=\dot X(0)=\frac{d\bold X_{q(t)}}{dt}\,\hbox{\vrule
height14pt depth8pt width 0.5pt}_{\ t=0}.\hskip -2em
\mytag{18.8}
$$
It is easy to understand that, when $q$ is fixed, $\bold Z_q$
is a spin-tensor from the space $S^\varepsilon_\eta\bar S^\sigma_\zeta
T^e_f(p,M)$ at the point $p=\pi(q)$. By varying $q\in N$, we find that 
the spin-tensors $\bold Z_q$ constitute a smooth extended spin-tensorial
field $\bold Z$. As a result of the above considerations we have
constructed a map
$$
\hskip -2em
D\!:\,\bold S(M)\to\bold S(M).
\mytag{18.9}
$$
It is easy to check up that the map \mythetag{18.9} defined by means 
of the formulas \mythetag{18.5}, \mythetag{18.6}, \mythetag{18.7}, and 
\mythetag{18.8} is a differentiation of the algebra of extended 
spin-tensorial fields $\bold S(M)$, i\.\,e\. $D\in\goth D_{\bold S}(M)$
(see the definition~\mythedefinition{12.1} above). Moreover, due to the 
formula \mythetag{18.5} this differentiation $D$ depends on the extended 
spin-tensorial field $\bold Y$. The easiest way to study this dependence
$D=D(\bold Y)$ is to write the equality \mythetag{18.8} in a local chart, 
i\.\,e\. in some local coordinates \mythetag{9.11} and \mythetag{9.12}:
$$
\gathered
\bold Z=\dsize\msums{2}{3}\Sb a_1,\,\ldots,\,a_\varepsilon\\
b_1,\,\ldots,\,b_\eta\\ \bar a_1,\,\ldots,\,\bar a_\sigma\\ 
\bar b_1,\,\ldots,\,\bar b_\zeta\\
c_1,\,\ldots,\,c_e\\ d_1,\,\ldots,\,d_f\endSb
\dsize\msums{2}{3}\Sb i_1,\,\ldots,\,i_\alpha\\
j_1,\,\ldots,\,j_\beta\\ \bar i_1,\,\ldots,\,\bar i_\nu\\ 
\bar j_1,\,\ldots,\,\bar j_\gamma\\
h_1,\,\ldots,\,h_m\\ k_1,\,\ldots,\,k_n\endSb
\left(Y^{i_1\ldots\,i_\alpha\bar i_1\ldots\,
\bar i_\nu h_1\ldots\,h_m}_{j_1\ldots\,j_\beta\bar j_1\ldots\,
\bar j_\gamma k_1\ldots\,k_n}
\ \frac{\partial X^{a_1\ldots\,a_\varepsilon\bar a_1\ldots\,
\bar a_\sigma c_1\ldots\,c_e}_{b_1\ldots\,b_\eta\bar b_1\ldots
\,\bar b_\zeta d_1\ldots\,d_f}}{\partial S^{\,i_1\ldots\,i_\alpha
\bar i_1\ldots\,\bar i_\nu h_1\ldots\,h_m}_{j_1\ldots\,j_\beta
\bar j_1\ldots\,\bar j_\gamma k_1\ldots\,k_n}[P]
\vphantom{\vrule height 11pt depth 0pt}}\,+\right.\\
+\left.\overline{Y^{i_1\ldots\,i_\alpha\bar i_1\ldots\,
\bar i_\nu h_1\ldots\,h_m}_{j_1\ldots\,j_\beta\bar j_1\ldots\,
\bar j_\gamma k_1\ldots\,k_n}}
\ \frac{\partial X^{a_1\ldots\,a_\varepsilon\bar a_1\ldots\,
\bar a_\sigma c_1\ldots\,c_e}_{b_1\ldots\,b_\eta\bar b_1\ldots
\,\bar b_\zeta d_1\ldots\,d_f}}{\partial\overline{S^{\,i_1\ldots\,i_\alpha
\bar i_1\ldots\,\bar i_\nu h_1\ldots\,h_m}_{j_1\ldots\,j_\beta
\bar j_1\ldots\,\bar j_\gamma k_1\ldots\,k_n}[P]
\vphantom{\vrule height 11pt depth 0pt}}}\right)
\boldsymbol\Psi^{b_1\ldots\,b_\eta\bar b_1\ldots\,\bar b_\zeta d_1
\ldots\,d_f}_{a_1\ldots\,a_\varepsilon\bar a_1\ldots\,\bar a_\sigma
c_1\ldots\,c_e}.
\endgathered
\mytag{18.10}
$$
Here $\alpha=\alpha_P$, $\beta=\beta_P$, $\nu=\nu_P$, $\gamma=\gamma_P$, 
$m=m_P$, and $n=n_P$. Two partial derivatives in \mythetag{18.10} behave
like the components of two extended spin-tensorial fields under a change
of a local chart. We denote these fields by $\vnabla\bold X$ and
$\Bar\vnabla\bold X$ respectively. Here $\vnabla$ is a special sign, 
the {\tencyr\char '074}double bar nabla{\tencyr\char '076}. It was used 
in \mycite{5} for the first time. Now we can write \mythetag{18.10} as
follows:
$$
\hskip -2em
\bold Z=C(\bold Y\otimes\vnabla\bold X)+C(\tau(\bold Y)\otimes
\Bar\vnabla\bold X).
\mytag{18.11}
$$
Comparing \mythetag{18.11} with \mythetag{18.2}, we see that the formula 
\mythetag{18.11} introduces two multivariate differentiations of the types
$(\beta,\alpha|\gamma,\nu|n,m)$ and $(\gamma,\nu|\beta,\alpha|n,m)$:
$$
\hskip -2em
\aligned
&\vnabla[P]\!:\,S^\alpha_\beta\bar S^\nu_\gamma T^m_n(M)
\to\goth D_{\bold S}(M),\\
\vspace{1ex}
&\Bar\vnabla[P]\!:\,S^\nu_\gamma\bar S^\alpha_\beta T^m_n(M)
\to\goth D_{\bold S}(M).
\endaligned
\mytag{18.12}
$$
In terms of these multivariate differentiations the formula
\mythetag{18.11} is rewritten as
$$
\hskip -2em
\bold Z=\vnabla_{\bold Y}\bold X+\Bar\vnabla_{\!\tau(\bold Y)}\bold X.
\mytag{18.13}
$$
In a local chart the multivariate differentiations \mythetag{18.12} are
represented by the partial derivatives taken from the formula
\mythetag{18.10}:
$$
\align
&\hskip -2em
\vnabla^{j_1\ldots\,j_\beta\bar j_1\ldots\,\bar j_\gamma
k_1\ldots\,k_n}_{i_1\ldots\,i_\alpha\bar i_1\ldots\,\bar i_\nu
h_1\ldots\,h_m}[P]=\frac{\partial}
{\partial S^{\,i_1\ldots\,i_\alpha\bar i_1\ldots\,\bar i_\nu h_1\ldots
\,h_m}_{j_1\ldots\,j_\beta\bar j_1\ldots\,\bar j_\gamma
k_1\ldots\,k_n}[P]
\vphantom{\vrule height 11pt depth 0pt}},
\mytag{18.14}\\
\vspace{1ex}
&\hskip -2em
\Bar\vnabla^{j_1\ldots\,j_\gamma\,\bar j_1\ldots\,\bar j_\beta\,
k_1\ldots\,k_n}_{i_1\ldots\,i_\nu\,\bar i_1\ldots\,\bar i_\alpha\,
h_1\ldots\,h_m}[P]=\frac{\partial}
{\partial\overline{S^{\,\bar i_1\ldots\,\bar i_\alpha\,i_1\ldots\,
i_\nu\,h_1\ldots\,h_m}_{\bar j_1\ldots\,\bar j_\beta\,j_1\ldots\,
j_\gamma\,k_1\ldots\,k_n}[P]
\vphantom{\vrule height 11pt depth 0pt}}},
\mytag{18.15}\\
\endalign
$$
Here $\alpha=\alpha_P$, $\beta=\beta_P$, $\nu=\nu_P$, $\gamma=\gamma_P$, 
$m=m_P$, and $n=n_P$. Following the tradition, we shall use the term 
{\it multivariate derivatives} for the differential operators representing 
the differentiation $\vnabla[P]$ and $\Bar\vnabla[P]$ in \mythetag{18.14}
and \mythetag{18.15}.
\mydefinition{18.2} The multivariate differentiations $\vnabla[P]$ and
$\Bar\vnabla[P]$ defined through the formulas \mythetag{18.5}, \mythetag{18.6}, \mythetag{18.7}, \mythetag{18.8}, \mythetag{18.13} and
represented by the formulas \mythetag{18.14} and \mythetag{18.15} in local
coordinates are called the {\it $P$-th canonical vertical multivariate 
differentiation} and the {\it barred $P$-th canonical\footnotemark\ 
vertical multivariate differentiation} respectively.
\enddefinition
\footnotetext{ \ Note that $\vnabla[P]$ and $\Bar\vnabla[P]$ are 
canonically associated with the bundle $N$, their definition does 
not require any auxiliary structures like metrics and
connections.}\adjustfootnotemark{-1}
    Let $\bold S[R]$ be $R$-th native extended spin-tensorial field
associated with the composite spin-tensorial bundle $N$. Then by means 
of direct calculations in local coordinates one can derive the following
formulas:
$$
\align
\hskip -2em
\vnabla_{\!\bold Y}[P]\bold S[R]
&=\cases \bold Y &\text{for \ }P=R,\\
0 &\text{for \ }P\neq R,\endcases
\mytag{18.16}\\
\hskip -2em
\Bar\vnabla_{\!\bold Y}[P]\tau(\bold S[R])
&=\cases \bold Y &\text{for \ }P=R,\\
0 &\text{for \ }P\neq R,\endcases
\mytag{18.17}
\endalign
$$
Apart from \mythetag{18.16} and \mythetag{18.17} we can write the
following formulas:
$$
\align
\hskip -2em
\vnabla_{\!\bold Y}[P]\tau(\bold S[R])&=0\text{\ \ for all $P$ and $R$},
\mytag{18.18}\\
\vspace{1ex}
\hskip -2em
\Bar\vnabla_{\!\bold Y}[P]\bold S[R]&=0\text{\ \ for all $P$ and $R$}.
\mytag{18.19}
\endalign
$$
These formulas \mythetag{18.18} and \mythetag{18.19} are also easily 
derived by means of direct calculations in local coordinates.\par
     Like covariant differentiations (see theorem~\mythetheorem{17.1} 
and definition~\mythedefinition{17.1}), multivariate differentiations are
associated with some lifts. However, unlike covariant differentiations,
they lift not vectors, but spin-tensors, though converting them into 
tangent vectors of the bundle $N$. In the case of the canonical multivariate 
differentiation $\vnabla[P]$ for each point $q\in N$ we 
have some $\Bbb C$-linear map
$$
\hskip -2em
f[P]\!:\,S^{\alpha_P}_{\beta_P}\bar S^{\nu_P}_{\gamma_P}
T^{m_P}_{n_P}(\pi(q),M)\to\Bbb CT_q(N),
\mytag{18.20}
$$ 
The map \mythetag{18.20} takes a spin-tensor $\bold Y\in
S^{\alpha_P}_{\beta_P}\bar S^{\nu_P}_{\gamma_P}T^{m_P}_{n_P}(\pi(q),M)$ to
the following vector in the tangent space $\Bbb CT_q(N)$ of the manifold
$N$ at the point $q$:
$$
f[P](\bold Y)=\dsize\msums{2}{3}\Sb i_1,\,\ldots,\,i_\alpha\\
j_1,\,\ldots,\,j_\beta\\ \bar i_1,\,\ldots,\,\bar i_\nu\\ 
\bar j_1,\,\ldots,\,\bar j_\gamma\\
h_1,\,\ldots,\,h_m\\ k_1,\,\ldots,\,k_n\endSb
Y^{i_1\ldots\,i_\alpha\,\bar i_1\ldots\,\bar i_\nu\,h_1\ldots
\,h_m}_{j_1\ldots\,j_\beta\bar j_1\ldots\,\bar j_\gamma\,k_1\ldots\,
k_n}\ \bold W^{j_1\ldots\,j_\beta\,\bar j_1\ldots\,\bar j_\gamma
\,k_1\ldots\,k_n}_{i_1\ldots\,i_\alpha\,\bar i_1\ldots\,\bar i_\nu
\,h_1\ldots\,h_m}[P].\quad
\mytag{18.21}
$$
Similarly, in the case of the barred canonical multivariate differentiation 
$\Bar\vnabla[P]$ for each point $q\in N$ we have some $\Bbb C$-linear map
$$
\hskip -2em
\bar f[P]\!:\,S^{\nu_P}_{\gamma_P}\bar S^{\alpha_P}_{\beta_P}
T^{m_P}_{n_P}(\pi(q),M)\to\Bbb CT_q(N)
\mytag{18.22}
$$ 
This map \mythetag{18.22} is determined by the formula
$$
\aligned
\bar f[P](\bold Y)=\dsize\msums{2}{3}\Sb i_1,\,\ldots,\,i_\alpha\\
j_1,\,\ldots,\,j_\beta\\ \bar i_1,\,\ldots,\,\bar i_\nu\\ 
\bar j_1,\,\ldots,\,\bar j_\gamma\\
h_1,\,\ldots,\,h_m\\ k_1,\,\ldots,\,k_n\endSb
Y^{i_1\ldots\,i_\nu\,\bar i_1\ldots\,\bar i_\alpha\, h_1\ldots
\,h_m}_{j_1\ldots\,j_\gamma\,\bar j_1\ldots\,
\bar j_\beta\,k_1\ldots\,k_n}\ \bar\bold W^{j_1\ldots\,
j_\gamma\,\bar j_1\ldots\,\bar j_\beta\,k_1\ldots\,k_n}_{i_1\ldots
\,i_\nu\,\,\bar i_1\ldots\,\bar i_\alpha\,h_1\ldots\,h_m}[P].
\endaligned\quad
\mytag{18.23}
$$
Here again $\alpha=\alpha_P$, $\beta=\beta_P$, $\nu=\nu_P$,
$\gamma=\gamma_P$, $m=m_P$, $n=n_P$, while the vectors 
$\bold W^{j_1\ldots\,j_\beta\,\bar j_1\ldots\,\bar j_\gamma\,
k_1\ldots\,k_n}_{i_1\ldots\,i_\alpha\,\bar i_1\ldots\,\bar i_\nu
\,h_1\ldots\,h_m}[P]$ and $\bar\bold W^{j_1\ldots\,
j_\gamma\,\bar j_1\ldots\,\bar j_\beta\,k_1\ldots\,k_n}_{i_1\ldots
\,i_\nu\,\,\bar i_1\ldots\,\bar i_\alpha\,h_1\ldots\,h_m}[P]$
in \mythetag{18.21} and \mythetag{18.23} are given by the formulas
\mythetag{9.19} and \mythetag{9.20}. In a coordinate-free form the 
formulas \mythetag{18.21} and \mythetag{18.23} can be interpreted
as follows: the vector $f[P](\bold Y)+\bar f[P](\tau(\bold Y))$ 
is the tangent vector of the parametric curve \mythetag{18.5} at its 
initial point $q=q(0)$.\par
     Let's consider the image of the $\Bbb C$-linear map \mythetag{18.20}.
We denote it $V_q[P](N)$. Then from \mythetag{18.21} one easily derives 
that $V_q[P](N)$ is a subspace within the vertical subspace $V_q(N)$ of 
the tangent space $\Bbb CT_q(N)$. Similarly, we denote by $\bar V_q[P](N)$
the image of the $\Bbb C$-linear map \mythetag{18.22}. This image also is
a subspace within the vertical subspace $V_q(N)$. Moreover, we have
$$
\hskip -2em
\gathered
V_q(N)=V_q[1](N)\oplus\ldots\oplus V_q[J+Q](N)\,+\\
\vspace{1ex}
+\,\bar V_q[1](N)\oplus\ldots\oplus \bar V_q[J+Q](N).
\endgathered
\mytag{18.24}
$$
The formula \mythetag{18.24} is a well-known fact, it follows from
\mythetag{9.7}. Due to the inclusions
$$
\hskip -2em
\aligned
\Img f[P]&=V_q[P](N)\subset V_q(N),\\
\Img\bar f[P]&=\bar V_q[P](N)\subset V_q(N)
\endaligned
\mytag{18.25}
$$
the multivariate differentiations \mythetag{18.12} both are vertical
differentiations.\par
\head
19. Horizontal covariant differentiations\\
and extended connections.
\endhead
    Let $\nabla$ be some horizontal covariant differentiation of 
the algebra of extended spin-tensorial fields $\bold S(M)$ and let $f$ be 
the horizontal lift of vectors associated with it (see the
definition~\mythedefinition{17.3}). The horizontality of $f$ means 
that the image of the linear map \mythetag{17.2} is some $4$-dimensional
subspace $H_q(N)$ within the tangent space $\Bbb CT_q(N)$. It is called a
{\it horizontal subspace}. Due to $\pi_*\compos f=\idop$ the mappings
$$
\xalignat 2
&\hskip -2em
f\!:\Bbb CT_{\pi(q)}(M)\to H_q(N),
&&\pi_*\!:H_q(N)\to\Bbb CT_{\pi(q)}(M)
\mytag{19.1}
\endxalignat
$$
are inverse to each other. Due to the same equality $\pi_*\compos f=\idop$
the sum of the vertical and horizontal subspaces is a direct sum:
$$
\hskip -2em
H_q(N)\oplus V_q(N)=\Bbb CT_q(N).
\mytag{19.2}
$$
\mytheorem{19.1} Defining a horizontal lift of vectors from $M$ to $N$
is equivalent to fixing some direct complement $H_q(N)$ of the vertical
subspace $V_q(N)$ within the tangent space $\Bbb CT_q(N)$ at each point
$q\in N$.
\endproclaim
\demo{Proof} Suppose that some horizontal lift of vectors $f$ is given.
Then the subspace $H_q(N)$ at the point $q$ is determined as the image
of the mapping \mythetag{17.2}, while the relationship \mythetag{19.2}
is derived from $\pi_*\compos f=\idop$ and from \mythetag{17.3}.\par
     Conversely, assume that at each point $q\in N$ we have a subspace
$H_q(N)$ complementary to $V_q(N)$. Then at each point $q\in N$ the
relationship \mythetag{19.2} is fulfilled. The kernel of the mapping
$\pi_*\!:\Bbb CT_q(N)\to\Bbb CT_{\pi(q)}(M)$ coincides with $V_q(N)$,
therefore, the restriction of $\pi_*$ to the horizontal subspace $H_q(N)$
is a bijection. The lift of vectors $f$ from $M$ to $N$ then can be defined
as the inverse mapping for $\pi_*\!:H_q(N)\to\Bbb CT_{\pi(q)}(M)$. If $f$
is defined in this way, then the mappings \mythetag{19.1} appear to be
inverse to each other and we get the equality $\pi_*\compos f=\idop$.
According to the definition~\mythedefinition{17.3}, it means that $f$ is a
horizontal lift of vectors. The theorem is completely proved.
\qed\enddemo
    Let's study a horizontal lift of vectors $f$ in a coordinate form. 
Upon choosing some local chart $U$ in $M$ equipped with a positively
polarized right orthonormal frame $\boldsymbol\Upsilon_0,\,
\boldsymbol\Upsilon_1,\,\boldsymbol\Upsilon_2,\,\boldsymbol\Upsilon_3$ 
and its associated spinor frame $\boldsymbol\Psi_1,\,\boldsymbol\Psi_2$ 
we can apply the lift $f$ to the frame vector fields $\boldsymbol\Upsilon_0,
\,\boldsymbol\Upsilon_1,\,\boldsymbol\Upsilon_2,
\,\boldsymbol\Upsilon_3$. As a result we get
$$
\gathered
f(\boldsymbol\Upsilon_j)=\bold U_j\,-\\
\vspace{2ex}
-\sum^{J+Q}_{P=1}
\dsize\msums{2}{3}\Sb i_1,\,\ldots,\,i_\alpha\\
j_1,\,\ldots,\,j_\beta\\ \bar i_1,\,\ldots,\,\bar i_\nu\\ 
\bar j_1,\,\ldots,\,\bar j_\gamma\\
h_1,\,\ldots,\,h_m\\ k_1,\,\ldots,\,k_n\endSb
\Gamma^{\,i_1\ldots\,i_\alpha\,\bar i_1\ldots\,\bar i_\nu\,h_1\ldots
\,h_m}_{\!j\,j_1\ldots\,j_\beta\,\bar j_1\ldots\,\bar j_\gamma\,k_1\ldots
\,k_n}[P]\ 
\bold W^{j_1\ldots\,j_\beta\,\bar j_1\ldots\,\bar j_\gamma\,k_1\ldots
\,k_n}_{i_1\ldots\,i_\alpha\,\bar i_1\ldots\,\bar i_\nu\,h_1\ldots
\,h_m}[P]\,-\\
\vspace{2ex}
-\sum^{J+Q}_{P=1}
\dsize\msums{2}{3}\Sb i_1,\,\ldots,\,i_\alpha\\
j_1,\,\ldots,\,j_\beta\\ \bar i_1,\,\ldots,\,\bar i_\nu\\ 
\bar j_1,\,\ldots,\,\bar j_\gamma\\
h_1,\,\ldots,\,h_m\\ k_1,\,\ldots,\,k_n\endSb
\bar\Gamma^{\,\bar i_1\ldots\,\bar i_\nu\,\,i_1\ldots
\,i_\alpha\,h_1\ldots\,h_m}_{\!j\,\bar j_1\ldots\,\bar j_\gamma\,j_1
\ldots\,j_\beta\,k_1\ldots\,k_n}[P]\ 
\bar\bold W^{\bar j_1\ldots\,\bar j_\gamma\,j_1\ldots\,j_\beta\,k_1
\ldots\,k_n}_{\bar i_1\ldots\,\bar i_\nu\,\,i_1\ldots\,i_\alpha\,h_1
\ldots\,h_m}[P].
\endgathered\quad
\mytag{19.3}
$$
Here $\alpha=\alpha_P$, $\beta=\beta_P$, $\nu=\nu_P$,
$\gamma=\gamma_P$, $m=m_P$, $n=n_P$, while the vectors $\bold U_j$,
$\bold W^{j_1\ldots\,j_\beta\,\bar j_1\ldots\,\bar j_\gamma\,k_1
\ldots\,k_n}_{i_1\ldots\,i_\alpha\,\bar i_1\ldots\,\bar i_\nu\,h_1
\ldots\,h_m}[P]$ and $\bar\bold W^{\bar j_1\ldots\,\bar j_\gamma\,j_1
\ldots\,j_\beta\,k_1\ldots\,k_n}_{\bar i_1\ldots\,\bar i_\nu\,\,i_1
\ldots\,i_\alpha\,h_1\ldots\,h_m}[P]$ are determined by the formulas 
\mythetag{9.26}, \mythetag{9.18}, \mythetag{9.19}, and \mythetag{9.20}.
The formula \mythetag{19.3} for $f(\bold E_j)$ follows from $\pi_*\compos
f=\idop$ due to the equalities 
$$
\align
&\pi_*\!\left(\bold U_j\right)=\bold E_j,\\
\vspace{1ex}
&\pi_*\!\left(\bold W^{j_1\ldots\,j_\beta\,\bar j_1\ldots\,\bar j_\gamma
\,k_1\ldots\,k_n}_{i_1\ldots\,i_\alpha\,\bar i_1\ldots\,\bar i_\nu\,h_1
\ldots\,h_m}[P]\right)=0,\\
\vspace{1ex}
&\pi_*\!\left(\bar\bold W^{\bar j_1\ldots\,\bar j_\gamma\,j_1\ldots\,
j_\beta\,k_1\ldots\,k_n}_{\bar i_1\ldots\,\bar i_\nu\,\,i_1\ldots\,
i_\alpha\,h_1\ldots\,h_m}[P]\right)=0.
\endalign
$$
The quantities $\Gamma^{\,i_1\ldots\,i_\alpha\,\bar i_1\ldots\,
\bar i_\nu\,h_1\ldots\,h_m}_{\!j\,j_1\ldots\,j_\beta\,\bar j_1\ldots
\,\bar j_\gamma\,k_1\ldots\,k_n}[P]$ and $\bar\Gamma^{\,\bar i_1\ldots
\,\bar i_\nu\,\,i_1\ldots\,i_\alpha\,h_1\ldots\,h_m}_{\!j\,\bar j_1\ldots
\,\bar j_\gamma\,j_1\ldots\,j_\beta\,k_1\ldots\,k_n}[P]$ in 
\mythetag{19.3} are called the {\it components} of a horizontal lift 
of vectors in a local chart $U$. If the lift $f$ is induced by some
horizontal covariant differentiation $\nabla$, then for its components
in the above formula \mythetag{19.3} we have
$$
\hskip -2em
\aligned
&\Gamma^{\,i_1\ldots\,i_\alpha\,\bar i_1\ldots\,\bar i_\nu\,h_1\ldots
\,h_m}_{\!j\,j_1\ldots\,j_\beta\,\bar j_1\ldots\,\bar j_\gamma\,k_1
\ldots\,k_n}[P]=-Z^{\,i_1\ldots\,i_\alpha\,\bar i_1\ldots\,\bar i_\nu
\,h_1\ldots\,h_m}_{j\,j_1\ldots\,j_\beta\,\bar j_1\ldots\,\bar j_\gamma
\,k_1\ldots\,k_n}[P],\\
\vspace{1em}
&\bar\Gamma^{\,\bar i_1\ldots\,\bar i_\nu\,\,i_1\ldots\,i_\alpha\,h_1
\ldots\,h_m}_{\!j\,\bar j_1\ldots\,\bar j_\gamma\,j_1\ldots\,j_\beta
\,k_1\ldots\,k_n}[P]=-\bar Z^{\,\bar i_1\ldots\,\bar i_\nu\,i_1\ldots
\,i_\alpha\,h_1\ldots\,h_m}_{j\,\bar j_1\ldots\,\bar j_\gamma\,j_1\ldots
\,j_\beta\,k_1\ldots\,k_n}[P].
\endaligned
\mytag{19.4}
$$
The quantities $Z^{\,i_1\ldots\,i_\alpha\,\bar i_1\ldots\,\bar i_\nu\,
h_1\ldots\,h_m}_{j\,j_1\ldots\,j_\beta\,\bar j_1\ldots\,\bar j_\gamma
\,k_1\ldots\,k_n}[P]$ and $\bar Z^{\,\bar i_1\ldots\,\bar i_\nu\,i_1\ldots
\,i_\alpha\,h_1\ldots\,h_m}_{j\,\bar j_1\ldots\,\bar j_\gamma\,j_1\ldots
\,j_\beta\,k_1\ldots\,k_n}[P]$ in \mythetag{19.4} are the same as in
\mythetag{15.6}, \mythetag{15.7}, \mythetag{15.11}, \mythetag{15.13},
and \mythetag{15.17}. As for the quantities $Z^i_j$ 
in \mythetag{15.5} and \mythetag{15.12}, in the case of a horizontal
covariant differentiation they are given by the Kronecker's delta-symbol:
$Z^i_j=\delta^i_j$. Substituting $Z^i_j=\tilde Z^i_j=\delta^i_j$ into
\mythetag{15.13} and \mythetag{15.17} and taking into account
\mythetag{19.4}, we derive the transformation rules for the components 
of a horizontal lift of spin-tensors in the formula \mythetag{19.3}:
$$
\gather
\gathered
\Gamma^{\,i_1\ldots\,i_\alpha\,\bar i_1\ldots\,\bar i_\nu\,h_1
\ldots\,h_m}_{j\,j_1\ldots\,j_\beta\,\bar j_1\ldots\,\bar j_\gamma
\,k_1\ldots\,k_n}[P]=\sum^3_{d=0}
\dsize\msums{2}{3}\Sb a_1,\,\ldots,\,a_\alpha\\ b_1,\,\ldots,\,b_\beta\\
\bar a_1,\,\ldots,\,\bar a_\nu\\ \bar b_1,\,\ldots,\,\bar b_\gamma\\
c_1,\,\ldots,\,c_m\\ d_1,\,\ldots,\,d_n\endSb
\goth S^{\,i_1}_{a_1}\ldots\,\goth S^{\,i_\alpha}_{a_\alpha}\ 
\goth T^{b_1}_{j_1}\ldots\,\goth T^{b_\beta}_{j_\beta}\ 
\overline{\goth S^{\,\bar i_1}_{\bar a_1}}
\ldots\,\overline{\goth S^{\,\bar i_\nu}_{\bar a_\nu}}\times\\
\vspace{1ex}
\times\,\overline{\goth T^{\,\bar b_1}_{\bar j_1}}
\ldots\,\overline{\goth T^{\,\bar b_\gamma}_{\bar j_\gamma}}\ 
S^{h_1}_{c_1}\ldots\,S^{h_m}_{c_m}\
T^{\,d_1}_{k_1}\ldots\,T^{\,d_n}_{k_n}\ T^d_j\ 
\tilde\Gamma^{\,a_1\ldots\,a_\alpha\,\bar a_1\ldots\,
\bar a_\nu\,c_1\ldots\,c_m}_{d\,b_1\ldots\,b_\beta\,\bar b_1\ldots\,
\bar b_\gamma\,d_1\ldots\,d_n}[P]\,+\\
\vspace{1ex}
+\sum^\alpha_{\mu=1}\sum^2_{v_\mu=1}\vartheta^{\,i_\mu}_{j\,v_\mu}
\ S^{\,i_1\ldots\,v_\mu\ldots\,i_\alpha\,\bar i_1\ldots\,\bar i_\nu\,
h_1\ldots\,h_m}_{\,j_1\ldots\,\ldots\,\ldots\,j_\beta\,\bar j_1\ldots\,
\bar j_\gamma\,k_1\ldots\,k_n}[P]-\sum^\beta_{\mu=1}\sum^2_{w_\mu=1}
\vartheta^{\,w_\mu}_{j\,j_\mu}\,\times\\
\times\,S^{\,i_1\ldots\,\ldots\,\ldots\,i_\alpha\,\bar i_1\ldots\,
\bar i_\nu\,h_1\ldots\,h_m}_{\,j_1\ldots\,w_\mu\,\ldots\,j_\beta\,
\bar j_1\ldots\,\bar j_\gamma\,k_1\ldots\,k_n}[P]
+\sum^\nu_{\mu=1}\sum^2_{v_\mu=1}\overline{\vartheta^{\,
\bar i_\mu}_{j\,v_\mu}}
\ S^{\,i_1\ldots\,i_\alpha\,\bar i_1\ldots\,v_\mu\ldots\,\bar i_\nu\,
h_1\ldots\,h_m}_{\,j_1\ldots\,j_\beta\,\bar j_1\ldots\,\ldots\,\ldots\,
\bar j_\gamma\,k_1\ldots\,k_n}[P]\,-\\
-\sum^\gamma_{\mu=1}\sum^2_{w_\mu=1}
\overline{\vartheta^{\,w_\mu}_{j\,\bar j_\mu}}
\ S^{\,i_1\ldots\,i_\alpha\,\bar i_1\ldots\,\ldots\,\ldots\,
\bar i_\nu\,h_1\ldots\,h_m}_{\,j_1\ldots\,j_\beta\,\bar j_1\ldots\,
w_\mu\,\ldots\,\bar j_\gamma\,k_1\ldots\,k_n}[P]
+\sum^m_{\mu=1}\sum^3_{v_\mu=0}\theta^{\,i_\mu}_{j\,v_\mu}\,\times\\
\times\,S^{\,i_1\ldots\,i_\alpha\,\bar i_1\ldots\,\bar i_\nu\,
h_1\ldots\,v_\mu\ldots\,h_m}_{\,j_1\ldots\,j_\beta\,\bar j_1\ldots\,
\bar j_\gamma\,k_1\ldots\,\ldots\,\ldots\,k_n}[P]
-\sum^n_{\mu=1}\sum^3_{w_\mu=0}\theta^{\,w_\mu}_{j\,j_\mu}
\ S^{\,i_1\ldots\,i_\alpha\,\bar i_1\ldots\,\bar i_\nu\,
h_1\ldots\,\ldots\,\ldots\,h_m}_{\,j_1\ldots\,j_\beta\,\bar j_1\ldots\,
\bar j_\gamma\,k_1\ldots\,w_\mu\,\ldots\,k_n}[P],
\endgathered\\
\vspace{2ex}
\gathered
\bar\Gamma^{\,\bar i_1\ldots\,\bar i_\nu\,i_1\ldots\,i_\alpha\,h_1\ldots
\,h_m}_{j\,\bar j_1\ldots\,\bar j_\gamma\,j_1\ldots\,j_\beta\,k_1
\ldots\,k_n}[P]=\sum^3_{d=0}
\dsize\msums{2}{3}\Sb a_1,\,\ldots,\,a_\alpha\\ b_1,\,\ldots,\,b_\beta\\
\bar a_1,\,\ldots,\,\bar a_\nu\\ \bar b_1,\,\ldots,\,\bar b_\gamma\\
c_1,\,\ldots,\,c_m\\ d_1,\,\ldots,\,d_n\endSb
\goth S^{\,\bar i_1}_{\bar a_1}\ldots\,
\goth S^{\,\bar i_\nu}_{\bar a_\nu}\ \goth T^{\bar b_1}_{\bar j_1}\ldots\,
\goth T^{\bar b_\gamma}_{\bar j_\gamma}\ \overline{\goth S^{\,i_1}_{a_1}}
\ldots\,\overline{\goth S^{\,i_\alpha}_{a_\alpha}}\,\times\\
\vspace{1ex}
\times\,\overline{\goth T^{\,b_1}_{j_1}}
\ldots\,\overline{\goth T^{\,b_\beta}_{j_\beta}}\ 
S^{h_1}_{c_1}\ldots\,S^{h_m}_{c_m}\,
T^{\,d_1}_{k_1}\ldots\,T^{\,d_n}_{k_n}\ T^d_j\ 
\tbGamma^{\bar a_1\ldots\,\bar a_\nu\,a_1\ldots\,a_\alpha\,c_1
\ldots\,c_m}_{d\,\bar b_1\ldots\,\bar b_\gamma\,b_1\ldots\,b_\beta
\,d_1\ldots\,d_n}[P]\,+\\
\vspace{1ex}
+\sum^\alpha_{\mu=1}\sum^2_{v_\mu=1}
\overline{\vartheta^{\,i_\mu}_{j\,v_\mu}}
\ \overline{S^{\,i_1\ldots\,v_\mu\ldots\,i_\alpha\,\bar i_1
\ldots\,\bar i_\nu\,h_1\ldots\,h_m}_{\,j_1\ldots\,\ldots\,\ldots\,
j_\beta\,\bar j_1\ldots\,\bar j_\gamma\,k_1\ldots\,k_n}[P]}
-\sum^\beta_{\mu=1}\sum^2_{w_\mu=1}
\overline{\vartheta^{\,w_\mu}_{j\,j_\mu}}\,\times\\
\times\ \overline{S^{\,i_1\ldots\,\ldots\,\ldots\,i_\alpha\,\bar i_1
\ldots\,\bar i_\nu\,h_1\ldots\,h_m}_{\,j_1\ldots\,w_\mu\,\ldots\,j_\beta
\,\bar j_1\ldots\,\bar j_\gamma\,k_1\ldots\,k_n}[P]} 
+\sum^\nu_{\mu=1}\sum^2_{v_\mu=1}\vartheta^{\,\bar i_\mu}_{j\,v_\mu}
\ \overline{S^{\,i_1\ldots\,i_\alpha\,\bar i_1\ldots\,v_\mu
\ldots\,\bar i_\nu\,h_1\ldots\,h_m}_{\,j_1\ldots\,j_\beta\,\bar j_1
\ldots\,\ldots\,\ldots\,\bar j_\gamma\,k_1\ldots\,k_n}[P]}\,-\\
-\sum^\gamma_{\mu=1}\sum^2_{w_\mu=1}\vartheta^{\,w_\mu}_{j\,\bar j_\mu}
\ \overline{S^{\,i_1\ldots\,i_\alpha\,\bar i_1\ldots\,\ldots\,
\ldots\,\bar i_\nu\,h_1\ldots\,h_m}_{\,j_1\ldots\,j_\beta\,\bar j_1
\ldots\,w_\mu\,\ldots\,\bar j_\gamma\,k_1\ldots\,k_n}[P]}
+\sum^m_{\mu=1}\sum^3_{v_\mu=0}\theta^{\,i_\mu}_{j\,v_\mu}\,\times\\
\times\ \overline{S^{\,i_1\ldots\,i_\alpha\,\bar i_1\ldots\,\bar i_\nu
\,h_1\ldots\,v_\mu\ldots\,h_m}_{\,j_1\ldots\,j_\beta\,\bar j_1\ldots\,
\bar j_\gamma\,k_1\ldots\,\ldots\,\ldots\,k_n}[P]}
-\sum^n_{\mu=1}\sum^3_{w_\mu=0}\theta^{\,w_\mu}_{j\,j_\mu}
\ \overline{S^{\,i_1\ldots\,i_\alpha\,\bar i_1\ldots\,\bar i_\nu
\,h_1\ldots\,\ldots\,\ldots\,h_m}_{\,j_1\ldots\,j_\beta\,\bar j_1\ldots
\,\bar j_\gamma\,k_1\ldots\,w_\mu\,\ldots\,k_n}[P]}.
\endgathered
\endgather
$$\par
     Other geometric structures associated with a horizontal covariant
differentiation $\nabla$ reveal when we apply $D=\nabla_{\bold Y}$ to
the modules $S^1_0\bar S^0_0 T^0_0(M)$, $S^0_0\bar S^1_0 T^0_0(M)$, and 
$S^0_0\bar S^0_0 T^1_0(M)$. The action of $D=\nabla_{\bold Y}$ within 
these modules is determined by the formulas \mythetag{13.23} and
\mythetag{13.22}, where $\hat{\boldsymbol\Upsilon}_i$, 
$\hat{\boldsymbol\Psi}_i$, and $\hat{\bPsi}_i$ are defined by the
formulas \mythetag{13.6}, \mythetag{13.8}, and \mythetag{13.10}
respectively. In the present case we can take 
$\bold Y=\hat{\boldsymbol\Upsilon}_j$ and write these formulas 
as follows:
$$
\align
&\hskip -2em
\nabla_{\hat{\boldsymbol\Upsilon}_j}\!\hat{\boldsymbol\Psi}_i
=\sum^2_{k=1}\Alpha^k_{j\,i}\,\boldsymbol\Psi_k,
\mytag{19.5}\\
&\hskip -2em
\nabla_{\hat{\boldsymbol\Upsilon}_j}\!\hat{\bPsi}_i
=\sum^2_{k=1}\bar{\Alpha}\vphantom{\Alpha}^k_{j\,i}\,\bPsi_k.
\mytag{19.6}\\
&\hskip -2em
\nabla_{\hat{\boldsymbol\Upsilon}_j}\!
\hat{\boldsymbol\Upsilon}_i=\sum^3_{k=0}\Gamma^k_{j\,i}\,\bold E_k.
\mytag{19.7}
\endalign
$$
The coefficients $\Alpha^k_{j\,i}$,
$\bar{\Alpha}\vphantom{\Alpha}^k_{j\,i}$, and $\Gamma^k_{j\,i}$
in \mythetag{19.5}, \mythetag{19.6}, and \mythetag{19.7} are the same 
as in the formulas \mythetag{15.9}, \mythetag{15.10}, and \mythetag{15.8}.
Since $Z^i_j=\delta^i_j$ for a horizontal covariant differentiation,
the transformation formulas \mythetag{15.15}, \mythetag{15.16}, and
\mythetag{15.14} for $\Alpha^k_{j\,i}$,
$\bar{\Alpha}\vphantom{\Alpha}^k_{j\,i}$, and $\Gamma^k_{j\,i}$
now reduce to the following ones:
$$
\align
&\hskip -2em
\Alpha^k_{j\,i}=\dsize\sum^2_{b=1}\sum^2_{a=1}\sum^3_{c=0}
\goth S^k_a\,\goth T^b_i\,T^c_j\ \tilde{\Alpha}\vphantom{\Alpha}^a_{c\,b}
+\vartheta^k_{j\,i},
\mytag{19.8}\\
&\hskip -2em
\bar{\Alpha}\vphantom{\Alpha}^k_{j\,i}=\sum^2_{b=1}\sum^2_{a=1}
\sum^3_{c=0}\overline{\goth S^k_a}\ \overline{\goth T^b_i}\,T^c_j\  
\tilde{\bar{\Alpha}}\vphantom{\Alpha}^a_b
+\overline{\vartheta^k_{j\,i}}.
\mytag{19.9}\\
\vspace{2ex}
&\hskip -2em
\Gamma^k_{j\,i}=\dsize\sum^3_{b=0}\sum^3_{a=0}\sum^3_{c=0}
S^k_a\,T^b_i\,T^c_j\ \tilde\Gamma^a_{c\,b}+\theta^k_{j\,i},
\mytag{19.10}
\endalign
$$
\mydefinition{19.1} Let $N$ be a composite spin-tensorial bundle
over the space-time manifold $M$. An extended affine connection 
is a geometric object such that in each local chart $U$ of $M$ 
equipped with a positively polarized right orthonormal frame
$\boldsymbol\Upsilon_0,\,\boldsymbol\Upsilon_1,\,\boldsymbol\Upsilon_2,
\,\boldsymbol\Upsilon_3$ and associated spinor frame 
$\boldsymbol\Psi_1,\,\boldsymbol\Psi_2$ it is represented by its
components $\Alpha^k_{j\,i}$, $\bar{\Alpha}\vphantom{\Alpha}^k_{j\,i}$,
$\Gamma^k_{j\,i}$ and such that its components are smooth functions
of the variables \mythetag{9.11} and \mythetag{9.12} transforming 
according to the formulas \mythetag{19.8}, \mythetag{19.9}, and
\mythetag{19.10} under a change of a local chart.
\enddefinition
\mytheorem{19.2} Any smooth paracompact space-time manifold $M$ 
admitting the spinor structure and equipped with a composite 
spin-tensorial bundle $N$ possesses at least one extended affine 
connection.
\endproclaim
     Compare the theorem~\mythetheorem{19.2} with the theorem~13.2
in \mycite{5}. See also Chapter~\uppercase\expandafter{\romannumeral 3} 
of the thesis \mycite{12} for ideas how to prove this theorem.
\mydefinition{19.2} A horizontal covariant differentiation $\nabla$ 
of the algebra of extended spin-tensorial fields $\bold S(M)$ is called 
a {\it spatial covariant differentiation} or a {\it spatial gradient} if
the operator $\nabla$ annuls all native extended spin-tensorial fields
and all complex conjugate fields for them:
$$
\nabla\bold S[P]=0\text{ \ \ and \ \ }\nabla\tau(\bold S[P])=0
\text{ \ \ for all \ \ }P=1,\,\ldots,\,J+Q.\quad
\mytag{19.11}
$$
\enddefinition
    Let's study the first equality \mythetag{19.11} in local coordinates.
For this purpose we use the formula \mythetag{15.11} substituting 
$Z^i_j=\delta^i_j$ and taking into account \mythetag{19.4}. As a result
from $\nabla_{\!\!j}S^{\,i_1\ldots\,i_\alpha\,\bar i_1\ldots\,
\bar i_\nu\,h_1\ldots\,h_m}_{j_1\ldots\,j_\beta\,\bar j_1\ldots\,
\bar j_\gamma\,k_1\ldots\,k_n}[P]=0$ we derive the equality
$$
\gather
\Gamma^{\,i_1\ldots\,i_\alpha\,\bar i_1\ldots\,\bar i_\nu\,h_1\ldots
\,h_m}_{j\,j_1\ldots\,j_\beta\,\bar j_1\ldots\,\bar j_\gamma\,k_1\ldots
\,k_n}[P]=
\sum^\alpha_{\mu=1}\sum^2_{v_\mu=1}\Alpha^{i_\mu}_{j\,v_\mu}\ 
S^{\,i_1\ldots\,v_\mu\,\ldots\,i_\alpha\,\bar i_1\ldots\,\bar i_\nu
\,h_1\ldots\,h_m}_{j_1\ldots\,\ldots\,\ldots\,j_\beta\,\bar j_1\ldots\,
\bar j_\gamma\,k_1\ldots\,k_n}[P]\,-\\
-\sum^\beta_{\mu=1}\sum^2_{w_\mu=1}\Alpha^{w_\mu}_{j\,j_\mu}\
S^{\,i_1\ldots\,\ldots\,\ldots\,i_\alpha\,\bar i_1\ldots\,\bar i_\nu
\,h_1\ldots\,h_m}_{j_1\ldots\,w_\mu\,\ldots\,j_\beta\,\bar j_1\ldots\,
\bar j_\gamma\,k_1\ldots\,k_n}[P]
+\sum^\nu_{\mu=1}\sum^2_{v_\mu=1}
\bar{\Alpha}\vphantom{\Alpha}^{\bar i_\mu}_{j\,v_\mu}\,\times\\
\times\,S^{\,i_1\ldots\,i_\alpha\,\bar i_1\ldots\,v_\mu\,\ldots\,
\bar i_\nu\,h_1\ldots\,h_m}_{j_1\ldots\,j_\beta\,\bar j_1\ldots\,\ldots
\,\ldots\,\bar j_\gamma\,k_1\ldots\,k_n}[P]-
\sum^\gamma_{\mu=1}\sum^2_{w_\mu=1}
\bar{\Alpha}\vphantom{\Alpha}^{w_\mu}_{j\,\bar j_\mu}\
S^{i_1\ldots\,i_\alpha\,\bar i_1\ldots\,\ldots\,\ldots\,\bar i_\nu
\,h_1\ldots\,h_m}_{j_1\ldots\,j_\beta\,\bar j_1\ldots\,w_\mu\,\ldots\,
\bar j_\gamma\,k_1\ldots\,k_n}[P]\,+\\
+\sum^m_{\mu=1}\sum^3_{v_\mu=0}\Gamma^{h_\mu}_{j\,v_\mu}\ 
S^{i_1\ldots\,i_\alpha\,\bar i_1\ldots\,\bar i_\nu
\,h_1\ldots\,v_\mu\,\ldots\,h_m}_{j_1\ldots\,j_\beta\,\bar j_1\ldots\,
\bar j_\gamma\,k_1\ldots\,\ldots\,\ldots\,k_n}[P]-
\sum^n_{\mu=1}\sum^3_{w_\mu=0}\Gamma^{w_\mu}_{j\,k_\mu}\,\times\\
\times\,S^{i_1\ldots\,i_\alpha\,\bar i_1\ldots\,\bar i_\nu\,
h_1\ldots\,\ldots\,\ldots\,h_m}_{j_1\ldots\,j_\beta\,\bar j_1\ldots\,
\bar j_\gamma\,k_1\ldots\,w_\mu\,\ldots\,k_n}[P]
\vphantom{\sum^n_{\mu=1}\sum^3_{w_\mu=0}}
\endgather
$$
expressing $\Gamma^{\,i_1\ldots\,i_\alpha\,\bar i_1\ldots\,\bar i_\nu
\,h_1\ldots\,h_m}_{j\,j_1\ldots\,j_\beta\,\bar j_1\ldots\,\bar j_\gamma
\,k_1\ldots\,k_n}[P]$ through the parameters $\Alpha^k_{j\,i}$,
$\bar{\Alpha}\vphantom{\Alpha}^k_{j\,i}$, and $\Gamma^k_{j\,i}$.
Similarly, from the second equality \mythetag{19.11}, using the formula 
\mythetag{4.6}, we derive the equality  $\nabla_{\!\!j}\overline{S^{\,
\bar i_1\ldots\,\bar i_\alpha\,i_1\ldots\,i_\nu\,h_1\ldots\,
h_m}_{\bar j_1\ldots\,\bar j_\beta\,j_1\ldots\,j_\gamma\,k_1\ldots\,
k_n}[P]}=0$. Applying \mythetag{15.11} to it, we get
$$
\gather
\bar\Gamma^{\,i_1\ldots\,i_\nu\,\bar i_1\ldots\,\bar i_\alpha\,h_1\ldots
\,h_m}_{j\,j_1\ldots\,j_\gamma\,\bar j_1\ldots\,\bar j_\beta\,k_1\ldots
\,k_n}[P]=\sum^\nu_{\mu=1}\sum^2_{v_\mu=1}\Alpha^{i_\mu}_{j\,v_\mu}
\ \overline{S^{\,\bar i_1\ldots\,\bar i_\alpha\,i_1\ldots\,v_\mu\,\ldots\,
i_\nu h_1\ldots\,h_m}_{\bar j_1\ldots\,\bar j_\beta\,j_1\ldots\,\ldots\,
\ldots\,j_\gamma k_1\ldots\,k_n}[P]}\,-\\
-\sum^\gamma_{\mu=1}\sum^2_{w_\mu=1}\Alpha^{w_\mu}_{j\,j_\mu}
\ \overline{S^{\,\bar i_1\ldots\,\bar i_\alpha\,i_1\ldots\,\ldots\,\ldots
\,i_\nu\,h_1\ldots\,h_m}_{\bar j_1\ldots\,\bar j_\beta\,j_1\ldots\,w_\mu
\,\ldots\,j_\gamma\,k_1\ldots\,k_n}[P]}+\sum^\alpha_{\mu=1}\sum^2_{v_\mu=1}
\bar{\Alpha}\vphantom{\Alpha}^{\bar i_\mu}_{j\,v_\mu}\,\times\\
\times\,\overline{S^{\,\bar i_1\ldots\,v_\mu\,\ldots\,\bar i_\alpha\,i_1
\ldots\,i_\nu\,h_1\ldots\,h_m}_{\bar j_1\ldots\,\ldots\,\ldots\,
\bar j_\beta\,j_1\ldots\,j_\gamma\,k_1\ldots\,k_n}[P]}
-\sum^\beta_{\mu=1}\sum^2_{w_\mu=1}
\bar{\Alpha}\vphantom{\Alpha}^{w_\mu}_{j\,\bar j_\mu}
\ \overline{S^{\,\bar i_1\ldots\,\ldots\,\ldots\,\bar i_\alpha\,
i_1\ldots\,i_\nu\,h_1\ldots\,h_m}_{\bar j_1\ldots\,w_\mu\,\ldots\,
\bar j_\beta\,j_1\ldots\,j_\gamma\,k_1\ldots\,k_n}[P]}\,+\\
+\sum^m_{\mu=1}\sum^3_{v_\mu=0}\Gamma^{h_\mu}_{j\,v_\mu}\ 
\overline{S^{\bar i_1\ldots\,\bar i_\alpha\,i_1\ldots\,i_\nu\,
h_1\ldots\,v_\mu\,\ldots\,h_m}_{\bar j_1\ldots\,\bar j_\beta\,j_1
\ldots\,j_\gamma\,k_1\ldots\,\ldots\,\ldots\,k_n}[P]}-
\sum^n_{\mu=1}\sum^3_{w_\mu=0}\Gamma^{w_\mu}_{j\,k_\mu}\,\times\\
\times\,\overline{S^{\bar i_1\ldots\,\bar i_\alpha\,i_1\ldots\,i_\nu
\,h_1\ldots\,\ldots\,\ldots\,h_m}_{\bar j_1\ldots\,\bar j_\beta
j_1\ldots\,j_\gamma\,k_1\ldots\,w_\mu\,\ldots\,k_n}[P]}.
\vphantom{\sum^n_{\mu=1}\sum^3_{w_\mu=0}}
\endgather
$$
This equality expresses $\bar\Gamma^{\,i_1\ldots\,i_\nu\,\bar i_1\ldots
\,\bar i_\alpha\,h_1\ldots\,h_m}_{j\,j_1\ldots\,j_\gamma\,\bar j_1\ldots
\,\bar j_\beta\,k_1\ldots\,k_n}[P]$ through the same parameters 
$\Alpha^k_{j\,i}$, $\bar{\Alpha}\vphantom{\Alpha}^k_{j\,i}$, and $\Gamma^k_{j\,i}$. Thus, due to the additional condition \mythetag{19.11} 
a spacial covariant differentiation depends on less number of parameters
than an arbitrary horizontal covariant differentiation. Substituting the above two expressions back into the formula \mythetag{5.11} and taking 
into account \mythetag{19.4} and $Z^i_j=\delta^i_j$, we derive the general formula for a spacial covariant derivative in local coordinates \pagebreak
\mythetag{9.11} and \mythetag{9.12}. Though being more huge, this formula is a specialization of the formula \mythetag{15.11}:
$$
\allowdisplaybreaks
\gather
\gathered
\nabla_{\!\!j}X^{a_1\ldots\,a_\varepsilon\bar a_1\ldots\,\bar a_\sigma c_1\ldots\,c_e}_{b_1\ldots\,b_\eta\bar b_1\ldots\,\bar b_\zeta
d_1\ldots\,d_f}
=\sum^3_{k=0}\Upsilon^k_j\,
\frac{\partial X^{a_1\ldots\,a_\varepsilon\bar a_1\ldots\,\bar a_\sigma c_1\ldots\,c_e}_{b_1\ldots\,b_\eta\bar b_1\ldots\,\bar b_\zeta d_1\ldots\,
d_f}}{\partial x^k}\,-\\
\vspace{1.5ex}
-\sum^{J+Q}_{P=1}\dsize\msums{2}{3}\Sb i_1,\,\ldots,\,i_\alpha\\
j_1,\,\ldots,\,j_\beta\\ \bar i_1,\,\ldots,\,\bar i_\nu\\ 
\bar j_1,\,\ldots,\,\bar j_\gamma\\
h_1,\,\ldots,\,h_m\\ k_1,\,\ldots,\,k_n\endSb
\left(\,\shave{\sum^\alpha_{\mu=1}\sum^2_{v_\mu=1}}\Alpha^{i_\mu}_{j\,v_\mu}
\ S^{\,i_1\ldots\,v_\mu\,\ldots\,i_\alpha\,\bar i_1\ldots\,\bar i_\nu
\,h_1\ldots\,h_m}_{j_1\ldots\,\ldots\,\ldots\,j_\beta\,\bar j_1\ldots\,
\bar j_\gamma\,k_1\ldots\,k_n}[P]\,-\right.\\
-\sum^\beta_{\mu=1}\sum^2_{w_\mu=1}\Alpha^{w_\mu}_{j\,j_\mu}
\ S^{\,i_1\ldots\,\ldots\,\ldots\,i_\alpha\,\bar i_1\ldots\,\bar i_\nu
\,h_1\ldots\,h_m}_{j_1\ldots\,w_\mu\,\ldots\,j_\beta\,\bar j_1\ldots\,
\bar j_\gamma\,k_1\ldots\,k_n}[P]
+\sum^\nu_{\mu=1}\sum^2_{v_\mu=1}
\bar{\Alpha}\vphantom{\Alpha}^{\bar i_\mu}_{j\,v_\mu}\,\times\\
\times\,S^{\,i_1\ldots\,i_\alpha\,\bar i_1\ldots\,v_\mu\,\ldots\,
\bar i_\nu\,h_1\ldots\,h_m}_{j_1\ldots\,j_\beta\,\bar j_1\ldots\,\ldots
\,\ldots\,\bar j_\gamma\,k_1\ldots\,k_n}[P]-
\sum^\gamma_{\mu=1}\sum^2_{w_\mu=1}
\bar{\Alpha}\vphantom{\Alpha}^{w_\mu}_{j\,\bar j_\mu}\,\times\kern 4em\\
\times\,S^{i_1\ldots\,i_\alpha\,\bar i_1\ldots\,\ldots\,\ldots\,\bar i_\nu
\,h_1\ldots\,h_m}_{j_1\ldots\,j_\beta\,\bar j_1\ldots\,w_\mu\,\ldots\,
\bar j_\gamma\,k_1\ldots\,k_n}[P]
+\sum^m_{\mu=1}\sum^3_{v_\mu=0}\Gamma^{h_\mu}_{j\,v_\mu}\,\times\\
\times\,S^{i_1\ldots\,i_\alpha\,\bar i_1\ldots\,\bar i_\nu
\,h_1\ldots\,v_\mu\,\ldots\,h_m}_{j_1\ldots\,j_\beta\,\bar j_1\ldots\,
\bar j_\gamma\,k_1\ldots\,\ldots\,\ldots\,k_n}[P]-
\sum^n_{\mu=1}\sum^3_{w_\mu=0}\Gamma^{w_\mu}_{j\,k_\mu}\,\times\kern 4em\\
\left.\vphantom{\shave{\sum^\alpha_{\mu=1}\sum^2_{v_\mu=1}}}
\times\,S^{i_1\ldots\,i_\alpha\,\bar i_1\ldots\,\bar i_\nu\,
h_1\ldots\,\ldots\,\ldots\,h_m}_{j_1\ldots\,j_\beta\,\bar j_1\ldots\,
\bar j_\gamma\,k_1\ldots\,w_\mu\,\ldots\,k_n}[P]\right)
\frac{\partial X^{a_1\ldots\,a_\varepsilon\bar a_1\ldots\,\bar a_\sigma c_1\ldots\,c_e}_{b_1\ldots\,b_\eta\bar b_1\ldots\,\bar b_\zeta d_1\ldots\,
d_f}}{\partial S^{\,i_1\ldots\,i_\alpha\,\bar i_1\ldots\,\bar i_\nu\,
h_1\ldots\,h_m}_{j_1\ldots\,j_\beta\,\bar j_1\ldots\,\bar j_\gamma\,
k_1\ldots\,k_n}[P]}\,-\\
\vspace{1.5ex}
-\sum^{J+Q}_{P=1}\dsize\msums{2}{3}\Sb i_1,\,\ldots,\,i_\alpha\\
j_1,\,\ldots,\,j_\beta\\ \bar i_1,\,\ldots,\,\bar i_\nu\\ 
\bar j_1,\,\ldots,\,\bar j_\gamma\\
h_1,\,\ldots,\,h_m\\ k_1,\,\ldots,\,k_n\endSb
\left(\shave{\sum^\nu_{\mu=1}\sum^2_{v_\mu=1}}\Alpha^{i_\mu}_{j\,v_\mu}
\ \overline{S^{\,\bar i_1\ldots\,\bar i_\alpha\,i_1\ldots\,v_\mu\,\ldots\,
i_\nu h_1\ldots\,h_m}_{\bar j_1\ldots\,\bar j_\beta\,j_1\ldots\,\ldots\,
\ldots\,j_\gamma k_1\ldots\,k_n}[P]}\,-\right.\\
-\sum^\gamma_{\mu=1}\sum^2_{w_\mu=1}\Alpha^{w_\mu}_{j\,j_\mu}
\ \overline{S^{\,\bar i_1\ldots\,\bar i_\alpha\,i_1\ldots\,\ldots\,\ldots
\,i_\nu\,h_1\ldots\,h_m}_{\bar j_1\ldots\,\bar j_\beta\,j_1\ldots\,w_\mu
\,\ldots\,j_\gamma\,k_1\ldots\,k_n}[P]}+\sum^\alpha_{\mu=1}\sum^2_{v_\mu=1}
\bar{\Alpha}\vphantom{\Alpha}^{\bar i_\mu}_{j\,v_\mu}\,\times\\
\times\,\overline{S^{\,\bar i_1\ldots\,v_\mu\,\ldots\,\bar i_\alpha\,i_1
\ldots\,i_\nu\,h_1\ldots\,h_m}_{\bar j_1\ldots\,\ldots\,\ldots\,
\bar j_\beta\,j_1\ldots\,j_\gamma\,k_1\ldots\,k_n}[P]}
-\sum^\beta_{\mu=1}\sum^2_{w_\mu=1}
\bar{\Alpha}\vphantom{\Alpha}^{w_\mu}_{j\,\bar j_\mu}\,\times\\
\times\,\overline{S^{\,\bar i_1\ldots\,\ldots\,\ldots\,\bar i_\alpha\,
i_1\ldots\,i_\nu\,h_1\ldots\,h_m}_{\bar j_1\ldots\,w_\mu\,\ldots\,
\bar j_\beta\,j_1\ldots\,j_\gamma\,k_1\ldots\,k_n}[P]}+
\sum^m_{\mu=1}\sum^3_{v_\mu=0}\Gamma^{h_\mu}_{j\,v_\mu}\,\times\\
\times\,\overline{S^{\bar i_1\ldots\,\bar i_\alpha\,i_1\ldots\,i_\nu\,
h_1\ldots\,v_\mu\,\ldots\,h_m}_{\bar j_1\ldots\,\bar j_\beta\,j_1
\ldots\,j_\gamma\,k_1\ldots\,\ldots\,\ldots\,k_n}[P]}-
\sum^n_{\mu=1}\sum^3_{w_\mu=0}\Gamma^{w_\mu}_{j\,k_\mu}\,\times\\
\left.\vphantom{\shave{\sum^\nu_{\mu=1}\sum^2_{v_\mu=1}}}
\times\,\overline{S^{\bar i_1\ldots\,\bar i_\alpha\,i_1\ldots\,i_\nu
\,h_1\ldots\,\ldots\,\ldots\,h_m}_{\bar j_1\ldots\,\bar j_\beta
j_1\ldots\,j_\gamma\,k_1\ldots\,w_\mu\,\ldots\,k_n}[P]}
\right)\frac{\partial X^{a_1\ldots\,a_\varepsilon\bar a_1\ldots\,
\bar a_\sigma c_1\ldots\,c_e}_{b_1\ldots\,b_\eta\bar b_1\ldots\,
\bar b_\zeta d_1\ldots\,d_f}}{\partial\overline{S^{\,i_1\ldots\,
i_\alpha\,\bar i_1\ldots\,\bar i_\nu\,h_1\ldots\,h_m}_{j_1\ldots\,
j_\beta\,\bar j_1\ldots\,\bar j_\gamma\,k_1\ldots\,k_n}[P]}\,}\,+
\endgathered
\mytag{19.12}\\
\gathered
\kern -9em
+\sum^\varepsilon_{\mu=1}\sum^2_{v_\mu=1}\Alpha^{a_\mu}_{j\,v_\mu}\ 
X^{a_1\ldots\,v_\mu\,\ldots\,a_\varepsilon\bar a_1\ldots\,\bar a_\sigma
c_1\ldots\,c_e}_{b_1\ldots\,\ldots\,\ldots\,b_\eta\bar b_1\ldots\,
\bar b_\zeta d_1\ldots\,d_f}\,-\\
\kern 9em-\sum^\eta_{\mu=1}\sum^2_{w_\mu=1}\Alpha^{w_\mu}_{j\,b_\mu}\
X^{a_1\ldots\,\ldots\,\ldots\,a_\varepsilon\bar a_1\ldots\,\bar a_\sigma
c_1\ldots\,c_e}_{b_1\ldots\,w_\mu\,\ldots\,b_\eta\bar b_1\ldots\,
\bar b_\zeta d_1\ldots\,d_f}\,+\\
\kern -9em
+\sum^\sigma_{\mu=1}\sum^2_{v_\mu=1}
\bar{\Alpha}\vphantom{\Alpha}^{\bar a_\mu}_{j\,v_\mu}\ 
X^{a_1\ldots\,a_\varepsilon\bar a_1\ldots\,v_\mu\,\ldots\,\bar a_\sigma
c_1\ldots\,c_e}_{b_1\ldots\,b_\eta\bar b_1\ldots\,\ldots\,\ldots\,
\bar b_\zeta d_1\ldots\,d_f}\,-\\
\kern 9em-\sum^\zeta_{\mu=1}\sum^2_{w_\mu=1}
\bar{\Alpha}\vphantom{\Alpha}^{w_\mu}_{j\,\bar b_\mu}\
X^{a_1\ldots\,a_\varepsilon\bar a_1\ldots\,\ldots\,\ldots\,\bar a_\sigma
c_1\ldots\,c_e}_{b_1\ldots\,b_\eta\bar b_1\ldots\,w_\mu\,\ldots\,
\bar b_\zeta d_1\ldots\,d_f}\,+\\
\kern -9em+\sum^e_{\mu=1}\sum^3_{v_\mu=0}\Gamma^{c_\mu}_{j\,v_\mu}\ 
X^{a_1\ldots\,a_\varepsilon\bar a_1\ldots\,\bar a_\sigma
c_1\ldots\,v_\mu\,\ldots\,c_e}_{b_1\ldots\,b_\eta\bar b_1\ldots\,
\bar b_\zeta d_1\ldots\,\ldots\,\ldots\,d_f}\,-\\
\kern 9em-\sum^f_{\mu=1}\sum^3_{w_\mu=0}\Gamma^{w_\mu}_{j\,b_\mu}\
X^{a_1\ldots\,a_\varepsilon\bar a_1\ldots\,\bar a_\sigma
c_1\ldots\,\ldots\,\ldots\,c_e}_{b_1\ldots\,b_\eta\bar b_1\ldots\,
\bar b_\zeta d_1\ldots\,w_\mu\,\ldots\,d_f}.
\endgathered\kern 4em
\endgather
$$
Any horizontal covariant differentiation differentiation is defined by
two independent geometric structures: 
\roster
\item a horizontal lift of vectors from $M$ to $N$;
\item an extended connection.
\endroster
In the case of a spacial differentiation these two structures are related
to each other. This result is formulated as the following theorem.
\mytheorem{19.3} Defining a spacial covariant differentiation in the
algebra of extended spin-tensorial fields $\bold S(M)$ is equivalent 
to defining an extended connection.
\endproclaim
\head
20. The structural theorem for differentiations.
\endhead
\mytheorem{20.1} Let $N$ be a composite tensor bundle over the 
space-time manifold $M$ in the sense of the formula \mythetag{9.7}. 
If\/ $M$ is equipped with some extended affine connection, then each 
differentiation $D$ of the algebra of extended spin-tensorial fields 
$\bold S(M)$ is uniquely expanded into a sum
$$
\hskip -2em
D=\nabla_{\bold X}+\sum^{J+Q}_{P=1}\vnabla_{\bold Y_{\!P}}[P]
+\sum^{J+Q}_{P=1}\Bar\vnabla_{\bar\bold Y_{\!P}}[P]+S,
\mytag{20.1}
$$
where $\nabla_{\bold X}$ is the spacial covariant differentiation along
some extended vector field $\bold X$, $\vnabla_{\bold Y_{\!P}}[P]$ is 
the $P$-th canonical vertical multivariate differentiation along some 
extended spin-tensorial field $\bold Y_{\!P}$, $\Bar\vnabla_{\bar\bold
Y_{\!P}}[P]$ is the $P$-th barred canonical vertical multivariate 
differentiation along some other extended spin-tensorial field $\bar
\bold Y_{\!P}$, and $S$ is some degenerate differentiation of the
algebra $\bold S(M)$.
\endproclaim
\demo{Proof} Let $D\in\goth D_{\Bbb C}(M)$. Then its restriction to
the module $S^0_0\bar S^0_0 T^0_0(M)$ is given by some vector field 
$\bold Z$ in $N$. The extended affine connection in $M$ determines 
some horizontal lift of vectors $f$ from $M$ to $N$. According to the 
theorem~\mythetheorem{19.1}, this lift of vectors determines the 
expansion of the tangent space $\Bbb CT_q(N)$ into a direct sum 
\mythetag{19.2} at each point $q\in N$. The vertical subspace $V_q(N)$ 
in \mythetag{19.2} has its own expansion \mythetag{18.24} into a 
direct sum. Combining \mythetag{19.2} and \mythetag{18.24}, we obtain
$$
\hskip -2em
\gathered
T_q(N)=H_q(N)\oplus V_q[1](N)\oplus\ldots\oplus V_q[J+Q](N)\,+\\
\vspace{1ex}
+\,\bar V_q[1](N)\oplus\ldots\oplus\bar V_q[J+Q](N).
\endgathered
\mytag{20.2}
$$
Then the vector field $\bold Z$ is expanded into a sum of vector fields
$$
\hskip -2em
\bold Z=\bold H+\bold V_{\!1}+\ldots+\bold V_{\!J+Q}
+\bar\bold V_{\!1}+\ldots+\bar\bold V_{\!J+Q}.
\mytag{20.3}
$$
uniquely determined by the expansion \mythetag{20.2}. Due to the maps 
\mythetag{18.20} and \mythetag{18.22} and due to the property
\mythetag{18.25} of these maps each vector field $\bold V_{\!P}$ in
\mythetag{20.3} is uniquely associated with some extended tensor field
$\bold Y_{\!P}$ and each vector field $\bar\bold V_{\!P}$ is uniquely
associated with some other extended tensor field $\bar\bold Y_{\!P}$.
Similarly, the vector field $\bold H$ is uniquely associated with some 
extended vector field $\bold X$ by means of the maps \mythetag{19.1}. 
Then we can consider the sum
$$
\hskip -2em
\tilde D=\nabla_{\bold X}+\sum^{J+Q}_{P=1}\vnabla_{\bold Y_{\!P}}[P]
+\sum^{J+Q}_{P=1}\Bar\vnabla_{\bar\bold Y_{\!P}}[P].
\mytag{20.4}
$$
The sum \mythetag{20.4} is a differentiation of $\bold S(M)$ such that
its restriction to $S^0_0\bar S^0_0 T^0_0(M)$ is given by the vector
field \mythetag{20.3}. Hence, $D-\tilde D$ is a differentiation of 
$\bold S(M)$ with identically zero restriction to $S^0_0\bar S^0_0
T^0_0(M)$. This means that $D-\tilde D$ is a degenerate differentiation 
(see the definition~\mythedefinition{14.1}). Denoting this degenerate
differentiation by $S$, from \mythetag{20.4} we derive the required
expansion \mythetag{20.1}. The theorem is proved.\qed\enddemo
    The theorem~\mythetheorem{20.1} is the structural theorem for
differentiations in the algebra of extended spin-tensorial fields 
$\bold S(M)$. It approves our previous efforts in studying the three 
basic types of differentiations which are used in the formula
\mythetag{20.1}. Their application to the description of real physical
fields is the subject for separate publications. The formula 
\mythetag{19.12} is the most important explicit formula for such
applications.\par
\head
21. Dedicatory.
\endhead
    This paper is dedicated to my aunt Abdurahmanova Nailya Muhamedovna
who came through the anxious times of the 20th century and had left us 
in the beginning of this new possibly not less anxious 21st century.


\Refs
\ref\myrefno{1}
\by Sharipov~R.~A.\book Quick introduction to tensor analysis
\publ free on-line textbook in Electronic Archive \myEarXivlink;
see \myhref{http://arXiv.org/abs/math/0403252}{math.HO/0403252}
and \myhref{http://www.geocities.com/r-sharipov/r4-b6.htm}
{r-sharipov/r4-b6.htm} in \myGeoCities
\endref
\ref\myrefno{2}
\by Sharipov~R.~A.\book Course of differential geometry\publ
Bashkir State University\publaddr Ufa\yr 1996\moreref see also
\myhref{http://arXiv.org/abs/math/0412421}{math.HO/0412421}
in Electronic Archive \myEarXivlink\ and 
\myhref{http://www.geocities.com/r-sharipov/r4-b3.htm}
{r-sharipov/r4-b3.htm} in \myhref{http://www.geocities.com}{Geo-}
\myhref{http://www.geocities.com}{Cities}
\endref
\ref\myrefno{3}
\by Sharipov~R.~A.\book Course of linear algebra and multidimensional
geometry\publ Bashkir State University\publaddr Ufa\yr 1996\moreref 
see also \myhref{http://arXiv.org/abs/math/0405323}{math.HO/0405323}
in Electronic Archive \myEarXivlink\ and 
\myhref{http://www.geocities.com/r-sharipov/r4-b2.htm}
{r-sharipov/r4-b2.htm} in \myGeoCities
\endref
\ref\myrefno{4}\by Sharipov~R.~A.\book Classical electrodynamics and
theory of relativity\publ Bashkir State University\publaddr Ufa\yr 1997
\moreref see also
\myhref{http://arXiv.org/abs/physics/0311011}{physics/0311011}
in Electronic Archive \myEarXivlink\ and 
\myhref{http://www.geocities.com/r-sharipov/r4-b5.htm}
{r-sharipov/r4-} \myhref{http://www.geocities.com/r-sharipov/r4-b5.htm}
{b5.htm} in \myGeoCities
\endref
\ref\myrefno{5}\by Sharipov~R.~A.\paper Tensor functions of tensors and 
the concept of extended tensor fields\publ e-print 
\myhref{http://arXiv.org/abs/math-ph/math/0503332/}{math.DG/0503332}
in Electronic Archive \myEarXivlink
\endref
\ref\myrefno{6}\by Sharipov~R.~A.\paper A note on the dynamics and 
thermodynamics of dislocated crystals\publ e-print 
\myhref{http://arXiv.org/abs/cond-mat/0504180}{cond-mat/0504180}
in Electronic Archive \myEarXivlink
\endref
\ref\myrefno{7}\by Berestetsky~V.~B., Lifshits E.~M., 
Pitaevsky~L.~P.\book Quantum Electrodynamics, {\rm Vol\.~\uppercase
\expandafter{\romannumeral 4} of} Theoretical Physics {\rm by 
L.~D.~Landau and E.~M.~Lifshits}\publ Nauka publishers
\publaddr Moscow\yr 1989
\endref
\ref\myrefno{8}\by Penrose~R., Rindler W.\book Spinors and space-time.
{\rm Vol\.~\uppercase\expandafter{\romannumeral 1}.} Two-spinor calculus
and relativistic fields\publaddr Cambridge University Press\yr 1984
\endref
\ref\myrefno{9}\by Dubrovin~B.~A., Novikov~S.~P., Fomenko~A.~T.\book
Modern geometry. {\rm Vol\.~\uppercase
\expandafter{\romannumeral 1}.} Methods and applications\publ Nauka 
publishers\publaddr Moscow\yr 1986
\endref
\ref\myrefno{10}
\by Mishchenko~A.~S.\book Vector bundles and their applications\publ
Nauka publishers\publaddr Moscow\yr 1984
\endref
\ref\myrefno{11}
\by Kobayashi~Sh., Nomizu~K\book Foundations of differential geometry,
{\rm Vol\.~\uppercase\expandafter{\romannumeral 1}}\publ Interscience
Publishers\publaddr New York, London\yr 1963\moreref\publ Nauka publishers
\publaddr Moscow\yr 1981
\endref
\ref\myrefno{12}
\by Sharipov~R.~A.\book Dynamical systems admitting the normal shift
\publ thesis for the degree of Doctor of Sciences in Russia\yr 2000
\moreref see \myhref{http://arXiv.org/abs/math.DG/0002202/}{math.DG/0002202}
in Electronic Archive \myEarXivlink
\endref
\endRefs
\enddocument
\end